\pdfoutput=1
\RequirePackage{ifpdf}
\ifpdf 
\documentclass[pdftex]{sigma}
\else
\documentclass{sigma}
\fi

\usepackage{mathtools}
\usepackage{bbm}
\usepackage{enumerate}
\usepackage{tikz}
\usetikzlibrary{decorations.pathreplacing}
\usepackage{array}
\usepackage{soul}
\usepackage{verbatim}
\usepackage{siunitx}
\usetikzlibrary{calc} \tikzset{>=latex}
\usetikzlibrary{calc,decorations.markings}
\usetikzlibrary{positioning,shapes,decorations.text}

\theoremstyle{plain}
\newtheorem{thm}{Theorem}[section]
\newtheorem{cor}[thm]{Corollary}
\newtheorem{Lemma}[thm]{Lemma}

\newtheorem{prop}[thm]{Proposition}

\newtheorem{assumption}[thm]{Assumption}
\theoremstyle{definition}
\newtheorem{defn}[thm]{Definition}
\newtheorem{rem}[thm]{Remark}
\newtheorem{eg}[thm]{Example}
\newtheorem{fact}[thm]{Fact}
\newtheorem{observe}[thm]{Observation}

\newcommand*{\medcup}{\mathbin{\scalebox{1.5}{\ensuremath{\cup}}}}

\newcommand{\rpm}{\sbox0{$1$}\sbox2{$\scriptstyle\pm$}
 \raise\dimexpr(\ht0-\ht2)/2\relax\box2 }

\tikzstyle{nd} = [anchor=base, inner sep=0pt]
\tikzstyle{ndpic} = [remember picture, baseline, every node/.style={nd}]

\DeclareMathOperator{\Pf}{Pf}
\def\beq{\begin{equation}}
\def\eeq{\end{equation}}
\def\ba{\begin{enumerate}[(a)]}
\def\bei{\begin{enumerate}[(i)]}
\def\be{\begin{enumerate}[(1)]}
\def\ee{\end{enumerate}}
\def\bi{\begin{itemize}}
\def\ei{\end{itemize}}
\def\beg{\begin{eg}}
\def\eeg{\end{eg}}
\def\bd{\begin{defn}}
\def\ed{\end{defn}}
\def\bt{\begin{thm}}
\def\et{\end{thm}}
\def\bl{\begin{Lemma}}
\def\el{\end{Lemma}}
\def\bfac{\begin{fact}}
\def\efac{\end{fact}}

\def\bc{\begin{cor}}
\def\ec{\end{cor}}
\def\bp{\begin{prop}}
\def\ep{\end{prop}}
\def\bo{\begin{observe}}
\def\eo{\end{observe}}
\def\bas{\begin{assumption}}
\def\eas{\end{assumption}}
\def\RR{\mathbb{R}}
\def\CC{\mathbb{C}}

\def\ZZ{\mathbb{Z}}
\def\NN{\mathbb{N}}

\def\beg{\begin{eg}}
\def\eeg{\end{eg}}

\def\tt{\mathrm{T}^{(2)}_{u,v}}

\makeindex

\numberwithin{equation}{section}

\begin{document}

\allowdisplaybreaks

\newcommand{\arXivNumber}{1705.05859}

\renewcommand{\PaperNumber}{092}

\FirstPageHeading

\ShortArticleName{Pfaffian Schur Correlation via Difference Operators}

\ArticleName{Correlation Functions of the Pfaffian Schur Process\\ Using Macdonald Difference Operators}

\Author{Promit GHOSAL}

\AuthorNameForHeading{P.~Ghosal}

\Address{Department of Statistics, Columbia University,\\ 1255 Amsterdam Avenue, New York, NY 10027, USA}
\Email{\href{mailto:pg2475@columbia.edu}{pg2475@columbia.edu}}

\ArticleDates{Received September 14, 2018, in final form November 19, 2019; Published online November 26, 2019}

\Abstract{We study the correlation functions of the Pfaffian Schur process. Borodin and Rains~[\textit{J.~Stat. Phys.} \textbf{121} (2005), 291--317] introduced the Pfaffian Schur process and de\-ri\-ved its correlation functions using a~Pfaffian analogue of the Eynard--Mehta theorem. We present here an alternative derivation of the correlation functions using Macdonald difference operators.}


\Classification{60C05; O5E05}

\vspace{-4mm}

\section{Introduction}
\vspace{-1mm}

In the last decade, we have seen a great success of the applications of the Schur process in probability and other related fields. Okounkov~\cite{O01} introduced the Schur measure. Later, Okounkov and Reshetikhin~\cite{OR03} generalized it to the Schur process. Schur process have been proved useful in number of different occasions, including, harmonic analysis of infinite symmetric group~\cite{BO03}, Fredholm determinant formula for Toeplitz determinants~\cite{BO00}, relative Gromov--Witten theory of the Riemann sphere $\mathbb{P}^1$~\cite{O03}, random domino tilings of the Aztec diamond~\cite{BBBC14, Bor11}, poly\-nuclear growth processes~\cite{KurtJ03}, anisotropic random growth models in $2+1$ dimensions~\cite{BorFer14Aniso} etc. See~\cite{BorGorLec} and the reference therein for other applications of the Schur process. Later, Borodin and Corwin~\cite{BorodinCorwinMac} further extended the scope by introducing the Macdonald processes where the underlying measure is defined using the Macdoanld $(q,t)$-polynomials~\cite{Mac15}. Macdonald processes have been applied in computing the asymptotics of the one-point marginals of O'Connell--Yor semi-discrete directed polymers \cite[Chapter~5]{BorodinCorwinMac}, \cite{BCF14}; log-gamma discrete directed polymer \cite[Chapter~5]{BorodinCorwinMac}, \cite{BCR13}; Kardar--Parisi--Zhang and stochastic heat equation~\cite{BCF14}; $q$-TASEP \cite{BorodinCorwinMac,BC15,BCS14}; $q$-PushTASEP \cite{BorPet16, CorPet15}; in showing the Gaussian free field fluctuations in $\beta$ Jacobi corner process~\cite{BorGor15} and in constructing a multilevel version of the Dyson Brownian motion~\cite{GS15}. We refer to~\cite{CorwinICM} for the discussion on many other aspects of the Macdonald processes.

 Schur process is a measure on a sequence of partitions defined using Schur functions. The correlation functions of the Schur process can be expressed as determinants of some positive de\-finite matrices~\cite{OR03}. Baik and Rains~\cite{BaiRa01Du} generalized the Schur measure for studying the longest increasing subsequences of symmetrized random permutations. Subsequently, \cite{R00} (see also~\cite{FNR06}) computed the correlation functions of a Pfaffian point process associated to symmetrized increasing subsequence problem and the geometric weight half-space last passage percolations (LPP). Later, generalizing the works of~\cite{R00, SmIm04}, \cite{BR05} introduced the Pfaffian Schur process. Unlike the Schur process, the correlation functions of the Pfaffian Schur process can be expressed as the Pfaffians of minors of a single infinite skew symmetric matrix. That matrix is referred to the \emph{kernel} of the correlation functions. In a recent work, \cite{BBCS16a} used the Pfaffian Schur process to determine the exact distribution of geometric LPP in the half space. Furthermore, they also extended the results of \cite{FNR06, SmIm04} to the exponential weight LPP and their results work for any choice of boundary parameters. They further broadened the scope by relating the exponential LPP in the half space with facilitated TASEP in~\cite{BBCS16b}. Recently, \cite{BBNV17} considered a further generalization of the Schur process which they named as free boundary Schur process. This new stochastic process comes with two free boundaries. From the free boundary Schur process, \cite{BBNV17} recovered the original Schur process by fixing both the boundaries and obtained the Pfaffian Schur process when one of the boundaries is kept fixed. We also refer to a recent work of~\cite{BBCW17} where they inroduced the half-space Macdonald processes to study the stochastic six vertex model in half-quadrant and the asymmetric simple exclusion process, the Kardar--Parisi--Zhang equation in half-space. We came across some interesting works which have few overlaps with ours since this paper had been posted. Here, we mention two of them. In~\cite{GuiBorCor17}, the authors studied the half-space Macdonald processes which is a generalization of the Pfaffian Schur process. They computed the moments and the Laplace transform formula for general half-space Macdonald observables using Mac\-donald difference operators. More recently, \cite{Betea18} found the correlation functions of the symplectic and the orthogonal Schur measures using two different approaches. One of his approach is based on Macdonald difference operators and uses some of the computations of our present work.

Okounkov and Reshetikhin \cite{O01} derived the correlation functions of the Schur process using the infinite wedge formalism. Borodin and Rains~\cite{BR05} provided another derivation of the correlation functions of the Schur process using the Eynard--Mehta theorem. Furthermore, they proposed a Pfaffian analogue of the Eynard--Mehta theorem and used that to find the correlation functions for the Pfaffian Schur process. In this paper, we present a new derivation of the correlation functions of the Pfaffian Schur process.

In their work, \cite{BR05} mainly used some tools from linear algebra. We will use Macdonald difference operators which emerges in the theory of symmetric functions. In the past few years, Macdonald difference operators has been proved immensely useful in the study of stochastic processes over the sequence of partitions. In Borodin and Corwin's work \cite{BorodinCorwinMac} (see also~\cite{BCGS16}), it plays a crucial role in deriving the Fredholm determinant formula of the $q$-exponential moments for the $q$-Whittaker process. Later, building on \cite[Remark~2.2.15]{BorodinCorwinMac}, \cite{A14} gave another derivation of the correlation functions of the Schur process using Macdonald difference operators. Apart from its application in probability, Macdonald difference operators have been widely used in algebraic combinatorics especially for deriving various identities of the symmetric functions namely, Kirillov--Noumi--Warnaar identity \cite{Warnaar08}; Lassalle--Schlosser identity \cite{Lm06,LmSm06} etc. See~\cite{BW16} and the reference therein for other useful identities derived using Macdonald difference operators and its connection to other fields.

\looseness=-1 One of the key properties of Macdonald difference operators is that they are diagonalized by the Macdonald polynomials (see~\cite{Mac15}) on the space of symmetric functions. One of the other important properties is that the action of the difference operators on the partition function of Macdonald processes has a nice representation in terms of the contour integrals (see~\cite{BorodinCorwinMac}). These two facts yield integral forms (amenable to further analysis) for the expected value of some of the interesting observables of the Macdonald processes (see~\cite{BCGS16}). From these formulas, in order to extract the form of the correlation functions, one needs to perform some algebraic manipulations inside the integrals \cite[Section~2.3]{A14}. To find the correlation functions of the Pfaffian Schur process, we take the same route. The main contributions of the present work are to derive the contour integral formulas for the action of the difference operators on the Pfaffian Schur partition function and then, to use those for finding the correlation kernel for the Pfaffian Schur process.

\section{Models and main result}

 In this section, we first define all the necessary terms and notations which we follow in the rest of the paper. Subsequently, we give a formal definition of the Pfaffian Schur process. Then, we proceed to state the main result of this paper.

\subsection{Partitions and symmetric functions}\label{PlanePartAndSymFunc}
 A \emph{partition} $\lambda$ is defined as a sequence of ordered non-negative integers $\lambda_1\geq \lambda_2\geq \cdots $ such that its \emph{weight} $|\lambda|:=\sum_{i=1} \lambda_i $ is finite. We call~$\lambda$ an \emph{even} partition when all of its components ($\lambda_i$'s) are even integers. Any partition $\lambda$ can also be represented graphically as a~\emph{Young diagram} with~$\lambda_1$ left justified box in the top row, $\lambda_2$ in the second row and so on. The set of all partitions (or Young diagrams) is denoted by~$\mathbb{Y}$. The Cartesian product of $m$ sets $\mathbb{Y}\times \cdots \times \mathbb{Y}$ will be denoted as~$\mathbb{Y}^m$. The \emph{transpose} of a Young diagram is denoted by $\lambda^\prime$ and defined by $\lambda^\prime_i:=|\{j,\lambda_j\geq i\}|$. The \emph{length} $l(\lambda)$ of a partition $\lambda$ is the number of nonzero entries in the sequence. We denote the set of all partitions whose weights are $n$ by $\mathbb{Y}_n$. Thus, we have $\mathbb{Y}=\cup_{n=0}^\infty \mathbb{Y}_n$. For any two partitions~$\lambda$ and~$\mu $, if we have $\lambda \supset \mu$ (as the set of boxes), we call $\lambda-\mu$ as a \emph{skew Young diagram}. Often $\lambda-\mu$ is denoted as $\lambda/\mu$.

Any polynomial of the variables $x_1,\dots ,x_n$ is called \emph{symmetric} if it is invariant under the action of the permutation group $\mathfrak{S}_n$. The set of all symmetric polynomials of the variables $x_1,\dots ,x_n$ is a subalgebra inside the algebra $\mathbb{C}[x_1,\dots ,x_n]$ of polynomials. For any $n$-tuple of integers $\alpha=(\alpha_1,\dots ,\alpha_n)\in \ZZ^n_{\geq 0}$, we can associate \emph{monomial symmetric polynomial} as
 \begin{gather*}
 m_\alpha(x_1,\dots ,x_n)= \sum_{\sigma \in \mathfrak{S}(n)}\prod_{i} x^{\alpha_{\sigma(i)}}_i.
\end{gather*}
For any Young diagram $\lambda$ with length less than $n$, one can define $m_{\lambda}$ in the same way as above by considering $\lambda$ as a $n$-tuple. The set of all $m_\lambda$ where $|\lambda|=k$ spans a subalgebra ${\rm Sym}^k_n$ inside the algebra $\mathbb{C}[x_1,\dots ,x_n]$. Furthermore, there exists an algebra homomorphism between ${\rm Sym}^k_{n}$ and ${\rm Sym}^k_{m}$ for any $m>n\in \NN$ by mapping $f(x_1,\dots ,x_n,\dots,x_m)$ to $f(x_1,\dots ,x_n, 0,\dots ,0)$. Thus, one can define the direct limit ${\rm Sym}^k$ of the sequence of subalgebras $\big\{{\rm Sym}^k_n\big\}_{n\in \NN}$ for each fixed $k$. The direct sum of all such subalgebras ${\rm Sym}:=\oplus_k {\rm Sym}^k$ is referred as the \emph{algebra of symmetric functions} with countably many variables. In fact, it can be noted that ${\rm Sym}$ is a $\ZZ_{\geq 0}$-graded algebra. Some typical examples of the symmetric functions includes, $(a)$ \emph{the elementary symmetric functions} $e_r:=m_{1^r}$; $(b)$ \emph{the complete homogeneous symmetric functions} $h_r:=\sum_{|\lambda|=r}m_{\lambda}$; $(c)$ \emph{the power sum symmetric functions} $p_r(x_1,\dots ,x_n):=\sum_{i}x^r_i$ and $p_\lambda:= \prod_{\lambda_i}p_{\lambda_i}$ etc. Let us also point out that $\{p_\lambda\}_{\lambda\in \mathbb{Y}}$ form a basis of the space ${\rm Sym}$. \emph{Schur} and \emph{skew Schur} functions are special kind of symmetric functions which are defined as
\begin{gather*}
s_\lambda:= \det (h_{\lambda_i-i+j} ), \qquad \text{and} \qquad s_{\lambda/\mu}:= \det (h_{\lambda_i-\mu_j-i+j} ).
\end{gather*}

 The \emph{topological completion} $\overline{{\rm Sym}}$ of the space ${\rm Sym}$ is defined as the set of all formal power series
 \begin{gather*}
 a=\sum_{k=0}^\infty a_k, \qquad a_k\in {\rm Sym}^k.
 \end{gather*}
For any $a\in \overline{{\rm Sym}}$, its \emph{lower degree} $\mbox{ldeg}(a)$ is defined as the maximal $K$ such that $a_k=0$ for all $k<K$. Furthermore, one can also consider a \emph{graded topology} in $\overline{{\rm Sym}}$ by defining $b:=\lim\limits_{n\to \infty}b_n $ when $\mbox{ldeg}(b-b_n)$ converges to $\infty$ as $n$ tends to $\infty$.

 We denote the algebra of symmetric functions of any set of variables $X=(x_1,x_2,\dots )$ as ${\rm Sym}(X)$. For any $f\in {\rm Sym}$, one can represent $f(X,Y)\in {\rm Sym}(X,Y)$ as a sum of products of symmetric polynomials of $X$ and symmetric polynomials of $Y$ where $X$ and $Y$ are two sets of variables and $(X,Y)$ denotes their union ($X\cup Y$). This defines a comultiplication from the space ${\rm Sym}$ to ${\rm Sym}\otimes {\rm Sym}$ which turns ${\rm Sym}$ into a bi-algebra. Moreover, the topological completion of the space ${\rm Sym}(X)\otimes {\rm Sym}(Y)$ is given by $\overline{{\rm Sym}(X)}\otimes \overline{{\rm Sym}(Y)}$.

Any algebra homomorphism from the space ${\rm Sym}$ to $\mathbb{C}$ is referred as \emph{specialization}. For instance, the evaluation of any $f\in {\rm Sym}$ at any set $X$ comprised of finitely many elements from~$\CC$ is an example of specialization. We often call these as finite specializations. For any two specializations $\rho_1$ and $\rho_2$, their disjoint union $\rho_1\cup \rho_2$ is defined by
 \begin{gather*}
 p_k(\rho_1\cup \rho_2) := p_k(\rho_1)+p_k(\rho_2),
 \end{gather*}
 which implies it is also a specialization. We will often denote $\rho_1\cup\rho_2$ by $(\rho_1,\rho_2)$.

\bd
Consider any specialization $\rho$. We call $\rho$ \emph{Schur non-negative} if $s_{\lambda/\mu}(\rho)$ is positive for all $\mu\subset \lambda$ where $\mu, \lambda\in \mathbb{Y}$.
\ed

 Schur functions are pairwise orthogonal with respect to a bilinear form given as
\begin{gather*}
 \langle p_\lambda(X),p_\mu(X) \rangle^{X} := \mathbbm{1}_{\lambda=\mu} \prod_{i\geq 1} i^{m_i(\lambda)} (m_i(\lambda))!,
\end{gather*}
where $\lambda = 1^{m_1} 2^{m_2}\cdots $. This defines an inner product in the space ${\rm Sym}(X)$. In particular, we have $\langle s_\lambda(X), s_\mu(X) \rangle^X = \mathbbm{1}_{\lambda=\mu}$ and for any two set of variables $X$ and $Z$, we get $\langle s_\lambda(Z,X), s_\mu (X)\rangle^X = s_{\lambda/\mu }(Z)$ (see \cite[Corollary~2.1.2]{A14}) where the inner product is considered in the space ${\rm Sym}(X)$. One can also restrict the inner product over the set of partitions with bounded weights. For instance, one can define
\begin{gather*}
 \langle p_\lambda(X),p_\mu(X) \rangle^{X}_{u} := \mathbbm{1}_{\lambda=\mu} \mathbbm{1}_{|\mu|\leq u}\prod_{i\geq 1} i^{m_i(\lambda)} (m_i(\lambda))!
\end{gather*}
for any $u\in \NN$. In terms of the finite inner product, the orthogonality of the Schur functions translates to $\langle s_\lambda(X), s_\mu(X) \rangle^X_u = \mathbbm{1}_{\lambda=\mu}\mathbbm{1}(|\mu|\leq u)$. Furthermore, one can also write
\begin{align}\label{eq:BilinearFormSkewSpecial}
\langle s_\lambda(Z,X), s_\mu (X)\rangle^X_u = s_{\lambda/\mu }(Z)\mathbbm{1}(|\mu|\leq u).
\end{align}
For more detailed expositions, see \cite[Chapter~2]{BC15}, \cite[Section~2]{BCGS16}, \cite[Chapter~1, Section~5]{Mac15}. Let us also point out that we will often use the notations $\vee$ and $\wedge$ to denote the maximum and minimum between any two real numbers.

\subsection{Schur and Pfaffian Schur process}\label{SchurPfaffSchur}
 For any $m\in \NN$, the \emph{Pfaffian Schur process} is a probability measure on a sequence of partitions
 \begin{gather}\label{eq:PartitionSeq}
 \varnothing \subset \lambda^{(1)}\supset \mu^{(1)} \subset \lambda^{(2)} \supset \cdots \subset \mu^{(n-1)} \subset \lambda^{(m)} \supset \varnothing
\end{gather}
given by the following product form
\begin{gather}
\mathcal{P}_{psp}(\bar{\lambda},\bar{\mu}; \rho^{+}, \rho^{-}):= c \tau_{\lambda^{(1)}}(\rho_0^{-})s_{\lambda^{(1)}/\mu^{(1)}}(\rho_1^{+})s_{\lambda^{(2)}/\mu^{(1)}}(\rho_1^{-})\cdots \nonumber\\
\hphantom{\mathcal{P}_{psp}(\bar{\lambda},\bar{\mu}; \rho^{+}, \rho^{-}):=}{}\times s_{\lambda^{(n)}/\mu^{(m-1)}}(\rho^{-}_{m-1})s_{\lambda^{(m)}}(\rho^{+}_{m}),\label{eq:PfaffianSchurMeasure}
\end{gather}
where $s_\lambda$, $s_{\lambda/\mu}$ are the Schur and skew Schur functions, $\rho^{\pm}_i$ are Schur non-negative specializations and $\tau$ is defined by $\tau_\lambda=\sum_{\mu^{\prime} \, \text{even}} s_{\lambda/\mu}$. Here, we have used the following shorthand notations $\rho^{\pm}:=\cup_{i}\rho^{\pm}_i$, $\bar{\lambda}=\big(\lambda^{(1)}, \dots ,\lambda^{(m)}\big)$ and $\bar{\mu}=\big(\mu^{(1)}, \dots ,\mu^{(m-1)}\big)$. The fact that all those specializations are non-negative implies that all the weights are non-negative. Furthermore, in order to define a measure using the weights in~\eqref{eq:PfaffianSchurMeasure}, we assume that the series $\sum_{(\bar{\lambda},\bar{\mu})}\mathcal{P}_{psp}(\bar{\lambda},\bar{\mu}; \rho^{+}, \rho^{-})$ is absolutely convergent. In that case, the value of the constant~$c$ in~\eqref{eq:PfaffianSchurMeasure} is equal to the inverse of the \emph{partition function},
\begin{gather*}
\begin{split}&
Z(\rho^{+}; \rho^{-}):= \sum_{(\bar{\lambda}, \bar{\mu})\in \mathbb{Y}^m\times \mathbb{Y}^{m-1}} \tau_{\lambda^{(1)}}(\rho_0^{-})s_{\lambda^{(1)}/\mu^{(1)}}(\rho_1^{+})s_{\lambda^{(2)}/\mu^{(1)}}(\rho_1^{-}) \cdots \\
& \hphantom{Z(\rho^{+}; \rho^{-}):= \sum_{(\bar{\lambda}, \bar{\mu})\in \mathbb{Y}^m\times \mathbb{Y}^{m-1}} }{} \times s_{\lambda^{(m)}/\mu^{(m-1)}}(\rho^{-}_{m-1})s_{\lambda^{(m)}}(\rho^{+}_{m}).\end{split}
\end{gather*}
The probability measure of any set $\bar{A}$ (containing sequences of partitions) under the Pfaffian Schur process is given by
\begin{gather*}
\mathbb{P}_{psp}(\bar{A}):= \sum_{(\bar{\lambda},{\bar{\mu}})\in \mathbb{Y}^{m}\times \mathbb{Y}^{m-1}} \mathcal{P}_{psp}(\bar{\lambda},\bar{\mu}; \rho^{+}, \rho^{-})\mathbbm{1}(\bar{\lambda}\in \bar{A}).
\end{gather*}
We denote the expectation of any function $f$ of sequence of partitions by $\mathbb{E}_{psp}(f)$ where \[\mathbb{E}_{psp}(f):= \sum_{(\bar{\lambda},{\bar{\mu}})\in \mathbb{Y}^{m}\times \mathbb{Y}^{m-1}} \mathcal{P}_{psp}(\bar{\lambda},\bar{\mu}; \rho^{+}, \rho^{-}) f(\bar{\lambda}).\]
To see the contrast between the Pfaffian Schur process and the Schur process, let us also look at the definition of the latter. Like in \emph{Pfaffian Schur process}, the \emph{Schur process} on the sequence of partitions considered in \eqref{eq:PartitionSeq} is proportional to
\begin{gather*}
s_{\lambda^{(1)}}(\rho_0^{-})s_{\lambda^{(1)}/\mu^{(1)}}(\rho_1^{+})s_{\lambda^{(2)}/\mu^{(1)}}(\rho_1^{-})\cdots s_{\lambda^{(n)}/\mu^{(m-1)}}(\rho^{-}_{m-1})s_{\lambda^{(m)}}(\rho^{+}_{m})
\end{gather*}
and the associated \emph{partition function} will be denoted as $F(\rho^{+};\rho^{-})$. If $\rho_1$ and $\rho_2$ are any two non-negative specialization, then the Schur process defined over the set of all single partitions is referred to the \emph{Schur measure} which is given as
\begin{gather*}
\mathcal{P}_{sm}(\rho_1,\rho_2)(\lambda):= \frac{s_{\lambda}(\rho_1)s_{\lambda}(\rho_2)}{F(\rho_1;\rho_2)}.
\end{gather*}
Likewise, the Pfaffian Schur measure for a single plane partition $\lambda$ is defined as
\begin{gather}\label{eq:DefPfaffianSchurMeasure}
\mathcal{P}_{psm}(\rho_1,\rho_2)(\lambda):= \frac{\tau_{\lambda}(\rho_2)s_{\lambda}(\rho_1)}{Z(\rho_1;\rho_2)}.
\end{gather}
We denote the probability measure of any set $A$ and the expectation of any function $f$ of partitions (under the Pfaffian Schur measure) by $\mathbb{P}_{psm}(A)$ and $\mathbb{E}_{psm}(f)$ respectively.
For any two sets $X$ and $Y$ where $X=(x_1,x_2, \dots )$ and $Y=(y_1,y_2,\dots)$, we define the following two functions
\begin{gather}\label{eq:HallLittleSchur}
H_0(X):=\prod_{i<j}\frac{1}{1-x_ix_j} \qquad \text{and} \qquad H(X;Y):= \prod_{i, j}\frac{1}{1-x_iy_j}.
\end{gather}
Using the expansion $(1-x)^{-1}= \sum_{k=0}^\infty x^k$, one can note that $H_0(X)\in \overline{{\rm Sym}(X)}$ and $H(X,Y)\in \overline{{\rm Sym}(X)}\otimes \overline{{\rm Sym}(Y)}$.
Furthermore, using the identity $(1-x)^{-1}=\exp\big(\sum^\infty_{k=1}k^{-1}x^k\big)$, we write
\begin{gather*}
H_0(X)=\exp\left(\sum_{k=1}^\infty \frac{p^2_k(X)- p_{2k}(X)}{2k}\right) \qquad \text{and} \qquad H(X;Y)=\exp\left(\sum_{k=1}^\infty \frac{p_k(X)p_k(Y)}{k}\right).
\end{gather*}
In general, for any two Schur non-negative specializations $\rho_1$ and $\rho_2$, we define
\begin{gather}\label{eq:SpecializedH0H1}
H_0(\rho_1):= \exp\left(\sum_{k=1}^\infty \frac{p^2_k(\rho_1)- p_{2k}(\rho_1)}{2k}\right),\qquad H(\rho_1;\rho_2):= \exp\left(\sum_{k=1}^\infty \frac{p_k(\rho_1)p_{k}(\rho_2)}{k}\right),
\end{gather}
when the respective series $\sum^\infty_{k=1}(2k)^{-1}(p^2_k(\rho_1)- p_{2k}(\rho_1))$ and $\sum_{k=1}^\infty k^{-1}p_k(\rho_1)p_{k}(\rho_2)$ converge absolutely.

\bp[\cite{BR05}]\label{PfaffianPartitionFunctionIdentity}
Fix any two sequence of Schur non-negative specializations $\rho^{+}=(\rho^{+}_1, \dots , \rho^{+}_m)$ and $\rho^{-}=(\rho^{-}_0, \dots , \rho^{-}_{m-1})$ such that $\sum_{(\bar{\lambda},\bar{\mu})}\mathcal{P}_{psp}(\bar{\lambda},\bar{\mu}; \rho^{+}, \rho^{-})$ converges absolutely.
Then, one can write the partition functions $Z(\rho^{+};\rho^{-})$ and $F(\rho^{+};\rho^{-})$ of the Pfaffian Schur process and the Schur process as
\begin{gather}
Z(\rho^{+}; \rho^{-}) = \prod_{i=1}^m H_0(\rho^{+}_i) \prod_{0\leq i<j\leq m} H\big(\rho^{-}_i; \rho^{+}_j\big), \label{eq:PfaffianSchurPartitionIdentity} \\
F(\rho^{+}; \rho^{-}) = \prod_{0\leq i<j\leq m} H\big(\rho^{-}_i; \rho^{+}_j\big). \label{eq:SchurPartitionIdentity}
\end{gather}
\ep

\subsection{Main result}\label{MainRes}
 We first define the correlation functions of the Pfaffian Schur process. To any sequence of partitions defined in~\eqref{eq:PartitionSeq}, one can associate a point configuration in $[1,m]\times \ZZ$ as
\begin{gather}\label{eq:PointConfiguration}
\big\{\big(1,\lambda^{(1)}_i-i\big)\big\}_{i\in \NN} \cup \big\{\big(2,\lambda^{(2)}_i-i\big)\big\}_{i\in \NN} \cup \cdots \cup \big\{\big(m,\lambda^{(m)}_i-i\big)\big\}_{i\in \NN},
\end{gather}
where $[i,j]:=\{i,i+1,\dots ,j\}$ for any two numbers $i<j\in \NN$ and $\big(\lambda^{(1)},\dots ,\lambda^{(m)}\big)$ is a set of $m$-tuple of partitions. Moreover, one can note that the point configuration in the display~\eqref{eq:PointConfiguration} uniquely determines the sequence $\big(\lambda^{(1)}, \dots ,\lambda^{(m)}\big)$.

\bd[correlation functions] Fix $m\in \NN$. For any set of the form $T:= \cup_{i=1}^m T_i$ where $T_i$'s are given as $T_{i}:= \{(i,t_{i,1}), \dots ,(i,t_{i,d_i})\}$, the correlation functions with respect to the Pfaffian Schur process is defined as
\begin{gather}\label{eq:DefCorrelationFunction}
\rho_{psp}(T):=\mathbb{P}_{psp}\big(t_{i,j}\in \big\{\lambda^{(i)}_\kappa-\kappa\big\}_{\kappa\in \NN}\,|\,1\leq i\leq m,\, 1\leq j\leq d_i\big).
\end{gather}
\ed
 The correlation functions of the Pfaffian Schur process are expressed in terms the Pfaffian of a skew symmetric matrix which is the prime object of the study in this paper. Recall that the Pfaffian of any skew symmetric matrix $A_{2d\times 2d}$ is defined as
 \begin{gather*}
 \Pf(A)=\frac{1}{2^d d!}\sum_{\sigma\in \mathfrak{S}(2d)} (-1)^{\sigma} \prod_{i=1}^d A_{\sigma(2i-1), \sigma(2i)},
 \end{gather*}
 where $\mathfrak{S}(2d)$ is the group of all permutations of the numbers in the set $[1,2d]$.
 In what follows, we present our main result. We consider the Pfaffian Schur process defined in \eqref{eq:PfaffianSchurMeasure} with two sets of finite Schur non-negative specializations $\rho^{+}:=\{\rho^{+}_1,\dots , \rho^{+}_m\}$ and $\rho^{-}:=\{\rho^{-}_0,\dots, \rho^{-}_{m-1}\}$. Define $\rho^{+}_{\mathcal{S}}:= \cup_{i\in \mathcal{S}} \rho^{+}_{i}$ and similarly, $\rho^{-}_{\mathcal{S}}$ for any set $\mathcal{S}\subseteq [1,m]$. We assume that the series $\sum_{(\bar{\lambda},\bar{\mu})}\mathcal{P}_{psp}(\bar{\lambda},\bar{\mu}; \rho^{+}, \rho^{-})$ is absolutely convergent.

 \bt\label{MainTheorem}
 Fix $m\in \NN$, $T_i = \{(i,t_{i,1}),\dots ,(i,t_{i,d_i})\}$ for $i=1,2, \dots ,m$ and define $T:=\cup_{i=1}^m T_i$, $d:=\sum_{i=1}^m d_i$. Then, the correlation functions $\rho_{psp}(T)$ is given by $\Pf\big(\mathcal{K}^{(T)}\big)$, i.e., Pfaffian of the $2d\times 2d$ skew symmetric matrix $\mathcal{K}^{(T)}$ $($see Remark~{\rm \ref{rem:Kmatrix}} for detailed description of $\mathcal{K}^{(T)})$. The coefficients of $\mathcal{K}^{(T)}$ are given explicitly by the contour integrals
\begin{gather}
\mathcal{K}^{(T)}_{1,1}(i,u;j,v) = \frac{1}{(2\pi{\rm i})^{2} }\oint_{\mathfrak{C}_1}\oint_{\mathfrak{C}_2} \frac{z-w}{\big(z^2-1\big)\big(w^2-1\big)(zw-1)}\nonumber\\
\hphantom{\mathcal{K}^{(T)}_{1,1}(i,u;j,v) =}{} \times\frac{H\big(\rho^{+}_{[i,m]}; z\big)H\big(\rho^{+}_{[j,m]}; w\big)}{H\big(\rho^{+}_{[1,m]}\cup \rho^{-}_{[0,i)}; z^{-1}\big)H\big(\rho^{+}_{[1,T]}\cup \rho^{-}_{[0,j)}; w^{-1}\big)} \frac{{\rm d}z {\rm d}w}{z^{t_{i,u}}w^{t_{j,v}}},\label{eq:PfaffianMatDes1stTerm}
\end{gather}
 where both $\mathfrak{C}_1$ and $\mathfrak{C}_2$ are circles with radius greater than $1$,
\begin{gather}
\mathcal{K}^{(T)}_{1,2}(i,u;j,v) = \frac{1}{(2\pi{\rm i})^{2} }\oint_{\mathfrak{C}_3}\oint_{\mathfrak{C}_4} \frac{z-w}{w\big(z^2-1\big)(zw-1)}\nonumber\\
\hphantom{\mathcal{K}^{(T)}_{1,2}(i,u;j,v) =}{} \times\frac{H\big(\rho^{+}_{[1,m]}\cup \rho^{-}_{[0,i)}; w\big)H\big(\rho^{+}_{[j,m]}; z\big)}{H\big( \rho^{+}_{[i,m]}; w^{-1}\big)H\big(\rho^{+}_{[1,m]}\cup \rho^{-}_{[0,j)}; z^{-1}\big)}\frac{{\rm d}z {\rm d}w}{z^{t_{i,u}}w^{t_{j,v}}},\label{eq:PfaffianMatDes2ndTerm}
\end{gather}
where the closed contours $\mathfrak{C}_3$ and $\mathfrak{C}_4$ satisfy $(a)$ $|z|>1$ and $|zw|<1$ for all $z\in \mathfrak{C}_3$ and $w\in \mathfrak{C}_4$ whenever $i\leq j$, and $(b)$ $|z|>1$, $|zw|>1$ for all $z\in \mathfrak{C}_3$ and $w\in \mathfrak{C}_4$ whenever $i> j$,
\begin{gather*}
\mathcal{K}^{(T)}_{2,1}(i,u;j,v) = - \mathcal{K}^{(T)}_{1,2}(i,v;j,u),
\end{gather*}
and
\begin{gather}
\mathcal{K}^{(T)}_{2,2}(i,u;j,v) = \frac{1}{(2\pi{\rm i})^{2} } \oint_{\mathfrak{C}_5}\oint_{\mathfrak{C}_6} \frac{z-w}{zw(zw-1)}\nonumber\\
\hphantom{\mathcal{K}^{(T)}_{2,2}(i,u;j,v) =}{}
\times\frac{H\big(\rho^{+}_{[1,m]}\cup \rho^{-}_{[0,i)}; z\big)H\big(\rho^{+}_{[1,m]}\cup \rho^{-}_{[0,i)}; w\big)}{H\big( \rho^{+}_{[i,m]}; z^{-1}\big)H\big(\rho^{+}_{[j,m]}; w^{-1}\big)}\frac{{\rm d}z {\rm d}w}{z^{t_{i,u}}w^{t_{j,v}}},\label{eq:PfaffianMatDes3rdTerm}
\end{gather}
where $\mathfrak{C}_5$ and $\mathfrak{C}_6$ are circles with radius less than $1$. In addition, all the contours considered in this theorem are oriented anticlockwise.
\et

\begin{rem}\label{rem:Kmatrix}The matrix $\mathcal{K}^{(T)}$ of Theorem~\ref{MainTheorem} can be partitioned into $m^2$ submatrices. The size of the $(i,j)$-submatrix of $\mathcal{K}^{(T)}$ is $2d_i\times 2d_j$ where $i,j\in [1,m]$. Moreover, each of the submatrices are composed of $2\times 2$ block matrices. For instance, there are $d_i\cdot d_j$ blocks in the $(i,j)$-submatrix. We denote the $(u,v)$-block matrix of the $(i,j)$-submatrix by $\mathcal{K}^{(T)}(i,u;j,v)$.
The matrix $\mathcal{K}^{(T)}$ is completely specified by knowing the values of all $\mathcal{K}^{(T)}(i,u;j,v)$ where
\begin{gather*}
\mathcal{K}^{(T)}(i,u;j,v):= \begin{pmatrix}
\mathcal{K}^{(T)}_{1,1}(i,u;j,v) & \mathcal{K}^{(T)}_{1,2}(i,u;j,v)\\
\mathcal{K}^{(T)}_{2,1}(i,u;j,v) & \mathcal{K}^{(T)}_{2,2}(i,u;j,v)
\end{pmatrix}
\end{gather*}
for $1\leq i,j\leq m$ and for any given $i,j\in [1,m]$, $u$ and $v$ vary in between $[1,d_i]$ and $[1,d_j]$ respectively.
\end{rem}

\begin{rem}\label{IMpRem}It should be noted that the expression of $\mathcal{K}^{(T)}_{22}$ in~\eqref{eq:PfaffianMatDes3rdTerm} is in disagreement with the expression given in \cite[Theorem~3.3]{BR05}. Inside the integral, \cite{BR05} has $(1-zw)$ instead of $(zw-1)$ in the denominator. However, their proof suggests that such a change of sign is indeed needed. This ambiguity has been pointed out in \cite[Remark~4.1]{BBCS16a} and our analysis gives a rigorous account in favor of this change of sign in~$\mathcal{K}^{(T)}_{22}$.
\end{rem}

\subsection{Outline of the proof}
We discuss here a brief sketch of the proof of Theorem~\ref{MainTheorem}. The organization of the rest of the paper will also be indicated concurrently. In Section~\ref{Diffoperators}, we introduce Macdonald difference ope\-ra\-tors. In the next section, we will derive contour integral formulas (see Corollary~\ref{RthMacOnPart} and Theorem~\ref{OneMainTheorem}) of the action of Macdonald $(q,q)$-difference operators on the partition function~$Z(\rho_1;\rho_2)$ of the Pfaffian Schur measure where~$\rho_1$ and~$\rho_2$ are two Schur nonnegative specializations. In Section~\ref{CorOfSinglePartition}, we find the correlation functions of the Schur measure over the set of partitions by leveraging on the contour integral formulas of Section~\ref{ContourIntegralFormula} and the fact that the Schur polynomials are the eigenfunction of the difference operators.
However, the same strategy does not work to find the correlation function of the Pfaffian Schur process over the set of all finite sequences of partitions. This is because the skew Schur functions do not always belong to the set of eigenfunctions of Macdonald difference operators. To get around this, we use the representation of the skew Schur functions in terms of the inner product of the Schur functions. Then, we let the difference operators act on the Schur functions inside the inner product. This brings in many contour integrals inside the inner products. Using bilinearity of the inner product, we pull the integrals outside and consequently, get a contour integral representation which yields the correlation functions of the Pfaffian Schur process in the same way as in the case of Pfaffian Schur measure. We elaborate on each of these steps more clearly in Section~\ref{MainStepOfProof}. We would like to point out that similar method has been also used in~\cite{A14} to deal with the case of the Schur process. Also, \cite{BCGS16} expressed the skew-Macdonald polynomials as the scalar inner product of the Macdonald polynomials to get the expectation of the multi-level observables. However, to deal with the Pfaffian version of the Schur case, we need to develop some additional tools (namely, Theorem~\ref{ConFormulaForAction}, Lemma~\ref{InnerProductLemma1} etc.) which might be of independent interests. Finally, we complete the proof of Theorem~\ref{MainTheorem} in Section~\ref{ProofOfMainTheorem}.

\section{Macdonald difference operators}\label{Diffoperators}

The main purpose of this section is to define some linear operators $\mathcal{D}$ in the space ${\rm Sym}_n$ such that $\mathcal{D}$ can be diagonalized using the Schur polynomials. More precisely, this new set of operators are expected to yield $\mathcal{D}s_\lambda=d_\lambda s_\lambda$ where $d_\lambda$ is some scalar. As a consequence, one may write
\begin{gather*}
\mathcal{D}_X Z(X;Y) =\sum_{\lambda}d_\lambda s_\lambda(X)\tau_\lambda(Y) = \mathbb{E}_{\mathrm{psm}}[d_{\lambda}].
\end{gather*}
 It was observed in \cite{Mac15} that Macdonald difference operators are diagonalized by the Schur polynomials. In what follows, we introduce Macdonald difference operators for a finite set of variables $X=\{x_1,\dots ,x_n\}$ and discuss its properties.
\bd\label{DefDiffOpp}
For any $q\in \CC$ and $1\leq i\leq n$, define the \emph{shift operators} $T_{q,i}$ by
\begin{gather*}
T_{q,i}(F)(x_1,\dots ,x_n)= F(x_1,\dots ,qx_i,\dots , x_n).
\end{gather*}
For any subset $I\subset \{1,\dots ,n\}$, define
\begin{gather*}
A_{I}(x;t)= \prod_{i\in I, j\notin I}\frac{tx_i-x_j}{x_i-x_j}.
\end{gather*}
Now, for $r=1,2,\dots , n$ define the \emph{Macdonald $(q,t)$-difference operators} as
\begin{gather*}
D^{r,q,t}_{n,X} := q^{r(r-1)/2}\sum_{|I|=r} A_I(x;t)\prod_{i\in I}T_{q,i},
\end{gather*}
where $X=\{x_1,\dots ,x_n\}$.
\ed

\begin{rem}Throughout the rest of paper, we only use Macdonald $(q,q)$-difference operators. For convenience, we refer Macdonald $(q,q)$-difference operators as Macdonald $q$-difference operators and denote $D^{r,q,q}_{n,X}$ by $D^{r,q}_{n,X}$.
\end{rem}

 The next lemma and the following discussion aim towards justifying the motivation behind the use of the Macdonald operators in stochastic processes.
\bl[{\cite[Chapter~6, Proposition~4.15]{Mac15}}] \label{lem:Mac}
For any partition $\lambda$ such that $l(\lambda)\leq n$, we have
\begin{gather*}
D^{r,q}_{n,X} s_\lambda(X)= e_r\big(q^{\lambda_1+n-1}, q^{\lambda_2+n-2},\dots , q^{\lambda_n}\big)s_{\lambda}(X),
\end{gather*}
 where $e_r$ is the $r$-th elementary symmetric function.
\el

Owing to Lemma~\ref{lem:Mac}, the action of $D^{1,q}_{n,X}$ on the partition function $Z(X;Y)$ of the Pfaffian Schur measure yields
\begin{gather}\label{eq:MacOnPart}
D^{1,q}_{n,X}Z(X;Y)= \sum_{\lambda }\tau_{\lambda}(Y)s_{\lambda}(X)\sum_{i=1}^{n} q^{\lambda_i+n-i}.
\end{gather}
Rewriting \eqref{eq:MacOnPart} in terms of the Pfaffian Schur measure, we get
\begin{gather}\label{eq:ExpecMacRelation}
\mathbb{E}_{psm}\left(\sum_{i=1}^n q^{\lambda_i +n-i}\right)= \frac{D^{1,q}_{n,X} Z(X;Y)}{Z(X;Y)}.
\end{gather}
By the successive actions of $D^{1,q}_{n,X}$ for multiple values of $q$ on the partition function $Z(X;Y)$, we arrive at
\begin{align}\label{eq:MacOnPartMultiple}
\left(\prod_{j=1}^m D^{1,q_j}_{n,X}\right)Z(X;Y) = \sum_{\lambda} \tau_\lambda(Y)s_\lambda(X) \prod_{j=1}^m \sum_{i=1}^{n} q^{\lambda_i+n-i}_{j}.
\end{align}
To find the correlation functions (see the definition in \eqref{eq:DefCorrelationFunction}) of some finite subset $T=\{t_j\}_{j\in [1,m]}$ of $\ZZ$, we need to sum the Pfaffian Schur measure of all those partitions $\lambda$ for which the set $\{\lambda_i-i\}_{i\in \NN}$ contains $T$. But, this sum is essentially given by the coefficient of $\prod_j q^{t_j+n}_j$ (up to some constant) in the right side of \eqref{eq:MacOnPartMultiple}. The goal of the next two sections is to find contour integral formulas of the coefficient of $\prod_j q^{t_j+n}_j$. We elaborate more on the correlation functions in Section~\ref{CorOfSinglePartition}.

\section{Contour integral formulas}\label{ContourIntegralFormula}
The main aim of this section is to derive a contour integral formula for the action of Macdonald difference operators on the partition function of the Pfaffian Schur measure. It can be anticipated from \eqref{eq:MacOnPartMultiple} that the action of Macdonald operators computes the expectation of $\prod_{j=1}^m \sum_{i=1}^{n} q^{\lambda_i+n-i}_{j}$ under the Pfaffian Schur measure. For similar exposition in the case of Macdonald process, see \cite[Section~2.2.3]{BorodinCorwinMac} and for Schur process, see \cite[Section~2.2]{A14}.

Consider a function $G\colon \CC^n\to \CC$ which has the following form
\begin{gather}\label{eq:PfaffianDecomposableFunstion}
G(X)= \prod_{i<j}f(x_ix_j) \prod_{i=1}^{n} g(x_i)
\end{gather}
 for $X:=(x_1,\dots ,x_n)$. We start with deriving a contour integral formula for the action of the Macdonald operators on $G(X)$ where the function $G$ has the form shown above. Similar formulas for the action of Macdonald $(q,t)$-operators recently appeared in \cite[Proposition~3.3]{GuiBorCor17}. Whenever~$f$ is the identity function, one can find an analogous formula in \cite[Proposition~2.2.11]{BorodinCorwinMac}.
\bt\label{ConFormulaForAction}
 Consider any function $G$ of the form shown in \eqref{eq:PfaffianDecomposableFunstion}. We assume $x_i\neq x_j$ and $qx_i\neq x_j$ if $i\neq j$. Let $f$ be a holomorphic function such that it is non-zero in some neighborhoods around $\{x_ix_j\}_{1\leq i,j\leq n}$ and it does not contain any poles near $\big\{q^2x_ix_j,qx_ix_j\big\}_{1\leq i,j\leq n}$ for some $q\in \CC$ with $|q|<1$. Furthermore, let us assume that $g$ is also a holomorphic function such that it is non-zero in the neighborhoods of $\{x_1,\dots, x_{n}\}$ and does not have any pole near $\{qx_1,\dots ,qx_{n}\}$. Then, we have
 \begin{gather}
 D^{r,q}_{n,X}G(X)= \frac{G(X)q^{r(r-1)/2}}{(2\pi{\rm i})^r r!(q-1)^{r}} \oint_{\mathcal{C}} \cdots \oint_{\mathcal{C}} \prod_{1\leq i\neq j\leq r}\frac{(z_i-z_j)}{(qz_i-z_j)f(qz_iz_j)}\prod_{1\leq i<j\leq r} f\big(q^2 z_iz_j\big)f(z_iz_j)\nonumber\\
 \hphantom{D^{r,q}_{n,X}G(X)=}{} \times\left(\prod_{i=1}^r \frac{f\big(z^2_i\big)}{f\big(qz^2_i\big)}\right)\prod_{i=1}^r\prod_{j=1}^{n} \frac{(qz_i - x_j)f(qz_ix_j)}{(z_i-x_j)f(z_ix_j)}\prod_{i=1}^{r} \frac{g(qz_i)}{g(z_i)}\frac{{\rm d}z_i}{z_i},
 \label{eq:MacActionContourInt}
\end{gather}
where in the right hand side, we have $r$-fold integral and the contour $\mathcal{C}$ is defined as the union $\mathcal{C}_1\cup \cdots \cup \mathcal{C}_n$ and $\mathcal{C}_i$ is a closed contour around the point $x_i$ satisfying
\begin{enumerate}\itemsep=0pt
\item[$(i)$] $q\mathcal{C}_i$ sits outside the contour $\mathcal{C}$,
\item[$(ii)$] $\mathcal{C}_i$ contains no other poles of the integrand
\end{enumerate}
for $i\in [1,n]$. For illustrations, see Fig.~{\rm \ref{fig:1}}.
\et

\begin{proof} We begin with a detailed description on the choice of the integration contours. We need $n$ non-intersecting positively oriented circles $\mathcal{C}_1,\dots ,\mathcal{C}_{n}$ centered at the points $\{x_1,\dots ,x_{n}\}$ respectively such that $q\mathcal{C}_i$ is not contained in any of the contours $\mathcal{C}_{j}$ for $j\in [1,n]$. To see why such contours exist, we choose $q\in \CC$ such that $|q|<1$ and then select $n$ positive real numbers ${s_1,\dots ,s_n}$ such that $q\mathcal{B}(x_i,s_i) \cap \mathcal{B}(x_j,s_j)$ is empty for all $i,j\in [1,n]$. This is possible since we have assumed $qx_i\neq x_j$ for $x_i\neq x_j$. Here, $\mathcal{B}(x,s)$ denotes the ball of radius $s$ centered around the point $x\in \mathbb{C}$. Now, we choose $\mathcal{C}_i$ to be the circle of radius $s_i$ centered around $x_i$. Recall that~$f$ has no zero near the points $x_ix_j$ and has no poles near at $\big\{q^2x_ix_j,x_ix_j\big\}_{1\leq i,j\leq n}$ and $g$ has no zero near the points $x_j$ and no poles near~$qx_j$. Thus, for all $1\leq i\leq n$, one can choose each of the elements of the set $\{s_1,\dots , s_n\}$ small enough such that $\mathcal{C}=\cup_i \mathcal{C}_i$ does not contain the zeros of $f(qz_iz_j)$, $f(z_ix_l)$ and $f(z_iy_l)$ for any $z_i\in \mathcal{C}_k $, $z_j\in \mathcal{C}_l$ where $k,l\in [1,n]$ and $1\leq i< j\leq r $.

Now, we turn to compute the integral in \eqref{eq:MacActionContourInt} and will show it exactly matches with the left side of \eqref{eq:MacActionContourInt}. First, we take the integral with respect to the variable $z_1$. The integral can be written as the sum of the residues at the poles enclosed by $\mathcal{C}$. Recall that $\mathcal{C}$ is union of $n$ contours $\{\mathcal{C}_i\}_i$. Moreover, $\mathcal{C}_i$ contains $x_i$ inside. But, there is no other singularity inside the contour $\mathcal{C}_i$. This is due to the facts that $(1)$ $q\mathcal{C}_j$ sits outside of $\mathcal{C}_i$ for all $j\in [1,n]$ (shows $qz_1$ lives outside $\mathcal{C}$) and $(2)$ $\mathcal{C}_i$ contains no other poles coming from the functions $f$ and $g$. Furthermore, each of the points in the set $\{x_i\}_{i\in [1,n]}$ is a simple pole of the integral. Henceforth, the integral which is the sum of $n$ residues is given as
\begin{gather}
2\pi{\rm i}(q-1)\sum_{k=1}^n\frac{g(qx_k)}{g(x_k)}\prod_{j=1,j\neq k}^n \frac{(qx_k-x_j)f(qx_kx_j)}{(x_k-x_j)f(x_kx_j)}\oint_{\mathcal{C}} \cdots \oint_{\mathcal{C}} \prod_{2\leq i\neq j\leq r}\frac{(z_i-z_j)}{(qz_i-z_j)f(qz_iz_j)} \nonumber\\
\qquad{} \times \prod_{2\leq i<j\leq r}f\big(q^2 z_iz_j\big)f(z_iz_j)\left(\prod_{i=2}^r \frac{f\big(z^2_i\big)}{f\big(qz^2_i\big)}\right)\prod_{i=2}^r\frac{(x_k-z_i)}{(qx_k-z_i)} \nonumber\\
\qquad{} \times \prod_{i=2}^r \prod_{j=1,j\neq k}^{n} \frac{(qz_i-x_j)f(qz_ix_j)}{(z_i-x_j)f(z_ix_j)}\prod_{i=2}^r\frac{g(qz_i)f\big(q^2z_ix_k\big)}{g(z_i)f(qz_ix_k)} \frac{{\rm d}z_i}{z_i}.\label{eq:IntegralAtFirstStep}
\end{gather}

 Note that for any $k\in [1,n]$, $k$-th term in the above sum no longer has any pole at the point~$x_k$ for the variables $z_i$ where $i$ now varies in $[2,r]$. Furthermore, the properties~$(i) $ and~$(ii)$ of the integration contours imply that in order to compute the remaining integrals in the $k$-th term in~\eqref{eq:IntegralAtFirstStep}, we only need to evaluate the residues at the points $\{x_1,\dots ,\hat{x}_k,\dots ,x_n\}$ where~$\hat{x}_k$ implies~$x_k$ is not present in the set. As we continue taking the integrals with respect to the other variables, this phenomenon would greatly simplifies the successive computation. For convenience, we define a set $\mathcal{S}:=\{(\mathcal{C}_{i_1},\dots ,\mathcal{C}_{i_r})\,|\,i_t\neq i_s \text{ for }t\neq s\}$. We observe that the integral in~\eqref{eq:MacActionContourInt} can be restricted over all such $r$-tuple of contours in $\mathcal{S}$. If we denote $\mathcal{A}:=\{i_1,\dots ,i_r\}$, then integral over the contour $\mathcal{C}_{i_1}\times \cdots \times\mathcal{C}_{i_r}$ with respect to the variables $z_1,\dots ,z_r$ will be equal to
\begin{gather}\label{eq:IntegralOverTheContour}
(q-1)^r(2\pi{\rm i})^r\prod_{i\in \mathcal{A}, j\notin \mathcal{A}, i<j}\frac{qx_i-x_j}{x_i-x_j}\frac{f(qx_ix_j)}{f(x_ix_j)} \prod_{i\in \mathcal{A}, j\in \mathcal{A}, i<j}\frac{f\big(q^2x_ix_j\big)}{f(x_ix_j)}\prod_{i\in \mathcal{A}} \frac{g(qx_i)}{g(x_i)}.
\end{gather}
We get the simplified form in \eqref{eq:IntegralOverTheContour} due to the telescopic cancellations of the terms at each of the stages of the integration. We demonstrate here how such cancellations come into play by computing the integrals with respect to~$z_1$ and~$z_2$. First, note that the integral with respect to~$z_1$ along~$\mathcal{C}_{i_1}$ is given by $i_1$-th term of the sum in \eqref{eq:IntegralAtFirstStep}. Furthermore, the integral with respect~$z_2$ along $\mathcal{C}_{i_2}$ will contribute a factor given as
\begin{gather*}
(q-1)(2\pi{\rm i})f\big(q^2x_{i_1}x_{i_2}\big)\frac{(x_{i_1}-x_{i_2})}{(qx_{i_1}-x_{i_2})f(qx_{i_1}x_{i_2})}\prod_{j=1,j\neq i_1,i_2}^n\frac{(qx_{i_2}-x_j)f(qx_{i_2}x_j)}{(x_{i_2}-x_j)f(x_{i_2}x_{j})}
\end{gather*}
outside the remaining integrals. Thus, the following term
\begin{gather*}
\frac{(qx_{i_1}-x_{i_2})f(qx_{i_1}x_{i_2})}{(x_{i_1}-x_{i_2})}
\end{gather*}
will get canceled from the product. This suggests that the term like
\begin{gather*}
\prod_{1\leq p<q\leq k}\frac{(qx_{i_p}-x_{i_q})f(qx_{i_p}x_{i_q})}{(x_{i_p}-x_{i_q})}
\end{gather*}
will be canceled after taking integrals with respect to $z_1,\dots ,z_k$ (for $k<r$) along the contours $\mathcal{C}_{i_1},\dots ,\mathcal{C}_{i_k}$ respectively. Owing to such repeated cancellations, we get~\eqref{eq:IntegralOverTheContour} as the end result of the integral. In fact, we get back the same value of the integral for $r!$ different choices of the contours from the set $\mathcal{S}$. This justifies the presence of the factor $r!$ in the right side of~\eqref{eq:MacActionContourInt}. Thus, the sum of the contributions like in~\eqref{eq:IntegralOverTheContour} when multiplied with $G(X)(q)^{r(r-1)/2}(2\pi{\rm i}(q-1))^{-r}(r!)^{-1}$ would exactly be equal to $D^{r,q}_{n}G(X)$. Hence, the claim follows.
\end{proof}

Owing to \eqref{eq:HallLittleSchur}, the partition function $Z(X;Y)$ for any two Schur non-negative finite specialization $X=(x_1,\dots ,x_n)$ and $Y=(y_1,\dots ,y_n)$ can be written as in \eqref{eq:PfaffianDecomposableFunstion} for $f(x)=(1-x)^{-1}$ and $g(x)= \prod_{i=1}^{n}(1-xy_i)^{-1}$. Using Theorem~\ref{ConFormulaForAction}, we get the following contour integral formula in the context of Pfaffian Schur measure.

\begin{figure}[t]\centering
\begin{tikzpicture}
 \draw[line width = 0.25mm, <->] (-1,-1) -- (7,-1) node {};
 \draw[line width = 0.25mm, <->] (0,-4) -- (0,4) node {};
 \draw[line width = 0.4mm] (0.5,0.8) circle (0.45cm);
 \draw[line width = 0.4mm] (2.5,0) circle (0.5cm);
 \node[below] at (2.5,0) (A) {$x_t$};
 \node[right] at (2.87,0.37) (B) {$\mathcal{C}_t$};
 \node[right] at (0.87,1.07) (C) {$q\mathcal{C}_t$};
 \node[below] at (0.5,0.81) (D) {$qx_t$};
 \fill (0.5,0.8) circle[radius=1pt];
 \fill (2.5,0) circle[radius=1pt];
 \draw[line width = 0.4mm] (4,2.8) circle (0.45cm);
 \draw[line width = 0.4mm] (6,2) circle (0.5cm);
 \node[below] at (6,2) (E) {$x_s$};
 \node[right] at (6.37,2.37) (F) {$\mathcal{C}_s$};
 \node[right] at (4.37,3.17) (G) {$q\mathcal{C}_s$};
 \node[below] at (4,2.81) (H) {$qx_s$};
 \fill (4,2.8) circle[radius=1pt];
 \fill (6,2) circle[radius=1pt];
\end{tikzpicture}
\caption{Plot of the contours $\{\mathcal{C}_t, q\mathcal{C}_t\}$ and $\{\mathcal{C}_s, q\mathcal{C}_s\}$.}\label{fig:1}
\end{figure}

\bc\label{RthMacOnPart}Consider the action of the Macdonald operators $D^{r,q}_{n,X}$ on the partition function $Z(X;Y)$. Fix two positive real numbers $t,s$. Assume that the following two sets of distinct non-negative real numbers $X=\{x_1,\dots ,x_{n}\}$ and $Y:=\{y_1,\dots ,y_{k}\}$ satisfy the relation
\begin{gather}\label{eq:XYRelations}
 t<\min_{1\leq j\leq n} \{x_j\}<\max_{1\leq j\leq n} \{x_j\}<s\leq \min_{1\leq j\leq n}\big\{|x_j|^{-1}\big\}\wedge \min_{1\leq j\leq n}\big\{|qy_j|^{-1}\big\}.
\end{gather}
Let us also assume that $\{x_1,\dots ,x_{n}\}$ are pairwise distinct. Thus, we get
\begin{gather}
D^{r,q}_{n,X} Z(X;Y)=\frac{Z(X;Y)(q)^{r(r-1)/2}}{r!(2\pi{\rm i})^r(q-1)^{r}}\oint_{\mathcal{C}}\cdots \oint_{\mathcal{C}}\prod_{1\leq i\neq j\leq r}\frac{(z_i-z_j)(1-qz_iz_j)}{(qz_i-z_j)}\prod_{i=1}^r \frac{1-qz^2_i}{1-z^2_i}\nonumber\\
 \times\prod_{1\leq i<j\leq r}\frac{1}{\big(1-q^2z_iz_j\big)(1-z_iz_j)}\prod_{i=1}^{r}\prod_{j=1}^{n}\frac{(qz_i-x_j)(1-z_ix_j)}{(z_i-x_j)(1-qz_ix_j)} \prod_{i=1}^r \prod_{j=1}^{k} \frac{1-z_iy_j}{1-qz_iy_j}\prod_{i=1}^r\frac{{\rm d}z_i}{z_i},\label{eq:MacOnPartition}
\end{gather}
where the contour $\mathcal{C}:=\mathcal{C}_1\cup \cdots \cup \mathcal{C}_n$ satisfies the properties mentioned in Theorem~{\rm \ref{ConFormulaForAction}}.
\ec
\begin{proof}
The condition in \eqref{eq:XYRelations} implies that there exist small neighborhoods around $\{x_1,\dots ,x_{n}\}$ which do not contain any element in the sets $\big\{(qy_1)^{-1},\dots , (qy_{k})^{-1}\big\}$ and $\big\{x_1^{-1},\dots ,x_{n}^{-1}\big\}$. Since we have assumed $x_i$'s to be distinct, one can choose non-intersecting balls around each member of the set $\{x_1,\dots ,x_{n}\}$ exactly in the same way as we did it in the last theorem. Therefore, one can get~\eqref{eq:MacOnPartition} using Theorem~\ref{ConFormulaForAction}.
\end{proof}

\bt\label{OneMainTheorem}
Consider the action of $D^{1,.}_{n, X}$ on the partition function $Z(X;Y)$ for $d$ different complex numbers namely, $q_1,\dots ,q_{d}$. For the sake of notational convenience, we assume~$X$ and~$Y$ has same number of elements $($equal to $n)$. Fix some positive real numbers $\upsilon,s$ and a set of ordered positive real numbers $\{r_{1}>\cdots >r_{d}\}$ such that $\{\upsilon,s\}\cup \{r_{1}>\cdots >r_{d}\}\cup \{q_1,\dots ,q_{d}\}$ satisfy the following relations
\begin{gather}\label{eq:XY1DMacOpRelation1}
{\rm dist} \bigg(\medcup_{\substack{i\in[1,n]\\ j\in [1,d]}}\big\{q_jx_i, q^{-1}_jx_i, (q_jx_i)^{-1}, (q_jy_i)^{-1}, x^{-1}_i\big\}, \{x_i\}_{i\in [1,n]}\bigg) \geq r_1+s
\end{gather}
and
 \begin{gather}\label{eq:XY1DMacOpRelation2}
r_1<s\upsilon^2\wedge s\min_{j\in [1,m]}\big\{|q_j|^{-1}\wedge |q_j|\big\},
\end{gather}
 where
 \begin{gather*}
 \upsilon:= \min_{\substack{1\leq j\leq m-1\\1\leq i\leq n}}\big\{\big|q^{\pm 1}_jx^{\pm 1}_i\big|\wedge x_i\big\}.
\end{gather*}
Here, ${\rm dist}(A,B)$ denotes the Euclidean distance between two sets $A$ and $B$, i.e., ${\rm dist}(A,B) = \min\limits_{x\in A,\, y\in B}\|x-y\|$. Moreover, we assume $x_i$'s are all distinct positive real numbers less than $1$ and $0\leq y_i\leq 1$ for all $i\in [1,n]$.
 Then, we have
\begin{gather}
\frac{\prod\limits_{j=1}^{d} D^{1,q_j}_{n,X} Z(X;Y)}{Z(X;Y)} = \frac{1}{(2\pi{\rm i})^{d}}\oint_{\mathcal{C}_1}\cdots \oint_{\mathcal{C}_d} \prod_{j=1}^{d} \frac{1}{q_jz_j- z_j}\prod_{j=1}^{d}\prod_{i=1}^{n}\frac{q_jz_j - x_i}{z_j-x_i}\nonumber\\
\hphantom{\frac{\prod\limits_{j=1}^{d} D^{1,q_j}_{n,X} Z(X;Y)}{Z(X;Y)} =}{} \times \prod_{1\leq j<k\leq d}\frac{(q_jz_j-q_kz_k)(z_j-z_k)}{(z_j-q_kz_k)(q_jz_j-z_k)}\prod_{1\leq j<k\leq d} \frac{(1-q_k z_kz_j)(1-q_jz_kz_j)}{(1-q_jq_k z_jz_k)(1-z_jz_k)} \nonumber\\
\hphantom{\frac{\prod\limits_{j=1}^{d} D^{1,q_j}_{n,X} Z(X;Y)}{Z(X;Y)} =}{} \times \prod_{j=1}^{d}\prod_{i=1}^{n}\frac{(1-z_jx_i)(1-z_jy_i)}{(1-q_jz_jx_i)(1-q_jz_jy_i)}\prod_{j=1}^{d}\frac{1-q_jz^2_j}{1-z_j^2}{\rm d}z_j,\label{eq:MCIntForMultipleAction}
\end{gather}
where the contour $\mathcal{C}_j$ consists small circles of radius $r_j$ around $\{x_1,\dots ,x_d\}$ for $1\leq j\leq d$.
\et
\begin{proof} In the case when $d=1$, the result follows from Corollary~\ref{RthMacOnPart}. To prove the other cases, we use induction. Let us assume that the above formula is valid for $d=m-1$. In order to complete the proof, let us define
 \begin{gather}
 Z^{(m-1)}(X;Y) := \frac{1}{(2\pi{\rm i})^{m-1}} \oint_{\mathcal{C}_1}\cdots \oint_{\mathcal{C}_{m-1}} Z(X;Y)\prod_{j=1}^{m-1} \frac{1}{q_jz_j- z_j}\prod_{j=1}^{m-1}\prod_{i=1}^{n}\frac{q_jz_j - x_i}{z_j-x_i}\nonumber\\
 \hphantom{Z^{(m-1)}(X;Y) :=}{} \times \prod_{1\leq j<k\leq m-1}\frac{(q_jz_j-q_kz_k)(z_j-z_k)}{(z_j-q_kz_k)(q_jz_j-z_k)}\prod_{1\leq j<k\leq m-1} \frac{(1-q_k z_kz_j)(1-q_jz_kz_j)}{(1-q_jq_k z_jz_k)(1-z_jz_k)} \nonumber\\
\hphantom{Z^{(m-1)}(X;Y) :=}{} \times \prod_{j=1}^{m-1}\prod_{i=1}^{n}\frac{(1-z_jx_i)(1-z_jy_i)}{(1-q_jz_jx_i)(1-q_jz_jy_i)}\prod_{j=1}^{m-1}\frac{1-q_jz^2_j}{1-z_j^2}{\rm d}z_j.\label{eq:StringFieldCoupling}
\end{gather}
We assume $\prod_{j=1}^{m-1}D^{1,q_j}_{n,X} Z(X;Y)$ is given by $Z^{(m-1)}(X;Y)$. One can note that the integrand in~\eqref{eq:StringFieldCoupling} looks like $cG(X)$ for some $G$ as in~\eqref{eq:PfaffianDecomposableFunstion} and the variable $c$ does not depend on $X$. To be precise, we can write down the integrand as the product of $\prod_{1\leq i<j\leq n}(1-x_ix_j)^{-1}$ ($=H_0(X)$) and $\prod_{i=1}^{n} g_1(x_i)$ ($=G(X)$) where
 \begin{gather*}
g_1(x):=\Xi \cdot\prod_{j=1}^{m-1}\frac{(1-z_jx)(q_jz_j-x)}{(z_j-x)(1-q_jz_jx)}\prod_{i=1}^{n}\frac{1}{1-y_ix},
 \end{gather*}
 where $\Xi$ depends on $z_1,\dots, z_{m-1}$. To compute $D^{1,q_d}_{n,X} Z^{(m-1)}_{n}$, we take the integral in~\eqref{eq:StringFieldCoupling} after letting the difference operators $D^{1,q_m}_{n}$ acts on the integrand. Note that such interchange of the integral and the summation is possible because there are only finitely many terms in the summation. Now, we need to get a contour integral representation of $D^{1,q_m}_{n,X}H_0(X)G(X)$. To begin with, we elaborate on the choice of the contours. Contour for the new integral must avoid the poles of $g_1(q_mX)/g_1(X)$ at
\begin{gather*}
\big\{(q_my_i)^{-1}\big\}_{i\in [1,n]} \cup \big\{z_j^{-1}\big\}_{j\in [1,m-1]} \cup \{q_jz_j\}_{j\in [1,m-1]} \\
\qquad{} \cup \big\{q_m^{-1} z_j\big\}_{j\in [1,m-1]} \cup \big\{(q_jz_j)^{-1}\big\}_{j\in [1,m-1]}
 \end{gather*}
 as well as the poles at $\big\{(q_mz_j)^{-1}\big\}_{j\in [1,m-1]}$ which correspond to the factor $H_0(X)$. Moreover, the contour should include the set of points $\{x_1,\dots ,x_n\}$ inside it. From the conditions, it follows that we can choose~$q_m$ such that for all $i\in [1,n]$, $q_mx_i$, $q^{-1}_mx_i$, $(q_mx_i)^{-1}$ and $(q_my_i)^{-1}$ stay outside at least at a distance $s>0$ (fixed a priori) away from the circles of radius $r_1$ centered around each of the points in the set $\{x_1,\dots ,x_n\}$. Furthermore, recall that $r_1<s\upsilon^2\wedge s\min\limits_{j\in [1,m]}\big\{q^{-1}_j\wedge q_j\big\} $.
Thus, if we would have chosen the contours of $\{z_j\}_{j\in [1,m-1]}$ to be $\{\mathcal{C}_j\}_{j\in [1,m-1]}$ where $\mathcal{C}_j$ is composed of circles of radius $r_j$ around the points~$x_i$, then none of the poles from the union
\[\big\{z_j^{-1}\big\}_{j\in [1,m-1]} \cup \{q_jz_j\}_{j\in [1,m-1]} \cup \big\{q_m^{-1} z_j\big\}_{j\in [1,m-1]} \cup \big\{(q_jz_j)^{-1}\big\}_{j\in [1,m-1]}\]
comes inside the contour $\mathcal{C}_m$ for the variable $z_m$. This is due to the following two facts: whenever $z_j=x_i+r_j{\rm e}^{{\rm i}\theta_j}$, then
\begin{enumerate}\itemsep=0pt
\item[(a)] $q_jz_j$ or $q^{-1}_jz_j$ will be within the distance $r_j\big(q_j\vee q^{-1}_j\big)$ from the point $q_jx_i$ or $q^{-1}_jx_i$, thus,
\begin{gather*}
|q_jz_j-q_jx_i|\vee \big|q^{-1}_jz_j-q^{-1}_jx_i\big|\leq r_j\big(q^{-1}_j\wedge q_j\big)\leq s,
\end{gather*}
\item[(b)] $(q_jz_j)^{-1}$ or $q_jz^{-1}_j$ will be within the distance $r_j\upsilon^{-2}$ from the point $(q_jx_i)^{-1}$ or $q_jx^{-1}_i$, thus,
\begin{gather*}
\big|(q_jz_j)^{-1}-(q_jx_i)^{-1}\big|\vee \big|q_jz^{-1}_j-q_jx^{-1}_i\big|\leq r_j\upsilon^{-2}\leq s.
\end{gather*}
\end{enumerate}
 At this point, using Theorem~\ref{ConFormulaForAction} for $r=1$ and $q=q_m$, we get
 \begin{gather*}
 D^{1,q_m}_{n,X} Z (X;Y)\prod_{i=1}^n\prod_{j=1}^{m-1}\frac{(1-z_jx_i)(q_jz_j-x_i)}{(z_j-x_i)(1-q_jz_jx_i)} \\
\qquad{}= \frac{1}{(2\pi{\rm i})}\oint_{\mathcal{C}_m} Z(X;Y)\frac{1}{q_mz_m-z_m}\prod_{j=1}^{m-1}\prod_{i=1}^n\frac{(1-z_jx_i)(q_jz_j-x_i)}{(z_j-x)(1-q_jz_jx_i)}\prod_{i=1}^{n}\frac{(q_mz_m-x_i)}{(z_m-x_i)}\\
\qquad\quad{} \times \prod_{j=1}^{m-1}\frac{(q_jz_j-q_mz_m)(z_j-z_m)}{(q_jz_j-z_m)(z_j-q_mz_m)}\prod_{j=1}^{m-1}\frac{(1-q_mz_mz_j)(1-q_jz_jz_m)}{(1-q_mq_jz_mz_j)(1-z_mz_j)}\\
\qquad\quad{} \times \frac{1-q_mz^2_m}{1-z^2_m}\prod_{i=1}^{n}\frac{(1-z_mx_i)(1-z_my_i)}{(1-q_mz_mx_i)(1-q_mz_my_i)}{\rm d}z_m.
\end{gather*}
Combining \eqref{eq:StringFieldCoupling} and \eqref{eq:MacActionContourInt} yields the contour formula for $\prod_{j=1}^{m} D^{1,q_j}_{n,X}Z(X;Y)$ as given in~\eqref{eq:MCIntForMultipleAction} with $d=m$. This completes the proof.
\end{proof}

One can also derive the analogues of Corollary~\ref{RthMacOnPart} and Theorem~\ref{OneMainTheorem} for the partition function of the \emph{Schur process}. For that, one needs to plug in $f(x)=1$ and $g(x)= \prod_{j=1}^{n} (1-xy_j)^{-1}$ in Theorem~\ref{ConFormulaForAction}. One can see \cite[Propositions~2.2.4 and 2.2.7]{A14} for the details of the derivation. For the purpose of using it later, we just state the result here.

\bt\label{ResultForOnlySchur}
 Consider the action of $D^{1,.}_{n,X}$ on the partition function $F(X;Y)$ for $d$ different complex numbers namely, $q_1, \dots, q_d$. Fix $s>0$ and a set of ordered positive real numbers $\{r_1>\cdots > r_d\}$ such that the union $\{s\}\cup\{r_1>\cdots >r_d\}\cup \{q_1,\dots ,q_d\}$ satisfies the following relations
 \begin{gather*}
{\rm dist}\bigg(\medcup_{\substack{i\in [1,d]\\j\in [1,m]}} \big\{(q_jy_i)^{-1},q_jx_i,q_j^{-1}x_i\big\}, \{x_i\}_{i\in [1,n]}\bigg)>r_1+s
 \end{gather*}
 and $\min\limits_{i\in [1,n]}\{|x_i|\}>r_1$.
 Assume that $\{x_i\}_{i\in [1,n]}$ and $\{y_i\}_{1\in [1,n]}$ are two sets of distinct positive real numbers less than $1$. Then,
 \begin{gather}
\frac{\prod\limits_{j=1}^{d} D^{1,q_j}_{n,X} F(X;Y)}{F(X;Y)} = \frac{1}{(2\pi{\rm i})^{d}}\oint_{\mathcal{C}_1}\cdots \oint_{\mathcal{C}_d} \prod_{j=1}^{d} \frac{1}{q_jz_j- z_j}\prod_{j=1}^{d}\prod_{i=1}^{n}\frac{q_jz_j - x_i}{z_j-x_i}\nonumber\\
\hphantom{\frac{\prod\limits_{j=1}^{d} D^{1,q_j}_{n,X} F(X;Y)}{F(X;Y)} =}{} \times \prod_{1\leq j<k\leq d}\frac{(q_jz_j-q_kz_k)(z_j-z_k)}{(z_j-q_kz_k)(q_jz_j-z_k)} \prod_{j=1}^{d}\prod_{i=1}^{n}\frac{(1-z_jy_i)}{(1-q_jz_jy_i)}\prod_{j=1}^d {\rm d}z_j,\label{eq:MCIntForMultipleActionOnSchurPartition}
\end{gather}
where the contour $\mathcal{C}_j$ consists small circles of radius $r_j$ around $\{x_1,\dots ,x_d\}$ for $1\leq j\leq d$.
\et

\section{Correlation functions of a single partition}\label{CorOfSinglePartition}
In this section, we find out the correlation functions of a single partition under the \emph{Pfaffian Schur measure}. We will primarily use finite Schur non-negative specializations in the definition of the Pfaffian Schur measure given in~\eqref{eq:DefPfaffianSchurMeasure}. For the sake of notational convenience, we assume both the specializations in \eqref{eq:DefPfaffianSchurMeasure} have the same length~$n$. To motivate the discussion of this section, we start with the definition of the correlation functions for the single partition. For any random partition $\lambda=(\lambda_1,\lambda_2,\dots)$ where $\lambda_1\geq \lambda_2\geq \cdots$ and a fixed set of integers $\{t_1,t_2, \dots ,t_j\}$, the correlation function $\rho_{psm}(t_{1}, t_{2},\dots ,t_{j})$ is defined as
\begin{gather}\label{eq:CorrelationFunction}
\rho_{psm}(t_{1},\dots ,t_{j}) = \mathbb{P} (\lambda_{i_1}-i_1=t_{1},\dots , \lambda_{i_j}-i_j=t_{j}\text{ for some } i_1,\dots ,i_j\in [1,n] ).
\end{gather}
 Note that $\rho_{psm}(t_{1},\dots , t_{d})$ is the sum of the coefficients of the terms like $\prod_{j=1}^d q^{t_{j}+n}_{j}$ in the r.h.s.\ of~\eqref{eq:MacOnPartMultiple}. In what follows, we present a simple recipe to obtain the correlation functions based on this eloquent fact.
\bt
Consider the Pfaffian Schur measure with two finite Schur non-negative specializations $X:=\{x_1,\dots ,x_n\} $ and $Y:=\{y_1,\dots ,y_n\}$ such that $\sum_{\lambda}s_\lambda(X)\tau_\lambda(Y)$ is absolutely convergent. Then, the correlation function $\rho_{psm}(T)$ of the set $T:=\{t_1,\dots ,t_d\}$ is the coefficient of $Q^{T}$ in the series expansion of $Z^{-1}(X;Y)\prod_{j=1}^{d}D^{1,q_j}_{n,X}Z(X;Y)$ where $Q^{T}= \prod_{j=1}^{d} q_j^{t_j+n}$. Furthermore, one can write
\begin{gather}\label{eq:ContourIntegralForCorrelation}
\rho_{psm}(T)= \frac{1}{(2\pi{\rm i})^d}\oint_{\mathcal{C}^{\prime}}\cdots \oint_{\mathcal{C}^\prime} Z^{-1}(X;Y)\prod_{j=1}^{d}D^{1,q_j}_{n,X}Z(X;Y) q_j^{-t_j-n-1}{\rm d}q_j,
\end{gather}
where the contour $\mathcal{C}^\prime$ is a circle of radius $r_{q}<1$ oriented in the anticlockwise direction around~$0$.
\et
\begin{proof}
The fact that $\rho_{psm}(T)$ is equal to the coefficient of $Q^{T}$ in $Z^{-1}(X,Y)\prod_{j=1}^{d}D^{1,q_j}_{n,X}Z(X;Y) $ follows by combining \eqref{eq:ExpecMacRelation} with~\eqref{eq:MacOnPartMultiple}.

Now, we turn to show the right hand side of \eqref{eq:CorrelationFunction} is equal to $\rho_{psm}(T)$. Since the contour $\mathcal{C}^{\prime}$ is a circle of radius $r_q<1$ oriented counterclockwise around $0$,
the absolute value of $\prod_{j=1}^d \sum_{k=1}^n q^{\lambda_k+n-k}_j$ is bounded above by $nr^{d}_q$ for any partition $\lambda$. By our assumption, $\sum_{\lambda}s_\lambda(X)\tau_\lambda(Y)$ is absolutely convergent. Therefore, the following series
\begin{gather}\label{eq:AbsoluteConvergenceOfSeries}
\sum_{\lambda\in \mathbb{Y}} s_\lambda(X)\tau_\lambda(Y)\prod_{j=1}^d \sum_{k=1}^n q^{\lambda_k+n-k}_j
\end{gather}
is also absolutely convergent. Now, we substitute $Z^{-1}(X,Y)\prod_{j=1}^{d}D^{1,q_j}_{n,X}Z(X;Y)$ by the sum in~\eqref{eq:AbsoluteConvergenceOfSeries} (since this is equal to $\mathbb{E}_{psm} \big[\prod_{j=1}^d \sum_{k=1}^n q^{\lambda_k+n-k}_j\big]$) inside the integrals of~\eqref{eq:AbsoluteConvergenceOfSeries} and interchange\footnote{This interchange is justified by the dominated convergence theorem.} the sum (over $\lambda$) and the integrals.
For any partition~$\lambda$, the integrals with respect to $q_1, \dots, q_d$ will be equal to $Z^{-1}(X;Y)s_\lambda(X)\tau_\lambda(Y)$ if $\lambda_{i_1}-i_1=t_{1},\dots , \lambda_{i_d}-i_d=t_{d}$ for some $i_1,\dots ,i_d\in [1,n]$ and otherwise, it will be equal to~$0$. Hence, the sum of those integrals over all partitions $\lambda \in \mathbb{Y}$ is equal to $\rho_{\mathrm{psm}}(T)$. This completes the proof.
\end{proof}

The rest of the discussion in this section is devoted to prove that the contour integral formula of the correlation functions in \eqref{eq:ContourIntegralForCorrelation} can be expressed as the Pfaffian of a $2d\times 2d$ skew symmetric matrix. For our next result, consider the Pfaffian Schur measure given in \eqref{eq:DefPfaffianSchurMeasure} with two Schur nonnegative specializations $X$ and $Y$ where $X=\{x_1,x_2,\dots ,x_{n}\}$ and $Y=\{y_1,\dots ,y_{n}\}$ are two sets of positive real numbers less than $1$ and the series $\sum_{\lambda}s_\lambda(X)\tau_\lambda(Y)$ is absolutely convergent.
\bt\label{TheoremForSinglePartition}
 Suppose $T=\{t_1,\dots ,t_{d}\}\subset \ZZ$ are pairwise disjoint integers. Assume $n$ is an integer greater than $\max\{d,d-\min T\}$. Then, we have $\rho_{psm}(T)=\Pf\big(\mathcal{K}^{(T)}\big)$ where $\mathcal{K}^{(T)}$ is a~$2d\times 2d$ skew symmetric matrix $($see Remark~{\rm \ref{rem:StructureOfK}} for more description on the structure of the matrix $\mathcal{K}^{(T)})$. The coefficients of $\mathcal{K}^{(T)}$ are given explicitly by the contour integrals shown as follows:
\begin{gather*}
\mathcal{K}^{(T)}_{11}(k,l):= \frac{1}{(2\pi{\rm i})^2}\oint_{\mathfrak{C}_1}\oint_{\mathfrak{C}_2} \frac{(z-w)}{\big(z^2-1\big)\big(w^2-1\big)(zw-1)} \\
\hphantom{\mathcal{K}^{(T)}_{11}(k,l):=}{} \times \prod_{i=1}^{n} \frac{\big(1-z^{-1}x_i\big)\big(1-w^{-1}x_i\big)\big(1-z^{-1}y_i\big)\big(1-w^{-1}y_i\big)}{(1-zx_i)(1-wx_i)} \frac{{\rm d}z{\rm d}w}{z^{t_k}w^{t_l}},
\end{gather*}
where $\mathfrak{C}_1$ and $\mathfrak{C}_2$ are positively oriented closed contours containing the singularities at $\big\{\!x^{-1}_i\!\big\}_{i\in [1,n]}$ and furthermore, $|z|>1$ and $|w|>1$ for all $z\in \mathfrak{C}_1$ and $w\in \mathfrak{C}_2$ respectively,
\begin{gather*}
\mathcal{K}^{(T)}_{21}(k,l):= -\frac{1}{(2\pi{\rm i})^2}\oint_{\mathfrak{C}_3} \oint_{\mathfrak{C}_4} \frac{(z-w)}{w\big(z^2-1\big)(zw-1)}\\
\hphantom{\mathcal{K}^{(T)}_{21}(k,l):=}{}\times\prod_{i=1}^{n}\frac{\big(1-w^{-1}x_i\big)\big(1-z^{-1}x_i\big)\big(1-z^{-1}y_i\big)}{(1-zx_i)(1-wx_i)(1-wy_i)}\frac{{\rm d}z {\rm d}w}{z^{t_k}w^{t_l}},\\
\mathcal{K}^{(T)}_{12}(k,l) :=- \mathcal{K}^{(T)}_{21}(l,k),
\end{gather*}
and
\begin{gather*}\label{eq:PfaffianMatrixExpression22}
\mathcal{K}^{(T)}_{22}(k,l):= \frac{1}{(2\pi{\rm i})^2}\oint_{\mathfrak{C}_5}\oint_{\mathfrak{C}_6} \frac{(z-w)}{zw(zw-1)}\prod_{i=1}^{n}\frac{\big(1-z^{-1}x_i\big)\big(1-w^{-1}x_i\big)}{(1-zx_i)(1-wx_i)(1-zy_i)(1-wy_i)}\frac{{\rm d}z{\rm d}w}{z^{t_k}w^{t_l}}.
\end{gather*}
All of $\mathfrak{C}_3$, $\mathfrak{C}_4$, $\mathfrak{C}_5$ and $\mathfrak{C}_6$ are closed, positively oriented contours. Furthermore, $\mathfrak{C}_3$, $\mathfrak{C}_4$ contain the singularities at $\big\{x^{-1}_i\big\}_{i\in [1,n]}$ inside. In addition, $\mathfrak{C}_3$ and $\mathfrak{C}_4$ in $\mathcal{K}^{(T)}_{12}(r,s)$ are such that $|z|>1$ and $|zw|>1$ are satisfied for all $(z,w)\in \mathfrak{C}_3\times \mathfrak{C}_4$ and in the case of $\mathcal{K}^{(T)}_{22}$, $z$ and $w$ must satisfy $|z|<1$ and $|w|<1$ for all pairs $(z,w)$ from the set $\mathfrak{C}_5\times \mathfrak{C}_6$.
\et

\begin{rem}\label{rem:StructureOfK}
The matrix $\mathcal{K}^{(T)}$ of Theorem~\ref{TheoremForSinglePartition} consists of a number $d^2$ of $2\times 2$ block matrices. For any $k,l\in [1,d]$, the $(k,l)$ block matrix of the matrix $\mathcal{K}^{(T)}$ is denoted as $\mathcal{K}^{(T)}\langle k,l\rangle$. One can completely determine $\mathcal{K}^{(T)}$ by specifying $\mathcal{K}^{(T)}\langle k,l\rangle$ for all $1\leq k,l \leq d$ where
\begin{gather*}
\mathcal{K}^{(T)}\langle k,l\rangle:=\begin{pmatrix}
\mathcal{K}^{(T)}_{11}(k,l) & \mathcal{K}^{(T)}_{12}(k,l)\\
\mathcal{K}^{(T)}_{21}(k,l) & \mathcal{K}^{(T)}_{22}(k,l)
\end{pmatrix}.
\end{gather*}
For instance, if $\mathcal{K}^{(T)}_{8\times 8}$ is given as follows
\begin{gather*}
\mathcal{K}^{(T)}=\left[
	\begin{array}{@{}cccc|ccccc@{}}
	* & * & \phantom{\Delta}* & \phantom{\Delta}* & \phantom{\Delta}* & \phantom{\Delta}* & \phantom{\Delta}* & \phantom{\Delta}* \\
	* & * & \phantom{\Delta}* & \phantom{\Delta}* & \phantom{\Delta}* & \phantom{\Delta}* & \phantom{\Delta}* & \phantom{\Delta}* \\
	 \Delta_{11} & \Delta_{12} & \phantom{\Delta}* & \phantom{\Delta}* & \phantom{\Delta}* & \phantom{\Delta}* & \phantom{\Delta}* & \phantom{\Delta}* \\
	 \Delta_{21} & \Delta_{22} & \phantom{\Delta}* & \phantom{\Delta}* & \phantom{\Delta}* & \phantom{\Delta}* & \phantom{\Delta}* & \phantom{\Delta}*\\
	\hline
	* & * & \phantom{\Delta}* & \phantom{\Delta}* & \phantom{\Delta}* & \phantom{\Delta}* & \phantom{\Delta}* & \phantom{\Delta}* \\
	* & * & \phantom{\Delta}* & \phantom{\Delta}* & \phantom{\Delta}* & \phantom{\Delta}* & \phantom{\Delta}* & \phantom{\Delta}* \\
	* & * & \phantom{\Delta}* & \phantom{\Delta}* & \phantom{\Delta}* & \phantom{\Delta}* & \phantom{\Delta}* & \phantom{\Delta}* \\
	* & * & \phantom{\Delta}* & \phantom{\Delta}* & \phantom{\Delta}* & \phantom{\Delta}* & \phantom{\Delta}* & \phantom{\Delta}*
	\end{array}
	\right],
\end{gather*}
then, $\mathcal{K}^{(T)}\langle 2,1\rangle= \begin{pmatrix}
\Delta_{11} & \Delta_{12}\\
\Delta_{21} & \Delta_{22}
\end{pmatrix}$.
\end{rem}

\begin{proof}[Proof of Theorem~\ref{TheoremForSinglePartition}]
To begin with, let us single out the following product
\begin{gather}\label{eq:PrinciplePfaffian}
\prod_{j=1}^d \frac{1-q_jz^2_j}{z_j- q_jz_j}\prod_{1\leq j<k\leq d}\frac{(q_jz_j-q_kz_k)(z_j-z_k)}{(z_j-q_kz_k)(q_jz_j-z_k)}\prod_{1\leq j<k\leq d}\frac{(1-q_kz_kz_j)(1-q_jz_kz_j)}{(1-q_jq_kz_jz_k)(1-z_jz_k)}
\end{gather}
from the integrand in \eqref{eq:MCIntForMultipleAction}. We claim that the product in~\eqref{eq:PrinciplePfaffian} is the Pfaffian of a skew symmetric matrix~$\mathcal{M}$ where
\begin{gather*}
\mathcal{M}\langle k,l\rangle= \begin{pmatrix}
\mathcal{M}_{11}(k,l) & \mathcal{M}_{12}(k,l)\\
\mathcal{M}_{21}(k,l) & \mathcal{M}_{22}(k,l)
\end{pmatrix} = \begin{pmatrix}
\dfrac{z_k-z_l}{1-z_kz_l} & \dfrac{1-q_lz_lz_k}{z_k-q_lz_l}\vspace{1mm}\\
-\dfrac{1-q_kz_kz_l}{z_l-q_kz_k} & \dfrac{q_kz_k-q_lz_l}{1-q_lq_kz_lz_k}
\end{pmatrix}.
\end{gather*}

 To see this, let us substitute\footnote{This substitution is suggested by Guillaume Barraquand.}
 \begin{gather*}
 z_j= u_{2j-1} \qquad \text{and} \qquad q_jz_j= u^{-1}_{2j},
 \end{gather*}
 where $U=\{u_1,\dots ,u_{2d}\}$ is a set of complex numbers. Under this substitution, the product in~\eqref{eq:PrinciplePfaffian} is transformed as
\begin{gather*}
\prod_{1\leq j<k\leq 2d} \frac{u_j-u_k}{1-u_ju_k}
\end{gather*}
and the matrix $\mathcal{M}$ is transformed to $\tilde{\mathcal{M}}$ where
\begin{gather*}
\tilde{\mathcal{M}}= \left(\left(\frac{u_j-u_k}{1-u_ju_k}\right)_{j,k\in[1,2d]}\right).
\end{gather*}
It remains to show that
\begin{gather}\label{eq:PfaffianIdentityToShow}
\Pf(\tilde{\mathcal{M}})= \prod_{1\leq j<k\leq 2d} \frac{u_j-u_k}{1-u_ju_k}.
\end{gather}
This follows from the Pfaffian Schur identity (see \cite[p.~225]{Schur1904}, \cite[p.~104]{STEM}, \cite{SunT}).

Recall that $\rho_{psm}(T)$ is expressed in Theorem~\ref{TheoremForSinglePartition} as a contour integral where the integrand is written in terms of $Z^{-1}(X,Y)\prod_{j=1}^{2d}D^{1,q_j}_{n,X}Z(X;Y)$. Substituting the contour integral formula for $Z^{-1}(X,Y)\prod_{j=1}^{2d}D^{1,q_j}_{n,X}Z(X;Y)$ \eqref{eq:MCIntForMultipleAction} inside the integral of the r.h.s.\ of~\eqref{eq:ContourIntegralForCorrelation} yields
\begin{gather}
\rho_{psm}(T)= \frac{1}{(2\pi{\rm i})^d}\oint_{\mathcal{C}^\prime}{\rm d}q_1 \cdots \oint_{\mathcal{C}^\prime}{\rm d} q_d \oint_{\mathcal{C}_1}{\rm d} z_1\cdots \oint_{\mathcal{C}_{d}} {\rm d} z_d\Pf(\mathcal{M})\prod_{j=1}^{d}\frac{q^{-t_j-n-1}_j}{1-z^2_j}\nonumber\\
\hphantom{\rho_{psm}(T)=}{} \times\prod_{j=1}^{d}\prod_{i=1}^{n}\frac{(q_jz_j-x_i)(1-z_jx_i)(1-z_jy_i)}{(z_j-x_i)(1-q_jz_jx_i)(1-q_jz_jy_i)}. \label{eq:CorrelationFurther}
\end{gather}
 One may ask whether the conditions of Theorem~\ref{OneMainTheorem} are satisfied in the above substitution. To satisfy the conditions in Theorem~\ref{OneMainTheorem}, we first restrict the contours $\mathcal{C}^\prime$ of the variables $q_j$. After restricting $\mathcal{C}^\prime$, we show that there exists contours $\{\mathcal{C}_j\}_{j\in [1,d]}$ such that those satisfy the conditions \eqref{eq:XY1DMacOpRelation1} and \eqref{eq:XY1DMacOpRelation2}. To see how it is possible, we first choose a real number $\eta>0$ such that
\begin{gather} \label{eq:ChoiceOfS}
 \eta<\frac{1}{3}\Big(\min_{i\neq j\in[1,n]}\big\{|x_i-x_j|\wedge \big|x_i-x^{-1}_j\big|\big\} \wedge \min_{i,j\in [1,n]}\big\{|x_i-y_j|\wedge \big|x_i-y^{-1}_j\big|\big\}\Big).
\end{gather}
where we have assumed that the right hand side of~\eqref{eq:ChoiceOfS} is strictly positive. Now, we define $\upsilon_0:= \min_{i\in [1,n]}\{|x_i|\wedge |x_i|^{-1}\wedge |y_i|\wedge |y_i|^{-1}\}$. Note that $\upsilon_0\leq 1$. We choose the contour $\mathcal{C}^\prime$ to be the circle of radius $1-\eta\upsilon_0$ around $0$. Using nonnegativity of the variables $\{x_i\}_{i\in [1,n]}$, we get
\begin{gather*}
|q_jx_k-x_i|\geq \max\{|\mbox{Re}(q_jx_k)-x_i|,|\mbox{Im}(q_jx_k)|\}.
\end{gather*}
Furthermore, we know $|\mbox{Im}(q_jx_k)|\geq (1-\eta\upsilon)x_k-|\mbox{Re}(q_jx_k)|>0$. As a consequence, we get $|q_jx_k-x_i|\geq 2^{-1}|(1-\eta\upsilon)x_k-x_i|$ owing to the inequality $\max\{|a|,|b|\}\geq 2^{-1}|a+b|$ for any $a,b\in \RR$. One can also see that we have $|(1-s\upsilon_0)x_k-x_i|\geq|x_k-x_i|-\eta\upsilon_0|x_k|>2\eta$ whenever $k\neq i$. Moreover, $q_jx_i$ and $x_i$ must be at least $\eta\upsilon_0|x_i|$ far apart from each other. Henceforth, on defining $s:= 2^{-1}\eta \min_{i}\{\upsilon_0 x_i\}$, we get
\begin{gather*}
{\rm dist}\bigg(\medcup_{\substack{j\in [1,d]\\i\in [1,n]}}\{q_jx_i\}, \{x_i\}_{i\in [1,n]}\bigg)>2s.
\end{gather*}
 Likewise, for the same choice of $s$ and $\upsilon_0$, one can prove a more general statement given as
 \begin{gather}\label{eq:FullConditionSatisfiability}
 {\rm dist} \bigg(\medcup_{\substack{i\in[1,n]\\ j\in [1,d]}}\big\{q_jx_i, q^{-1}_jx_i, (q_jx_i)^{-1}, (q_jy_i)^{-1}, x^{-1}_i\big\}, \{x_i\}_{i\in [1,n]}\bigg) \geq 2s.
 \end{gather}
 Now, note that the minimum value of $|q^{\pm 1}_jx^{\pm 1}_i|$ over all values of $j\in [1,d]$ and $i\in [1,n]$ lies very close to $\upsilon_0$. Now we fix a positive real number $r_1$ such that $r_1<s\upsilon^2_0$. This choice of $r_1$ satisfies the condition~\eqref{eq:XY1DMacOpRelation2}. Furthermore, owing to~\eqref{eq:FullConditionSatisfiability} and that fact that $\upsilon_0<1$, the condition~\eqref{eq:XY1DMacOpRelation1} is also satisfied. Now, we are ready to choose the contours $\{\mathcal{C}_j\}^{d}_{j=1}$. As in Theorem~\ref{OneMainTheorem}, we choose $\mathcal{C}_j$ for $j\in [1,d]$ composed of $n$ non-intersecting circles around the points $\{x_1,\dots ,x_n\}$ of radius $r_j$ from the set $\{r_1,\dots ,r_{d}\}$ where the choice of $r_1$ is given above and $\{r_2, \dots , r_{d}\}$ are chosen from $(0,r_1)$ under the restriction $r_d<r_{d-1}<\cdots <r_1$. Now, we focus on $\Pf(\mathcal{M})$. We may expand $\Pf(\mathcal{M})$ as
\begin{gather}\label{eq:PfaffianExpanOfMmatrix}
\Pf(\mathcal{M})=\frac{1}{2^d d!}\sum_{\sigma\in \mathfrak{S}(2d)}(-1)^{\sigma}\prod_{j=1}^{d} \mathcal{M}(\sigma(2j-1),\sigma(2j)).
\end{gather}
Substituting \eqref{eq:PfaffianExpanOfMmatrix} inside the integral of \eqref{eq:CorrelationFurther}, we see that the integrand can be written as
\begin{gather}
\frac{1}{2^d d!}\sum_{\sigma\in \mathfrak{S}(2d)}(-1)^{\sigma}\prod_{j=1}^{d} \mathcal{M}(\sigma(2j-1), \sigma(2j))\prod_{j=1}^{d}\frac{q^{-t_j-n-1}_j}{z^2_j-1}\nonumber\\
\qquad{} \times \prod_{j=1}^{d}\prod_{i=1}^{n}\frac{(q_jz_j-x_i)(1-z_jx_i)(1-z_jy_i)}{(z_j-x_i)(1-q_jz_jx_i)(1-q_jz_jy_i)}.\label{eq:ExapnsionWithOthertermsInt}
\end{gather}

Let us fix a permutation $\sigma$ and consider the term in the above summation indexed by $\sigma$. Let us make a substitution $\tilde{z}_j=(q_jz_j)^{-1}$.
Notice that each term of the sum~\eqref{eq:ExapnsionWithOthertermsInt} can be written in terms of the following set of variables $\{z_1,\dots ,z_{d}, \tilde{z}_{1}, \dots , \tilde{z}_{d}\}$. The Jacobian of the transformation is given by $\prod_{j=1}^{d}\big({-}q^{-2}_j\big)z^{-1}_j$. We need to multiply the integrand with the inverse of the Jacobian which we do after rewriting $\prod_{j=1}^{d}\big({-}q^{-2}_j\big)z^{-1}_j$ as $\prod_{j=1}^{d}\big({-}q^{-1}_j\big)\tilde{z}_j$. This will kill the factor $q_j^{-1}$ of the term $q_j^{-t_j-n-1}$ in \eqref{eq:ExapnsionWithOthertermsInt}. Furthermore, we note
\begin{gather*}
\frac{q_jz_j-x_i}{z_j-x_i}=q_j\frac{1-\tilde{z}_j x_i}{1-z^{-1}_j x_i}, \qquad \forall\, i\in [1,n],\quad j\in [1,d].
\end{gather*}
This piles up a factor $q^{n}_j$ in the integrand and thus, accounts for the cancellation of the term $\prod_{j=1}^{d}q^{-n}_j$ from the product in the r.h.s.\ of \eqref{eq:ExapnsionWithOthertermsInt}. Now, we substitute $\prod_{j=1}^{d} q_j^{-t_j}= \prod_{j=1}^d z^{t_j}_j\tilde{z}^{t_j}_j$ and as consequence, the r.h.s.\ of \eqref{eq:ExapnsionWithOthertermsInt} is expressed completely in terms of $\{z_1,\dots ,z_{d}, \tilde{z}_{1}, \dots , \tilde{z}_{d}\}$. Our next task is to derive the elements of the matrix $\mathcal{K}^{(T)}$ by manipulating the terms in the product in the r.h.s.\ of~\eqref{eq:ExapnsionWithOthertermsInt}. For notational convenience, we denote $\lceil \sigma(t)/2\rceil $ by $\sigma_2(t)$ for any $t\in [1,2d]$ in the rest of the proof. If both $\sigma(2j-1)$ and $\sigma(2j)$ are odd integers, then $\mathcal{M}(\sigma(2j-1),\allowbreak \sigma(2j))= (z_{\sigma_2(2j-1)}-z_{\sigma_2(2j)})(1-z_{\sigma_2(2j-1)}z_{\sigma_2(2j)})^{-1}$. Putting together $\mathcal{M}(\sigma(2j-1),\sigma(2j))$ and other terms involving $z_{\sigma_2(2j-1)}$ and $z_{\sigma_2(2j)}$ from the product, we arrive at
\begin{gather}
-\oint \oint \frac{{\rm d}z_{\sigma_2(2j-1)}{\rm d}z_{\sigma_2(2j)}}{z_{\sigma_2(2j-1)}^{-t_{\sigma_2(2j-1)}} z_{\sigma_2(2j)}^{-t_{\sigma_2(2j)}}}\frac{(z_{\sigma_2(2j-1)}-z_{\sigma_2(2j)})}{\big(1-z_{\sigma_2(2j-1)}^2\big) \big(1-z_{\sigma_2(2j)}^2\big)(1-z_{\sigma_2(2j-1)}z_{\sigma_2(2j)})}\nonumber\\
\qquad{} \times \prod_{i=1}^{n} \frac{(1-z_{\sigma_2(2j-1)}x_i)(1-z_{\sigma_2(2j)}x_i)(1-z_{\sigma_2(2j-1)}y_i)(1-z_{\sigma_2(2j)}y_i)}{\big(1-z_{\sigma_2(2j-1)}^{-1}x_i\big)\big(1-z_{\sigma_2(2j)}^{-1}x_i\big)} , \label{eq:FormOfTheFirstEntry}
\end{gather}
where the contours of the integral are oriented anticlockwise and contain the singularities at the points $\{x_i\}_{i\in [1,n]}$. In addition, the following two conditions $|z_{\sigma_2(2j-1)}|<1$ and $|z_{\sigma_2(2j)}|<1$ are satisfied over all region on the contours. This allows us to avoid the poles at $\pm 1$ and $z^{-1}_{\sigma_2(2j-1)}$ (or~$z^{-1}_{\sigma_2(2j)}$) for the variable $z_{\sigma_2(2j)}$ (or $z_{\sigma_2(2j-1)}$). Furthermore, the $(-)$ sign which appears in~\eqref{eq:FormOfTheFirstEntry} comes from the Jacobian of the transformations $(z_j,q_j)\mapsto \big(z_j,(q_jz_j)^{-1}\big)$. In the next step, we substitute the variables $z_{\sigma_2(2j-1)}$ with $z^{-1}$ and $z_{\sigma_2(2j)}$ to $w^{-1}$. This leads to the expression
\begin{gather}
\mathcal{K}^{(T)}_{11}(\sigma_2(2j-1), \sigma_2(2j))= \oint_{\mathfrak{C}_1}\oint_{\mathfrak{C}_2} \frac{{\rm d}z{\rm d}w}{z^{t_{\sigma_2(2j-1)}}w^{t_{\sigma_2(2j)}}}\frac{(z-w)}{\big(z^2-1\big)\big(w^2-1\big)(zw-1)}\nonumber \\
\hphantom{\mathcal{K}^{(T)}_{11}(\sigma_2(2j-1), \sigma_2(2j))=}{} \times\prod_{i=1}^{n} \frac{(1-z^{-1}x_i)\big(1-w^{-1}x_i\big)\big(1-z^{-1}y_i\big)\big(1-w^{-1}y_i\big)}{(1-zx_i)(1-wx_i)} .\label{eq:FormOfTheFirstEntryModified}
\end{gather}
Due to the above substitution, the contours $\mathfrak{C}_1$ and $\mathfrak{C}_2$ are closed containing all the singularities at $\big\{x^{-1}_i\big\}_{i\in [1,n]}$ and furthermore, $|z|>1$ and $|w|>1$ must be satisfied for all $(z,w)\in \mathfrak{C}_1\times \mathfrak{C}_2$.

Next, we consider the case when $\sigma(2j-1)$ is even and $\sigma(2j)$ is odd. One can also consider the opposite scenario namely, $\sigma(2j-1)$ is odd and $\sigma(2j)$ is even. In both of these two cases, the analysis would be similar. For brevity, we would only discuss the first case. Rewriting $\mathcal{M}(\sigma(2j-1),\sigma(2j))$ in terms of $\{z_1,\dots ,z_{d}, \tilde{z}_{1}, \dots , \tilde{z}_{d}\}$ yields \[\mathcal{M}(\sigma(2j-1),\sigma(2j))=\frac{1-z_{\sigma_2(2j)}\tilde{z}^{-1}_{\sigma_2(2j-1)}}{z_{\sigma_2(2j)}-\tilde{z}^{-1}_{\sigma_2(2j-1)}}.\]
Owing to this new form of $\mathcal{M}(\sigma(2j-1),\sigma(2j))$, $\mathcal{K}^{(T)}_{21}(\sigma_2(2j-1), \sigma_2(2j))$ can be written as
\begin{gather*}
-\oint \oint \frac{{\rm d}z_{\sigma_2(2j-1)}{\rm d}z_{\sigma_2(2j)}}{z_{\sigma_2(2j-1)}^{-t_{\sigma_2(2j-1)}}z_{\sigma(2j)}^{-t_{\sigma_2(2j)}}} \frac{\big(1-z_{\sigma_2(2j)}\tilde{z}^{-1}_{\sigma_2(2j-1)}\big)}{\tilde{z}_{\sigma_2(2j-1)}\big(1-z^2_{\sigma_2(2j)}\big) \big(z_{\sigma_2(2j)}-\tilde{z}^{-1}_{\sigma_2(2j-1)}\big)}\nonumber\\
\qquad{} \times \prod_{i=1}^n\frac{(1-\tilde{z}_{\sigma_2(2j-1)}x_i)(1-z_{\sigma_2(2j)}x_i)(1-z_{\sigma_2(2j)}y_i)}{\big(1-z^{-1}_{\sigma_2(2j)}x_i\big) \big(1-\tilde{z}^{-1}_{\sigma_2(2j-1)}x_i\big)\big(1-\tilde{z}^{-1}_{\sigma_2(2j-1)}y_i\big)},\label{eq:FormOfTheSecondCumThirdEntry}
\end{gather*}
where the contours for the integrals are closed and contain inside all the singularities of the integrand except those at $z_{\sigma_2(2j)}=\pm 1$ and $z_{\sigma_2(2j)}\tilde{z}_{\sigma_2(2j-1)}=1$. Not that these conditions will be satisfied for any two positively oriented closed contours such that $|z_{\sigma_2(2j)}|<1$ and $|z_{\sigma_2(2j)}\tilde{z}_{\sigma_2(2j-)}|<1$. We now make a substitution $\tilde{z}_{\sigma_2(2j-1)}= z^{-1}$ and $z_{\sigma_2(2j)}=w^{-1}$. Rewriting $\mathcal{K}^{(T)}_{21}(\sigma_2(2j-1), \sigma_2(2j))$ in terms of $z$ and $w$, we get
\begin{gather}
\mathcal{K}^{(T)}_{21}(\sigma_2(2j-1), \sigma_2(2j))= -\oint_{\mathfrak{C}_3}\oint_{\mathfrak{C}_4}\frac{{\rm d}z{\rm d}w}{z^{t_{\sigma_2(2j-1)}}w^{t_{\sigma_2(2j)}}} \frac{(z-w)}{z\big(w^2-1\big)(zw-1)} \nonumber\\
\hphantom{\mathcal{K}^{(T)}_{21}(\sigma_2(2j-1), \sigma_2(2j))=}{} \times\prod_{i=1}^{n} \frac{\big(1-z^{-1}x_i\big)\big(1-w^{-1}x_i\big)\big(1-z^{-1}y_i\big)}{(1-zx_i)(1-wx_i)(1-zy_i)}, \label{eq:FormOfTheSecondCumThirdEntryModified}
\end{gather}
where $\mathfrak{C}_3$ and $\mathfrak{C}_6$ are closed positively oriented contours containing the singularities at $\big\{\!x^{-1}_i\!\big\}_{\!i{\in}[1,n]}$. Furthermore, $|z|>1$ and $|zw|>1$ are satisfied for all $(z,w)\in \mathfrak{C}_3\times \mathfrak{C}_4$.

The last case is when both $\sigma(2j-1)$ and $\sigma(2j)$ are even. Substituting the expression of $\mathcal{M}(\sigma(2j-1), \sigma(2j))$ inside the integral of $\mathcal{K}^{(T)}_{22}(\sigma_2(2j-1),\sigma_2(2j))$ yields
\begin{gather}
- \oint \oint \frac{{\rm d}z_{\sigma_2(2j-1)}{\rm d}z_{\sigma_2(2j)}}{z_{\sigma_2(2j-1)}^{-t_{\sigma_2(2j-1)}}z_{\sigma_2(2j)}^{-t_{\sigma_2(2j)}}} \frac{\big(\tilde{z}^{-1}_{\sigma_2(2j-1)}-\tilde{z}^{-1}_{\sigma_2(2j)}\big)}{\tilde{z}_{\sigma_2(2j-1)} \tilde{z}_{\sigma_2(2j)}\big(1-\tilde{z}^{-1}_{\sigma_2(2j-1)}\tilde{z}^{-1}_{\sigma_2(2j)}\big)}\nonumber\\
\qquad{} \times \prod_{i=1}^n\frac{(1-\tilde{z}_{\sigma_2(2j-1)}x_i)(1-\tilde{z}_{\sigma_2(2j)}x_i)}{\big(1-\tilde{z}^{-1}_{\sigma_2(2j-1)}x_i\big) \big(1-\tilde{z}^{-1}_{\sigma_2(2j)}x_i\big)\big(1-\tilde{z}^{-1}_{\sigma_2(2j-1)}y_i\big)\big(1-\tilde{z}^{-1}_{\sigma_2(2j)}y_i\big)}, \label{eq:FormOfTheLastEntry}
\end{gather}
where the contours of the integral must contain inside the singularities at $\{x_i\}_{i\in [1,n]}$, but avoid the pole at $\tilde{z}_{\sigma_2(2j-1)}=\tilde{z}^{-1}_{\sigma_2(2j)}$. These conditions will be fulfilled for any two positively oriented contours which contain inside $\{x_i\}_{i\in [1,n]}$ and further, satisfies $|\tilde{z}_{\sigma_2(2j-1)}|>1$ and $|\tilde{z}_{\sigma_2(2j)}|>1$. Substituting $\tilde{z}_{\sigma_2(2j-1)}=z^{-1}$ and $\tilde{z}_{\sigma_2(2j)}=w^{-1}$ in \eqref{eq:FormOfTheLastEntry}, we arrive at
\begin{gather}
\mathcal{K}^{(T)}_{22}(\sigma_2(2j-1), \sigma_2(2j))= \oint_{\mathfrak{C}_5}\oint_{\mathfrak{C}_6} \frac{{\rm d}z{\rm d}w}{z^{t_{\sigma_2(2j-1)}}w^{t_{\sigma_2(2j)}}}\frac{(z-w)}{z(w^2-1)(zw-1)} \nonumber\\
\hphantom{\mathcal{K}^{(T)}_{22}(\sigma_2(2j-1), \sigma_2(2j))=}{} \times\prod_{i=1}^{n} \frac{\big(1-z^{-1}x_i\big)\big(1-w^{-1}x_i\big)}{(1-zx_i)(1-wx_i)(1-zy_i)(1-qy_i)}, \label{eq:FormOfTheLastEntrymodified}
\end{gather}
where the contours $\mathfrak{C}_5$ and $\mathcal{C}_6$ are closed, positively oriented and contain the poles at $\{x_i\}_{1\in [1,n]}$. Furthermore, for all $(z,w)\in \mathfrak{C}_5\times \mathfrak{C}_6$, we have $|z|<1$ and $|w|<1$.

Now, note that $\sigma$ is an arbitrary permutation in the symmetric group $\mathfrak{S}(2d)$. Therefore, combining \eqref{eq:ExapnsionWithOthertermsInt}, \eqref{eq:FormOfTheFirstEntryModified}, \eqref{eq:FormOfTheSecondCumThirdEntryModified} and \eqref{eq:FormOfTheLastEntrymodified}, one observes that
\begin{gather*}
\rho_{\mathrm{psm}}(T)= \frac{1}{2^d d!}\sum_{\sigma\in \mathfrak{S}(2d)} (-1)^{\sigma}\prod_{j=1}^{d} \mathcal{K}^{(T)}(\sigma(2j-1), \sigma(2j)),
\end{gather*}
where the right hand side is equal to $\Pf(\mathcal{K}^{(T)})$. This completes the proof.
\end{proof}

\section{Correlation function of the Pfaffian Schur process}\label{MainStepOfProof}
The main goal of this section is to prove Theorem~\ref{MainTheorem} which we do in the following subsection. Our proof will require some additional results on the nested inner products of the Schur functions whose proofs are deferred to the subsequent subsections.

\subsection{Proof of Theorem~\ref{MainTheorem}}\label{ProofOfMainTheorem}

Recall that $T_i=\{(i,t_{i, 1}),\dots ,(i,t_{i,d_i})\}$ for $1\leq i\leq m$. Consider $\rho^{+}_{1},\dots ,\rho^{+}_{m}$ and $\rho^{-}_0,\dots ,\rho^{-}_{m-1}$ two sets of finite Schur non-negative specializations of symmetric polynomials such that
\begin{gather}\label{eq:RhoPlusDef}
\rho^{+}_{i}= \{x_{i,1}, x_{i,2}, \dots , x_{i,n_i}\}
\end{gather}
and
\begin{gather}\label{eq:RhoMinusDef}
\rho^{-}_{i-1}=\{y_{i,1}, y_{i,2}, \dots , y_{i,n_i}\}
\end{gather}
for $1\leq i\leq m$. We assume that $\{x_{i,j}\}_{i\in [1,m], j\in [1,d_i]}$ and $\{y_{i,j}\}_{i\in [1,m], j\in [1,d_i]}$ satisfy the following
\begin{gather*}
\min|\rho^{-}_{i-1}|:= \min_{j\in [1,n_i]}\{|x_{i,j}|\}< \max |\rho^{-}_{i-1} | := \max_{j\in [1,n_i]}\{|x_{i,j}|\}<1,\\
\min|\rho^{+}_i|:= \min_{j\in [1,n_i]}\{|x_{i,j}|\}< \max|\rho^{+}_i| := \max_{j\in [1,n_i]}\{|x_{i,j}|\}<1
\end{gather*}
for all $i\in [1,m]$.
 Using analytic continuation, one can prove Theorem~\ref{MainTheorem} for any other finite Schur non-negative specializations which satisfies the assumption of the theorem. Denote $\rho^{+}=\{\rho^{+}_1,\dots ,\rho^{+}_m\}$ and $\rho^{-}=\{\rho^{-}_0, \dots , \rho^{-}_{m-1}\}$. Let us also fix two sets of positive integers $u_1:=\{u_{i,1}\}_{i\in [1,m-1]}$ and $u_2:=\{u_{i,2}\}_{i\in [1,m-1]}$. We will often denote $u_1$ and $u_2$ together as $u$.
For a finite non-negative specialization $X=\{x_1,\dots ,x_n\}$, we introduce the following shorthand notations:\vspace{-0.3cm}
\begin{gather}
H^{(1)}_{q}(X;z) := \prod_{i=1}^{n} \frac{1-zx_i}{1-qzx_i},\label{eq:DefineH1}\\
H^{(2)}_{q}(X;z) := \prod_{i=1}^{n} \frac{1-(qz)^{-1}x_i}{1-z^{-1}x_i}.\label{eq:DefineH2}
\end{gather}
 We consider the following generalization of the r.h.s.\ in \eqref{eq:MacOnPartMultiple}
\begin{gather}
 C(\rho^{+}; \rho^{-}; u):= \sum_{\mathclap{(\bar{\lambda},\bar{\mu})\in \mathbb{Y}^{m}\times \mathbb{Y}^{m-1}}} \mathcal{P}_{psp}(\bar{\lambda}, \bar{\mu}; \rho^{+}, \rho^{-}) \prod_{i=1}^{m-1}\mathbbm{1}(|\mu^{(i)}|\leq u_{i,1}\wedge u_{i,2}) \nonumber\\
 \hphantom{C(\rho^{+}; \rho^{-}; u):=}{} \times \prod_{i=1}^{m-1}\prod_{j=1}^{d_i} \sum_{k=1}^{n_i+u_{i,1}-1} q^{\lambda^{(i)}_k+n_i-k}_{i,j} \prod_{j=1}^{d_m} \sum_{k=1}^{n_m-1} q^{\lambda^{(m)}_k+n_m-k}_{m,j},\label{eq:ActionOnPfaffian}
\end{gather}
 where $\mathcal{P}_{psp}(\bar{\lambda}, \bar{\mu};\rho^{+}, \rho^{-})$ is an abbreviation used to denote the Pfaffian Schur process
\begin{gather*}
Z^{-1}(\rho^{+};\rho^{-})\tau_{\lambda^{(1)}}(\rho_0^{-})s_{\lambda^{(1)}/\mu^{(1)}}(\rho_1^{+})s_{\lambda^{(2)}/\mu^{(1)}}(\rho_1^{-})\cdots s_{\lambda^{(m)}/\mu^{(m-1)}}(\rho^{-}_{m-1})s_{\lambda^{(m)}}(\rho^{+}_{m}).
\end{gather*}
We divide the rest of the argument into four main steps.
 In Step I, we will show that the correlation function $\rho_{\mathrm{psp}}(T)$ is a limit of some contour integral formula involving $C(\rho^{+}; \rho^{-}; u)$. In Step~II, we state a contour integral formula of the limit of $C(\rho^{+}; \rho^{-}; u)$. The proof of this last formula will be deferred to Sections~\ref{sec:NIP}--\ref{sec:CIFofC}. In Steps~II and~III, we show that $\rho_{\mathrm{psp}}(T)$ can be written as the Pfaffian of some skew symmetric matrix of size $\big(\sum_{i=1}^{m}2d_i\big)\times \big(\sum_{i=1}^{m}2d_i\big)$ and identify the corresponding matrix respectively.

{\bf Step I.}
Here, we claim and prove that 
\begin{gather}\label{eq:CorrelationFunction1}
\rho_{psp}(T)= \lim_{u\to \infty}\frac{1}{(2\pi {\rm i})^{\sum_{i=1}^m d_i}}\oint_{\mathcal{C}^\prime} \cdots \oint_{\mathcal{C}^\prime} C(\rho^{+}; \rho^{-}; u) \prod_{i=1}^m \prod_{j=1}^{d_i} q^{-t_{i,j}-n_i-1}_{i,j} {\rm d}q_{i,j},
\end{gather}
where the contours of $q_{i,j}$ are anti-clockwise oriented circles $|z|=\eta$ for some $\eta<1$ and by $\lim_{u\to \infty}$, we mean
\begin{gather}\label{eq:SuccessiveLimitingNotations}
\lim_{u\to \infty}:=\lim_{u_{m-1,2}\to \infty}\lim_{u_{m-1,1}\to \infty}\cdots \lim_{u_{1,2}\to \infty}\lim_{u_{1,1}\to \infty}.
\end{gather}

\begin{proof}[Proof of \eqref{eq:CorrelationFunction1}] Under the assumption $\max\limits_{i,j}|q_{i,j}|<1$, it can be noted that the series in \eqref{eq:ActionOnPfaffian} is absolutely convergent when $\sum_{\bar{\lambda},\bar{\mu}} \mathcal{P}_{psp}(\bar{\lambda}, \bar{\mu};\rho^{+}, \rho^{-})$ converges absolutely. Moreover, we observe that the coefficient of $\prod_{i=1}^m \prod_{j=1}^{d_i}q^{t_{i,j}+n_i}_{i,j}$ in $C(\rho^{+}; \rho^{-}; u)$ is given by
\begin{gather}
\sum_{\mathclap{(\bar{\lambda},\bar{\mu}) \in \mathbb{Y}^m\times \mathbb{Y}^{m-1}}} \mathcal{P}_{psp}(\bar{\lambda}, \bar{\mu}; \rho^{+}, \rho^{-}) \prod_{i=1}^{m-1}\mathbbm{1}\big(\big|\mu^{(i)}\big|\leq u_{i,1}\wedge u_{i,2}\big) \prod_{i=1}^{m-1} \prod_{j=1}^{d_i} \mathbbm{1}\big(t_{i,j}\in \big\{\lambda^{(i)}_k-k\big\}_{1\leq k\leq n_i+u_{i,1}-1}\big)\nonumber\\
\qquad{} \times\prod_{j=1}^{d_m} \mathbbm{1}\big(t_{m,j}\in \big\{\lambda^{(m)}_k-k\big\}_{1\leq k\leq n_i-1}\big).\label{eq:RestrictedCoeffcient}
\end{gather}
Using the absolute convergence of the series $\sum_{\bar{\lambda},\bar{\mu}} \mathcal{P}_{psp}(\bar{\lambda}, \bar{\mu};\rho^{+}, \rho^{-})$, we write~\eqref{eq:RestrictedCoeffcient} as
\begin{gather*}
\left(\frac{1}{(2\pi{\rm i})}\right)^{\sum_{i=1}^m d_i} \oint_{\mathcal{C}^\prime} \cdots \oint_{\mathcal{C}^\prime}C(\rho^{+}; \rho^{-}; u) \prod_{i=1}^{m} \prod_{j=1}^{d_i} q^{-t_{i,j}-n_i-1}_{i,j} {\rm d}q_{i,j},
\end{gather*}
where the contour $\mathcal{C}^{\prime}$ is a circle of radius $\eta$ centered around $0$. Notice that the correlation function $\rho_{psp}(T)$ of the Pfaffian Schur process is the limit of \eqref{eq:RestrictedCoeffcient} when $u_{i,1}\wedge u_{i,2}$ tend to $\infty$ for all $i\in [1,m-1]$. This proves the claim.
\end{proof}

{\bf Step II.} Here, we will derive a contour integral formula of $\lim\limits_{u\to \infty} C(\rho^{+};\rho^{-};u)$.
\bt\label{ContourIntegralRepOfC}
Consider the Pfaffian Schur measure given in \eqref{eq:DefPfaffianSchurMeasure} where $\rho^{+}=\{\rho^{+}_1,\dots ,\rho^{+}_m\}$ and $\rho^{-}=\{\rho^{-}_0,\dots ,\rho^{-}_{m-1}\}$ are two sets of finite Schur non-negative specialization of the plane partition $($see~\eqref{eq:RhoPlusDef}--\eqref{eq:RhoMinusDef}$)$. Let $\{q_{i,j}\}_{j\in [1,d_i], i\in [1,m]}$ be a sequence of complex numbers such that
\begin{gather}\label{eq:ConstrainONQ}
0<\min_{i,j}\{|q_{i,j}|\} \leq \max_{i,j}\{|q_{i,j}|\} <1
\end{gather}
holds. Let $\mathcal{C}$ be a contour which is composed of a clockwise circle $|z|=r$ and a counterclockwise circle $|z|=R$ where $r,R$ satisfy
 \begin{gather}\label{eq:ConstraintOnRandr}
 r<\min_{i\in [1,m]} \{\min|\rho^{+}_i|\} \leq\max_{i\in [1,m]} \{\max|\rho^{+}_i|\} < R <1.
 \end{gather}
 Then, we have
 \begin{gather}
 \lim_{u\to \infty} C(\rho^{+};\rho^{-}; u)= \frac{Z(\rho^{+};\rho^{-})}{(2\pi{\rm i})^{\sum_{i=1}^m d_i}}\oint_{\mathcal{C}}\cdots \oint_{\mathcal{C}} \prod_{i=1}^{m} \prod_{j=1}^{d_i} \frac{q^{n_i}_{i,j}}{(q_jz_j-z_j)}\nonumber\\
 \hphantom{\lim_{u\to \infty} C(\rho^{+};\rho^{-}; u)=}{}\times \prod_{i=1}^m\prod_{ j<k}^{d_i}\frac{(q_{i,k}z_{i,k}-q_{i,j}z_{i,j})(z_{i,j}-z_{i,k})}{(z_{i,j}-q_{i,k}z_{i,k})(q_{i,j}z_{i,j}-z_{i,k})}\nonumber\\
 \hphantom{\lim_{u\to \infty} C(\rho^{+};\rho^{-}; u)=}{} \times \prod_{i=1}^m\left[\prod_{1\leq j< k \leq d_i}\frac{(1-q_{i,j}z_{i,j}z_{i,k})(1-q_{i,k}z_{i,k}z_{i,j})}{(1-q_{i,k}q_{i,j}z_{i,k}z_{i,j})(1-z_{i,j}z_{i,k})}
 \prod_{j=1}^{d_i}\frac{1-q_{i,j}z^2_{i,j}}{1-z^2_{i,j}}\right]\nonumber\\
\hphantom{\lim_{u\to \infty} C(\rho^{+};\rho^{-}; u)=}{} \times \prod_{k=1}^{m} \prod_{i=1}^{m}\prod_{j=1}^{d_{k}}H^{(1)}_{q_{k,j}}(\rho^{+}_{i};z_{k,j})\nonumber\\
\hphantom{\lim_{u\to \infty} C(\rho^{+};\rho^{-}; u)=}{} \times \prod_{1\leq m_1< m_2\leq m}\prod_{j=1}^{d_{m_1}}\prod_{k=1}^{d_{m_2}}\frac{(1-q_{m_2,j}z_{m_2,k}z_{m_1,j})(1-q_{m_1,j}z_{m_1,j}z_{m_2,k})}{(1-q_{m_1,j}z_{m_1,j}q_{m_2,k}z_{m_2,k})(1-z_{m_2,k}z_{m_1,j})}\nonumber\\
\hphantom{\lim_{u\to \infty} C(\rho^{+};\rho^{-}; u)=}{} \times \prod_{k=1}^{m} \prod_{i=0}^{k-1}\prod_{j=1}^{d_{k}}H^{(1)}_{q_{k,j}}(\rho^{-}_{i};z_{k,j})\nonumber\\
\hphantom{\lim_{u\to \infty} C(\rho^{+};\rho^{-}; u)=}{} \times \prod_{1\leq m_1<m_2\leq m}\prod_{j=1}^{d_{m_1}}\prod_{k=1}^{d_{m_2}}\frac{(q_{m_1,j}z_{m_1,j}-q_{m_2,j}z_{m_2,k})(z_{m_1,j}-z_{m_2,k})}{(z_{m_1,j}-q_{m_2,k}z_{m_2,k})(q_{m_1,j}z_{m_1,j}-z_{m_2,k})}
\nonumber\\
\hphantom{\lim_{u\to \infty} C(\rho^{+};\rho^{-}; u)=}{} \times \prod_{k=1}^{m}\prod_{i=k}^{m}\prod_{j=1}^{d_{k}}H^{(2)}_{q_{k,j}}(\rho^{+}_{i};z_{k,j}) \prod_{i=1}^{m} \prod_{j=1}^{d_i} {\rm d}z_{i,j},\label{eq:ContourIntegral12}
 \end{gather}
 where $\lim\limits_{u\to \infty}$ is defined in the same way as in~\eqref{eq:SuccessiveLimitingNotations}.
\et

The proof of Theorem~\ref{ContourIntegralRepOfC} is rather long and will span over three subsections. In the first of these three subsections (Section~\ref{sec:NIP}), we represent $C(\rho^{+};\rho^{-};u)$ as a nested inner products of Schur functions. In the second subsection (Section~\ref{sec:LNIP}), we will compute the limit of the nested inner products as $u\to \infty$ and the final subsection (Section~\ref{sec:CIFofC}) will contain the derivation of a~contour integral formula of the corresponding limit. In the following two steps, we assume the result of Theorem~\ref{ContourIntegralRepOfC} to complete the proof of Theorem~\ref{MainTheorem}.

 {\bf Step III.} Here, we claim and prove that
\begin{gather}
\rho_{psp}(T)= \frac{1}{2^d d!}\sum_{\sigma\in \mathfrak{S}(2d)}(-1)^{\sigma} \mathcal{K}^{(T)}_{\tilde{\sigma}(1),\tilde{\sigma}(2)}\big(\sigma(1)_1,\sigma(1)_2;\sigma(2)_1,\sigma(2)_2\big)\times \cdots \nonumber\\
\hphantom{\rho_{psp}(T)=}{} \times \cdots\times\mathcal{K}^{(T)}_{\tilde{\sigma}(2d-1),\tilde{\sigma}(2d)}\big(\sigma(2d-1)_1,\sigma(2d-1)_2;\sigma(2d)_1,\sigma(2d)_2\big)\label{eq:ReArrangement}
\end{gather}
for some skew symmetric matrix $\mathcal{K}^{T}$ where
\begin{gather*}
\sigma(\zeta)_1 := i, \qquad \text{if} \quad 2d_1+\cdots +2d_{i-1}< \sigma(\zeta)\leq 2d_1+\dots +2d_{i},\\
\sigma(\zeta)_2 = \sigma(\zeta) - \sum_{k=1}^{\sigma(\zeta)_1-1} 2d_k.
\end{gather*}
Furthermore, for any $\sigma\in \mathfrak{S}(2d)$ and $i\in [1,2d]$, we have used the notation $\tilde{\sigma}(i)$ in \eqref{eq:ReArrangement} to denote $1$ in the case when $\sigma(i)$ is an odd integer and $2$ when it is an even integer.

\begin{proof}[Proof of \eqref{eq:ReArrangement}]
 For showing \eqref{eq:ReArrangement}, we use the formula in \eqref{eq:CorrelationFunction1}. The limiting value of \allowbreak $C(\rho^{+}; \rho^{-}; u)$ as $u\to \infty$ (in the sense of \eqref{eq:SuccessiveLimitingNotations}) will be computed in Theorem~\ref{ContourIntegralRepOfC} which we use here. In order to use the limiting value of $C(\rho^{+}; \rho^{-}; u)$ given in \eqref{eq:ContourIntegral12}, one need to justify the interchange of the limit and the integral in \eqref{eq:CorrelationFunction1}. To carry out this interchange, one may use the dominated convergence theorem (DCT) and to apply DCT, it suffices to show that the integrand is uniformly bounded. To this end, we choose $\eta$ such that
\begin{gather}\label{eq:ConstraintOnEta}
\max_{i\in [0,m-1]}\{\max |\rho^{-}_{i}|\}<\eta^{dm^2}< 1.
\end{gather}
For the proof of interchange of limit and integral, we refer to \cite[Lemma~3.1.2]{A14}. In the present context, the assumption in \cite[Lemma~3.1.2]{A14} translates to the constraint on $\eta$ mentioned in~\eqref{eq:ConstraintOnEta}.

So, we are now entitled to interchange the limit and the integral in \eqref{eq:CorrelationFunction1}. After substituting $\lim\limits_{u\to \infty} C(\rho^{+};\rho^{-};u)$ from \eqref{eq:ContourIntegral12} inside the integral of~\eqref{eq:SuccessiveLimitingNotations}, we simplify the integrand using the Pfaffian representation of~\eqref{eq:PrinciplePfaffian} (using~\eqref{eq:PfaffianIdentityToShow}) in the proof of Theorem~\ref{TheoremForSinglePartition}. As a result, we get
\begin{gather}
\rho_{psp}(T)= \frac{1}{(2\pi {\rm i})^{2d}}\oint \cdots \oint \mathrm{Pf}\left(\mathcal{M}_{m}\right) \prod_{k=1}^{m} \prod_{i=1}^{m}\prod_{j=1}^{d_{k}}H^{(1)}_{q_{k,j}}(\rho^{+}_{i};z_{k,j})\prod_{k=1}^{m} \prod_{i=0}^{k-1}\prod_{j=1}^{d_{k}}H^{(1)}_{q_{k,j}}(\rho^{-}_{i};z_{k,j})\nonumber\\
\hphantom{\rho_{psp}(T)=}{} \times\prod_{k=1}^{m}\prod_{j=1}^{d_k}\frac{q^{-t_{k,j}-1}_{k,j}}{z^2_{k,j}-1}\prod_{k=1}^{m}\prod_{i=k}^{m}\prod_{j=1}^{d_{k}}H^{(2)}_{q_{k,j}}(\rho^{+}_{i};z_{k,j})\prod_{i=1}^{m} \prod_{j=1}^{d_i} {\rm d}z_{i,j}\prod_{i=1}^{m} \prod_{j=1}^{d_i} {\rm d}q_{i,j},\label{eq:FinalPfaffianRep}
\end{gather}
where $d=\sum_{i=1}^m d_i$ and $\mathcal{M}_m$ is an $2d\times 2d$ skew-symmetric matrix composed of the following $2\times 2$ block matrix
\begin{gather}
\mathcal{M}_{m}(i,u;j,v):= \begin{pmatrix}
\dfrac{z_{i,u}-z_{j,v}}{1-z_{i,u}z_{j,v}} & \dfrac{1-q_{j,v}z_{j,v}z_{i,u}}{z_{i,u}-q_{j,v}z_{j,v}}\vspace{1mm}\\
-\dfrac{1-q_{i,u}z_{i,u}z_{j,v}}{z_{j,v}-q_{i,u}z_{i,u}} & \dfrac{q_{i,u}z_{i,u}-q_{j,v}z_{j,v}}{1-q_{i,u}q_{j,v}z_{i,u}z_{j,v}}
\end{pmatrix}\nonumber\\
\hphantom{\mathcal{M}_{m}(i,u;j,v)}{} =\begin{pmatrix}
\mathcal{M}^{11}_{m}(i,u;j,v) & \mathcal{M}^{12}_{m}(i,u;j,v)\\
\mathcal{M}^{21}_{m}(i,u;j,v) & \mathcal{M}^{22}_{m}(i,u;j,v)
\end{pmatrix},\label{eq:PfaffianMatrixDescription}
\end{gather}
where the matrix on the far right hand side of \eqref{eq:PfaffianMatrixDescription} is written to introduce our shorthand notations for the elements of $\mathcal{M}_m$. The number of integrals in \eqref{eq:PfaffianMatrixDescription} is equal to $2\sum_{k=1}^{m}d_i$. The contour for the integral with respect to $q_{k,j}$ is the circle $|z|=\eta$ where $\eta$ satisfies~\eqref{eq:ConstraintOnEta}. The contour for the integral with respect to $z_{k,j}$ are composed of two circles $|z|=r$ and $|z|=R$ where~$r$,~$R$ satisfies \eqref{eq:ConstraintOnRandr}. Due to~\eqref{eq:ConstraintOnEta} and~\eqref{eq:ConstraintOnRandr}, there is no pole of the integrand for $z_{k,j}$ inside the circle $|z|=r$ except a simple pole at $0$. Since the residue of the integrand with respect~$z_{i,j}$ at~$0$ has a factor of $q_{k,j}^{-t_{i,j}-1}$, the integral of the residue with respect to $q_{k,j}$ will be equal to~$0$. Thus, we can restrict the contours of the variables $\{z_{k,j}\colon 1\leq k\leq m,\, 1\leq j\leq d_k\}$ in~\eqref{eq:PfaffianMatrixDescription} to the circle $|z|=R$.

Let us now elaborate bit more on the structure of $\mathcal{M}_m$. The matrix $\mathcal{M}_m$ can be partitioned into $m^2$ submatrices. The size of the $(i,j)$ submatrix of $\mathcal{M}_m$ is given by $2d_i\times 2d_j$. In fact, the $(u,v)$ block matrix (of order $2\times 2$) of the $(i,j)$ submatrix is given by $\mathcal{M}_m(i,u; j,v)$.

We expand $\Pf(\mathcal{M}_m)$ inside the integral in \eqref{eq:FinalPfaffianRep} as
\begin{gather}
\frac{1}{2^d d!}\sum_{\sigma\in \mathfrak{S}(2d)}(-1)^{\sigma}\mathcal{M}^{\tilde{\sigma}(1),\tilde{\sigma}(2)}_m\big(\sigma(1)_1,\sigma(1)_2;\sigma(2)_1,\sigma(2)_2\big)\times \cdots\nonumber\\
\qquad {} \times \cdots \times\mathcal{M}^{\tilde{\sigma}(2d-1),\tilde{\sigma}(2d)}_m\big(\sigma(2d-1)_1,\sigma(2d-1)_2;\sigma(2d)_1,\sigma(2d)_2\big).
\label{eq:PfaffianExpansionOfM}
\end{gather}
Next, we do a transformation of the variables $(\{z_{i,u}\}_{i,u}, \{q_{i,u}\}_{i,u})\mapsto(\{z_{i,u}\}_{i,u}, \{\tilde{z}_{i,u}\}_{i,u})$ by de\-fi\-ning $\tilde{z}_{i,k}=(q_{i,u}z_{i,u})^{-1}$ for $i\in [1,m]$ and $u\in [1,d_i]$ inside the integral \eqref{eq:FinalPfaffianRep}. By incorpora\-ting~$H^{(1)}_{q_{k,j}}$ and $H^{(2)}_{q_{k,j}}$ with the elements of the matrix $\mathcal{M}_{m}$ and taking the sum \eqref{eq:PfaffianExpansionOfM} corresponding to $\mathrm{Pf}(\mathcal{M}_{m})$ outside the integral of~\eqref{eq:FinalPfaffianRep} we arrive at~\eqref{eq:ReArrangement}.
\end{proof}

{\bf Step IV.}
 Now, it remains to identify the elements of the matrix $\mathcal{K}^{(T)}$. To achieve this, we fix a permutation $\sigma\in \mathfrak{S}(2d)$. We divide the following discussions into three main sub-steps which are given as follows.

(1) First we consider the case when $\sigma(2\zeta-1)$ and $\sigma(2\zeta)$ both are odd integers for some $\zeta\in [1,d]$. Then the matrix element $\mathcal{K}^{(T)}_{\tilde{\sigma}(2\zeta-1),\tilde{\sigma}(2\zeta)}\left(\sigma(2\zeta-1)_1,\sigma(2\zeta-1)_2;\sigma(2\zeta)_1,\sigma(2\zeta)_2\right)$ is given by $\mathcal{K}^{(T)}_{1,1}(i,u;j,v)$
where
\begin{gather}\label{eq:PermutationIdentification}
\sigma(2\zeta-1)_1=i ,\qquad \sigma(2\zeta-1)_2=u, \qquad \sigma(2\zeta)_1 = j, \qquad \sigma(2\zeta)_2=v.
\end{gather}
We retrieve $\mathcal{K}^{(T)}_{1,1}(i,u;j,v)$ by combining all the terms related to the variables $z_{i,u}$ and $z_{j,v}$ inside the integral \eqref{eq:FinalPfaffianRep} with $\mathcal{M}^{11}_m (i,u;j,v )$. By doing so, we arrive at
\begin{gather}
\mathcal{K}^{(T)}_{1,1} (i,u;j,v ):= - \oint \oint \frac{{\rm d}z_{i,u} {\rm d}z_{j,v}} {z_{i,u}^{-t_{i,u}}z_{j,v}^{-t_{j,v}}}\frac{(z_{i,u}-z_{j,v})}{\big(1-z_{i,u}^2\big)\big(1-z_{i,v}^2\big)(1-z_{i,u}z_{j,v})}\nonumber\\
\hphantom{\mathcal{K}^{(T)}_{1,1} (i,u;j,v ):=}{} \times \frac{H\big(\rho^{+}_{[i,m]};z^{-1}_{i,u}\big)H\big(\rho^{+}_{[j,m]};w^{-1}_{j,v}\big)}{H(\rho^{+}_{[1,m]}\cup \rho^{-}_{[0,i)};z_{i,u})H(\rho^{+}_{[1,m]}\cup \rho^{-}_{[0,j)};z_{j,v})} , \label{eq:K11AtLast}
\end{gather}
where the contours are positively oriented circles $|z|=R$ for $R<1$. Substituting $z=z^{-1}_{i,u}$ and $w=z^{-1}_{j,v}$, we get~\eqref{eq:PfaffianMatDes1stTerm}.

(2) When $\sigma(2\zeta-1)$ is an odd (or even) integer and $\sigma(2\zeta)$ is an even (or odd) integer, then $\mathcal{K}^{(T)}_{\tilde{\sigma}(2\zeta-1),\tilde{\sigma}(2\zeta)}(i,u;j,v)$ is given by $\mathcal{K}^{(T)}_{1,2}(i,u;j,v)$ (or $\mathcal{K}^{(T)}_{2,1}(i,u;j,v)$) where the choice of~$i$,~$j$,~$u$,~$v$ are same as in~\eqref{eq:PermutationIdentification}. To retrieve $\mathcal{K}^{(T)}_{1,2}(i,k;j,l)$, we combine all the terms related to the variables~$z_{i,u}$ and $\tilde{z}_{j,v}$ with $\mathcal{M}^{12}_m(i,u;j,v)$. Henceforth, we get
 \begin{gather}
\mathcal{K}^{(T)}_{1,2}(i,u;j,v):= \oint \oint \frac{{\rm d}z_{i,u}{\rm d}\tilde{z}_{j,v}}{z_{i,u}^{-t_{i,u}}\tilde{z}_{j,v}^{-t_{j,v}}} \frac{\big(1-z_{i,u}\tilde{z}^{-1}_{j,v}\big)}{\tilde{z}_{j,v}\big(1-z_{i,u}^2\big)\big(z_{i,u}-\tilde{z}^{-1}_{j,v}\big)}\nonumber\\
\hphantom{\mathcal{K}^{(T)}_{1,2}(i,u;j,v):=}{}
\times \frac{H\big(\rho^{+}_{[1,m]}\cup \rho^{-}_{[0,i)};\tilde{z}^{-1}_{j,v}\big)H\big(\rho^{+}_{[j,m]};z^{-1}_{i,u}\big)}{H(\rho^{+}_{[i,m]};\tilde{z}_{j,v})H(\rho^{+}_{[1,m]}\cup \rho^{-}_{[0,j)};z_{i,u})}, \label{eq:K21AtLast}
\end{gather}
where the contour of the variable $z_{i,u}$ is the positively oriented circle $|z|=R<1$. But, there is a dichotomy in the choice of the contours for the variable $\tilde{z}_{j,v}$. We know $\tilde{z}_{j,v}= (q_{j,v}z_{j,v})^{-1}$. To this point, let us recall from Theorem~\ref{ContourIntegralRepOfC} that the contour of $z_{i,u}$ contains $q_{j,v}z_{j,v}$ as a pole when $j>i$. Thus, the contour of $\tilde{z}_{j,v}$ is such that $(a)$ $|z_{i,k}\tilde{z}_{j,v}|<1$ is satisfied when $i\leq j$ and $(b)$ $|z_{i,u}\tilde{z}_{j,v}|>1$ holds when $i>j$. Now, we get \eqref{eq:PfaffianMatDes2ndTerm} after substituting $z=z^{-1}_{i,u}$ and $w=\tilde{z}^{-1}_{j,v}$.

(3) In the case when $\sigma(2\zeta-1)$ and $\sigma(2\zeta)$ both are even integers, then $\mathcal{K}^{(T)}_{\tilde{\sigma}(2\zeta-1),\tilde{\sigma}(2\zeta)}(i,u;j,v)$ is given by $\mathcal{K}^{(T)}_{2,2}(i,u;j,v)$ (see \eqref{eq:PermutationIdentification} for $i,j,u,v$) which can be obtained essentially in the same way from $\mathcal{M}^{22}_m(i,u;j,v)$. Thus, we get
\begin{gather}
\mathcal{K}^{(T)}_{2,2}(i,u;j,v):= -\oint \oint \frac{{\rm d}\tilde{z}_{i,u}{\rm d}\tilde{z}_{j,v}}{\tilde{z}_{i,u}^{-t_{i,u}}\tilde{z}_{j,v}^{-t_{j,v}}}\frac{\big(\tilde{z}^{-1}_{i,u}-\tilde{z}^{-1}_{j,v}\big)} {\tilde{z}_{i,u}\tilde{z}_{j,v}\big(1-\tilde{z}^{-1}_{i,u}\tilde{z}^{-1}_{j,v}\big)}\nonumber\\
\hphantom{\mathcal{K}^{(T)}_{2,2}(i,u;j,v):=}{} \times \frac{H\big(\rho^{+}_{[1,m]}\cup \rho^{-}_{[0,i)};\tilde{z}^{-1}_{i,u}\big)H\big(\rho^{+}_{[1,m]}\cup \rho^{-}_{[0,j)};\tilde{z}^{-1}_{j,v}\big)}{H(\rho^{+}_{[i,m]};\tilde{z}_{i,u})H(\rho^{+}_{[j,m]};\tilde{z}_{j,v})}, \label{eq:K22AtLast}
\end{gather}
where the contours are positively oriented circles around $0$ of radii greater than $1$. Now, under the substitution $z=\tilde{z}^{-1}_{i,u}$ and $w=\tilde{z}^{-1}_{j,v}$, \eqref{eq:K22AtLast} translates to \eqref{eq:PfaffianMatDes3rdTerm}.

Combining \eqref{eq:ReArrangement} with \eqref{eq:K11AtLast}, \eqref{eq:K21AtLast} and \eqref{eq:K22AtLast} completes the proof of Theorem~\ref{MainTheorem}.

The remainder of this section is devoted to the proof of Theorem~\ref{ContourIntegralRepOfC}.

\subsection{Nested inner products}\label{sec:NIP}
The key component in the proof of Theorem~\ref{MainTheorem} was indeed Theorem~\ref{ContourIntegralRepOfC} which gives a contour integral representation of the limit of $\mathcal{C}(\rho^{+};\rho^{-};u)$ as $u$ tends to $\infty$. In the case of the Pfaffian Schur measure, we have used Macdonald difference operators for obtaining similar contour integral formulas (see Theorem~\ref{OneMainTheorem}). But, to use the machinary of the difference operators, we need Schur functions in the place of skew Schur functions. In order to resolve this difficulty, we use the rule of finite inner product given in \eqref{eq:BilinearFormSkewSpecial}. In the following result, we represent the product $\mathcal{P}_{psp}(\bar{\lambda},\bar{\mu}; \rho^{+}, \rho^{-})\prod_{i=1}^{m-1}\mathbbm{1}\big(\big|\mu^{(i)}\big|\leq u_{i,1}\wedge u_{i,2}\big)$ in terms of the nested inner products of the Schur functions.

\bp Consider the following sets of Schur non-negative specializations $\xi^{+}_1, \dots , \xi^{+}_{m-1}$ and $\xi^{-}_1, \dots ,\xi^{-}_{m-1}$. Each of the sets contains countably many formal variables. Let $\{u_{i,1}\}_{i\in [1,m-1]}$ and $\{u_{i,2}\}_{i\in [1,m-1]}$ be two subsets of $\NN$. For any $1\leq i\leq (m-1)$, we denote the first $u$ elements of the set $\xi^{+}_i$ and $\xi^{-}_i$ by $\xi^{+,[1,u]}_{i}$ and $\xi^{-,[1,u]}_{i}$ respectively. Recall the definition of the bilinear form $\langle \cdot,\cdot\rangle$ from~\eqref{eq:BilinearFormSkewSpecial}. Then, one has
\begin{gather}
\mathcal{P}_{psp}(\bar{\lambda},\bar{\mu};\rho^{+}, \rho^{-}) \prod_{i=1}^{m-1}\mathbbm{1}\big(\big|\mu^{(i)}\big|\leq u_{i,1}\wedge u_{i,2}\big)\nonumber\\
\qquad {}= \tau_{\lambda^{(1)}}(\rho^{-}_0)\prod_{i=2}^{m} \big\langle s_{\lambda^{(i)}}\big(\rho^{-}_{i-1},\xi^{+,[1,u_{i-1,2}]}_{i-1}\big), s_{\mu^{(i-1)}}\big(\xi^{+,[1,u_{i-1,2}]}_{i-1}\big)\big\rangle^{\xi^{+,[1,u_{i-1,2}]}_{i-1}}_{u_{i-1,2}}\nonumber\\
 \qquad\quad{} \times \prod_{i=1}^{m-1} \big\langle s_{\lambda^{(i)}}\big(\rho^{+}_i, \xi^{-,[1,u_{i,1}]}_i\big), s_{\mu^{(i)}}\big( \xi^{-,[1,u_{i,1}]}_i\big)\big \rangle^{\xi^{-,[1,u_{i,1}]}_i}_{u_{i,1}} s_{\lambda^{(m)}}(\rho^{+}_m). \label{eq:ScalerProdRep}
\end{gather}
\ep
\begin{proof}
 We will show that the r.h.s.\ of \eqref{eq:ScalerProdRep} is equal to $\mathcal{P}_{psp}(\bar{\lambda},\bar{\mu};\rho^{+}, \rho^{-}) \prod_{i=1}^{m-1}\mathbbm{1}\big(\big|\mu^{(i)}\big|\leq u_{i,1}\wedge u_{i,2}\big)$. To this end, by using $\langle s_{\lambda}(X,Z^{[1,u]}),s_\mu(Z^{[1,u]})\rangle^{Z^{[1,u]}}_{u} = s_{\lambda/\mu}(X)\mathbbm{1}(|\mu|\leq u)$, we simplify the inner products in the r.h.s.\ of~\eqref{eq:ScalerProdRep}. Now, the result follows by writing
\begin{gather*} \mathbbm{1}\big(\big|\mu^{(i)}\big|\leq u_{i,1}\wedge u_{i,2}\big)=\mathbbm{1}\big(\big|\mu^{(i)}\big|\leq u_{i,1}\big) \mathbbm{1}\big(\big|\mu^{(i)}\big|\leq u_{i,2}\big).\tag*{\qed} \end{gather*}\renewcommand{\qed}{}
\end{proof}

Now, we use the representation of \eqref{eq:ScalerProdRep} to write $C(\rho^{+};\rho^{-};u)$ as a nested inner products of partition functions of the Schur and the Pfaffian Schur process.
\bp\label{eq:NestedInnerProductSummation}
 Consider the series $C(\rho^{+};\rho^{-}; u)$ defined in~\eqref{eq:ActionOnPfaffian} where $u$ denotes following two subsets $\{u_{i,1}\}_{i\in [1,m-1]}$, $\{u_{i,2}\}_{i\in [1,m-1]}$ from $\ZZ_{>0}$. Recall the definition of~$Z$ and~$F$ from~\eqref{eq:PfaffianSchurPartitionIdentity} and~\eqref{eq:SchurPartitionIdentity} respectively. Then, we have
\begin{gather}
 C^\prime (\rho^{+}; \rho^{-};u) = \bigg\langle \bigg\langle \cdots\bigg\langle \bigg\langle \prod_{j=1}^{d_1} D^{1, q_{1,j}}_{(n_1+u_{1,1}),[\rho^{+}_{1},\xi^{-,[1,u_{1,1}]}_1]} Z\big(\rho^{+}_{1},\xi^{-,[1,u_{1,1}]}_1; \rho^{-}_{0}\big),\nonumber\\
 \hphantom{C^\prime (\rho^{+}; \rho^{-};u) =}{} F\big(\xi^{-,[1,u_{1,1}]}_{1};\xi^{+,[1,u_{1,2}]}_{1}\big)\bigg\rangle^{\xi^{-,[1,u_{1,1}]}_1}_{u_{1,1}}, \prod_{j=1}^{d_{2}}D^{1,q_{2,j}}_{(n_2+u_{2,1}),[\rho^{+}_{2},\xi^{-,[1,u_{2,1}]}_{2}]}\nonumber\\
\hphantom{C^\prime (\rho^{+}; \rho^{-};u) =}{} F\big(\rho^{+}_{2},\xi^{-,[1,u_{2,1}]}_{2}; \rho^{-}_{1},\xi^{+,[1,u_{2,1}]}_{1}\big) \bigg\rangle^{\xi^{+,[1,u_{2,1}]}_{1}}_{u_{2,1}}\cdots \bigg\rangle^{\xi^{-,[1,u_{m-1,1}]}_{m-1}}_{u_{m-1,1}}, \nonumber\\
\hphantom{C^\prime (\rho^{+}; \rho^{-};u) =}{} \prod_{j=1}^{d_{m}}D^{1,q_{m,j}}_{n,\rho^{+}_{m}}F(\rho^{+}_{m};\rho^{-}_{m-1}, \xi^{+,[1,u_{m-1,2}]}_{m-1})\bigg\rangle^{\xi^{+,[1,u_{m-1,2}]}_{m-1}}_{u_{m-1,2}},\label{eq:ProofOfRep}
\end{gather}
where
\begin{gather*}\label{eq:DefineCPrime}
C^\prime(\rho^{+}; \rho^{-};u):= C(\rho^{+}; \rho^{-};u)\prod_{i=1}^{m-1} \prod_{j=1}^{d_i}q^{u_{i,1}}_{i,j}.
\end{gather*}
Here, $D^{1,q}_{n,X}$ denotes the difference operator which acts on the space ${\rm Sym}$ specialized at~$X$, where~$X$ is of length~$n$.
\ep
\begin{proof}
First, we show how the innermost inner product in \eqref{eq:ProofOfRep} simplifies. We start with recalling the action of $\prod_{j=1}^{d_1}D^{1,q_{1,j}}_{(n_1+u_{1,1}),[\rho^{+}_{1},\xi^{-,[1,u_{1,1}]}_1]}$ on $Z\big(\rho^{+}_{1},\xi^{-,[1,u_{1,1}]}_1; \rho^{-}_{0}\big)$:
\begin{gather*}
\prod_{j=1}^{d_1}D^{1,q_{1,j}}_{(n_1+u_{1,1}),[\rho^{+}_{1},\xi^{-,[1,u_{1,1}]}_1]} Z\big(\rho^{+}_{1},\xi^{-,[1,u_{1,1}]}_1; \rho^{-}_{0}\big)\\
\qquad{} =\sum_{\lambda^{(1)}} \tau_{\lambda^{(1)}}(\rho^{-}_{0}) s_{\lambda^{(1)}}\big(\rho^{+}_{1},\xi^{-,[1,u_{1,1}]}_1\big)\prod_{j=1}^{d_1}\sum_{k=1}^{(n_1+u_{1,1})-1} q^{\lambda^{(1)}_k+(n_1+u_{1,1})-k}_{1,j} .
\end{gather*}
Expanding $F\big(\xi^{-,[1,u_{1,1}]}_{1};\xi^{+,[1,u_{1,2}]}_{1}\big)$ in terms of Schur functions, applying the above equation and computing the inner products of the Schur function as prescribed in~\eqref{eq:BilinearFormSkewSpecial} yields
\begin{gather}
 \left\langle \prod_{j=1}^{d_1} D^{1, q_{1,j}}_{(n_1+u_{1,1}),[\rho^{+}_{1},\xi^{-,[1,u_{1,1}]}_1]} Z\big(\rho^{+}_{1},\xi^{-,[1,u_{1,1}]}_1; \rho^{-}_{0}\big),F\big(\xi^{-,[1,u_{1,1}]}_{1};\xi^{+,[1,u_{1,2}]}_{1}\big)\right\rangle^{\xi^{-,[1,u_{1,1}]}_1}_{u_{1,1}}\nonumber\\
 \qquad{} = \sum_{\mu^{(1)}\subset \lambda^{(1)}} \tau_{\lambda^{(1)}}(\rho^{-}_{0})\mathbbm{1}\big(\big|\mu^{(1)}\big|\leq u_{1,1}\big) s_{\lambda^{(1)}/\mu^{(1)}}(\rho^{+}_{1})s_{\mu^{(1)}}\big(\xi^{+,[1,u_{1,2}]}_{1}\big)\nonumber\\
\qquad\quad{} \times\prod_{j=1}^{d_{1}}\sum_{k=1}^{(n_1+u_{1,1})-1} q^{\lambda^{(1)}_k +(n_1+u_{1,1})-k}_{1,j}.\label{eq:InnerProdDecoding}
\end{gather}
Now, we focus on computing the next inner product. To do this, we first note
\begin{gather}
 D^{1,q_{2,j}}_{(n_2+u_{2,1}), [\rho^{+}_{2},\xi^{-,[1,u_{2,1}]}_{2}]} F\big(\rho^{+}_{2},\xi^{-,[1,u_{2,1}]}_{2}; \rho^{-}_{1},\xi^{+,[1,u_{1,2}]}_{1}\big) \nonumber\\
 \qquad{} = \sum_{\lambda^{(2)}}s_{\lambda^{(2)}}\big(\rho^{+}_{2},\xi^{-,[1,u_{2,1}]}_{2}\big)s_{\lambda^{(2)}}\big(\rho^{-}_{1},\xi^{+,[1,u_{1,2}]}_{1}\big)
 \prod_{j=1}^{d_{2}}\sum_{k=1}^{(n_2+u_{2,1})-1} q^{\lambda^{(2)}_k+(n_{1}+u_{2,1})-k}_{2,j}.\label{eq:SecondTermDecoding}
\end{gather}
Taking the inner product of the r.h.s.\ of \eqref{eq:InnerProdDecoding} with the r.h.s.\ of \eqref{eq:SecondTermDecoding} yields
\begin{gather*}
\sum_{\lambda^{(1)}\supset\mu^{(1)}\subset \lambda^{(2)}} \tau_{\lambda^{(1)}}(\rho^{-}_{0}) \mathbbm{1}\big(\big|\mu^{(1)}\big|\leq u_{1,1}\wedge u_{1,2}\big) s_{\lambda^{(1)}/\mu^{(1)}}(\rho^{+}_{1})s_{\lambda^{(2)}/\mu^{(1)}}(\rho^{-}_{1})\\
\qquad{} \times s_{\lambda^{(2)}}\big(\rho^{+}_{2},\xi^{-,[1,u_{2,1}]}_{2}\big)\prod_{i=1}^{2}\prod_{j=1}^{d_i}\sum_{k=1}^{(n_i+u_{i,1})-1} q^{\lambda^{(i)}_k +(n_i+u_{i,1})-k}_{i,j}.
\end{gather*}
Using induction, we get the following expression on taking the successive inner products until the last one
\begin{gather*}
{\rm (I)}:=\sum_{(\lambda\times \mu)\in \mathbb{Y}^m\times \mathbb{Y}^{m-1}} \tau_{\lambda^{(1)}}(\rho^{-}_{0}) \prod_{i=1}^{m-1}\mathbbm{1}\big(\big|\mu^{(i)}\big|\leq u_{i,1}\wedge u_{i,2}\big) s_{\lambda^{(1)}/\mu^{(1)}}(\rho^{+}_{1})s_{\lambda^{(2)}/\mu^{(1)}}(\rho^{-}_{1})\\
\hphantom{{\rm (I)}:=}{} \times \cdots \times s_{\mu^{(m-1)}}\big(\xi^{+,[1,u_{m-1,2}]}_{m-1}\big)\prod_{i=1}^{m-1}\prod_{j=1}^{d_{i}}\sum_{k=1}^{(n_i+u_{i,1})-1} q^{\lambda^{(i)}_k +(n_i+u_{i,1})-k}_{i,j}.
\end{gather*}
Now, we focus on computing the last inner product on the r.h.s.\ of \eqref{eq:ProofOfRep}. For this, we require to rewrite $\prod_{j=1}^{d_m} D^{1,q_{1,j}}_{n_m,\rho^{+}_{m}} F\big(\rho^{+}_{m};\rho^{-}_{m-1},\xi^{+,[1,u_{m-1,2}]}_{m-1}\big)$ as \begin{gather*}
{\rm (II)}:=\sum_{\lambda^{(m)}}s_{\lambda^{(m)}}(\rho^{+}_m)s_{\lambda^{(m)}}\big(\rho^{-}_{m-1},\xi^{+,[1,u_{m-1,2}]}_{m-1}\big)\prod_{j=1}^{d_m}\sum_{k=1}^{n_m-1} q^{\lambda^{(1)}_k+n_m-k}_{m,j}.
\end{gather*}
On taking the inner product $\langle {\rm (I)}, {\rm (II)}\rangle_{u_{m-1,2}}$ in the space ${\rm Sym}\big(\xi^{+,[1,u]}_{m-1}\big)$, we get $C(\rho^{+};\rho^{-};u)$ modulo a multiplicative factor of $\prod^{m-1}_{i=1}\prod_{j=1}^{d_i} q^{u_{i,1}}_{i,j}$. This completes the proof.
\end{proof}

\subsection{Limit of the nested inner products}\label{sec:LNIP}
The goal of this subsection is to obtain some preliminary result to compute the limit of \allowbreak $C(\rho^{+}; \rho^{-};u)$ as $u \to \infty$ in the sense of~\eqref{eq:ULimit}. To this end, we use the contour integral formulas of the action of Macdonald difference operators (from~\eqref{eq:MCIntForMultipleAction} and~\eqref{eq:MCIntForMultipleActionOnSchurPartition}) and substitute those into the nested inner product formula of Proposition~\ref{eq:NestedInnerProductSummation}. For computing the limit of the inner products of the contour integral formulas, our first step is to compute the limit of the inner products of the integrand. This last task is the main content of this subsection and will be performed in Theorem~\ref{LimitRep}. The result of Theorem~\ref{LimitRep} will be used to complete the computation of $\lim\limits_{u\to \infty} C(\rho^{+}; \rho^{-};u)$ in the next subsection.

Before stating any result, we recollect some elementary information on the contour integral formula for the action of Macdonald operators on~$F(X;Y)$ or~$Z(X;Y)$ where~$F(X;Y)$ and $Z(X;Y)$ are the partition functions of the Schur and Pfaffian Schur measure. Recall the definitions of~$H^{(1)}_{q}$ and~$H^{(2)}_{q}$ from~\eqref{eq:DefineH1} and~\eqref{eq:DefineH2} respectively. Notice that $H^{(2)}_q(X;z) = H^{(1)}_{q}\big(X;(qz)^{-1}\big)$. Note that the integrand in the integral formula
\eqref{eq:MCIntForMultipleActionOnSchurPartition} is a product of two terms given as
\begin{gather*}
\mathcal{A}_1 := F(X;Y)\prod_{j=1}^{d} q^{n}_j\prod_{j=1}^{d} H^{(1)}_{q_j}(Y;z_j)\prod_{j=1}^{d} H^{(2)}_{q_j}(X;z_j),\\
\mathcal{B}_1 := \prod_{j=1}^{d}\frac{1}{q_jz_j-z_j}\prod_{1\leq j<k\leq d} \frac{(q_kz_k-q_jz_j)(z_k-z_j)}{(z_j-q_kz_k)(q_jz_j-z_k)}.
\end{gather*}
 Similarly, one can check that the inner integrand inside the integral of~\eqref{eq:MCIntForMultipleAction} as the product of
 \begin{gather*}
 \mathcal{A}_{2} :=Z(X;Y)\prod_{j=1}^{d}q^{n}_j\prod_{j=1}^{d} H^{(1)}_{q_j}(X,Y;z_j)\prod_{j=1}^{d} H^{(2)}_{q_j}(X;z_j),\\
 \mathcal{B}_{2} := \prod_{j=1}^{d}\frac{\big(1-q_jz^2_j\big)}{(q_jz_j-z_j)\big(1-z^2_j\big)}\prod_{1\leq j<k\leq d} \frac{(q_kz_k-q_jz_j)(z_k-z_j)}{(z_j-q_kz_k)(q_jz_j-z_k)}\\
 \hphantom{\mathcal{B}_{2} :=}{}\times \prod_{1\leq j< k\leq d} \frac{(1-q_jz_jz_k)(1-q_kz_kz_j)}{(1-q_kq_jz_kz_j)(1-z_kz_j)}.
\end{gather*}
After writing the actions of the difference operators inside the inner products in \eqref{eq:ProofOfRep} in terms of the contour integrals, one can see that taking $u\to \infty$ (in the sense of \eqref{eq:SuccessiveLimitingNotations}) of $C(\rho^{+};\rho^{-};u)$ boils down to first computing the inner products of the terms like $\mathcal{A}_1$ and $\mathcal{A}_2$ for some choices of the specializations and then letting $u\to \infty$. 
The main result of this section which we state as follows formalizes the above mentioned heuristic.

\bt\label{LimitRep}
Consider the set of complex numbers $\{z_{i,j}\}_{1\leq j\leq d_i,1\leq i\leq m}$ such that $|z_{i,j}|\leq \Upsilon$ for positive real number $\Upsilon<1$. Let us denote $u:=\{u_{i,1}\}_{i\in [1,m-1]}\cup \{u_{i,2}\}_{i\in [1,m-1]}$.
Then, the limiting value
\begin{gather}
\lim_{u\to \infty} \bigg\langle \bigg\langle \cdots \bigg\langle \bigg\langle Z\big(\rho^{+}_{1},\xi^{-,[1,u_{1,1}]}_1; \rho^{-}_{0}\big)\prod_{j=1}^{d_1} H^{(1)}_{ q_{1,j}}\big(\rho^{+}_{1},\xi^{-,[1,u_{1,1}]}_1, \rho^{-}_{0}; z_{1,j}\big) \nonumber\\
\qquad{}\times H^{(2)}_{q_{1,j}}(\rho^{+}_{1},\xi^{-,[1,u_{1,1}]}_1; z_{1,j}), F\big(\xi^{-,[1,u_{1,1}]}_{1};\xi^{+,[1,u_{1,2}]}_{1}\big)\bigg\rangle^{\xi^{-,[1,u_{1,1}]}_1}_{u_{1,1}},\nonumber\\
\qquad{}F\big(\rho^{+}_{2},\xi^{-,[1,u_{2,1}]}_{2}; \rho^{-}_{1},\xi^{+,[1,u_{1,2}]}_{1}\big) \prod_{j=1}^{d_{2}}H^{(1)}_{q_{2,j}}\big(\rho^{-}_{1},\xi^{+,[1,u_{1,2}]}_{1};z_{2,j}\big)\nonumber\\
\qquad{} \times\prod_{j=1}^{d_2}H^{(2)}_{q_{2,j}}\big(\rho^{+}_{2},\xi^{-,[1,u_{2,1}]}_{2}; z_{2,j}\big)\bigg\rangle^{\xi^{+,[1,u_{1,2}]}_1}_{u_{1,2}}\cdots ,\nonumber\\
\qquad{} F\big(\xi^{-,[1,u_{m-1,1}]}_{m-1};\xi^{+,[1,u_{m-1,2}]}_{m-1}\big)\bigg\rangle^{\xi^{-,[1,u_{m-1,1}]}_{m-1}}_{u_{m-1,1}},
F\big(\rho^{+}_{m};\rho^{-}_{m-1},\xi^{+,[1,u_{m-1,2}]}_{m-1}\big) \nonumber\\
\qquad{}\times \prod_{j=1}^{d_{m}}H^{(1)}_{q_{m,j}}\big(\rho^{-}_{m-1}, \xi^{+,[1,u_{m-1,2}]}_{m-1};z_{m,j}\big)H^{(2)}_{q_{m,j}}(\rho^{+}_{m};z_{m,j})\bigg\rangle^{\xi^{+,[1,u_{m-1,2}]}_{m-1}}_{u_{m-1,2}}
\label{eq:ULimit}
\end{gather}
is given by
\begin{gather*}
H_0\big(\rho^{+}_{[1,m]}\big) \prod_{j=1}^{m-1}\prod_{k\geq j+1}^{m} F(\rho^{+}_{k}; \rho^{-}_{j}) \prod_{i=1}^m\prod_{j=2}^{d_i}\frac{1-q_{i,j}z^2_{i,j}}{1-z^2_{i,j}}\\
\qquad{}\times \prod_{i=2}^m\prod_{1\leq j<k\leq d_i}\frac{(1-q_{i,j}z_{i,j}z_{i,k})(1-q_{i,k}z_{i,k}z_{i,j})}{(1-q_{i,j}z_{i,j}q_{i,k}z_{i,k})(1-z_{i,j}z_{i,k})}\nonumber\\
\qquad{} \times \prod_{k=1}^{m} \prod_{i=1}^{m}\prod_{j=1}^{d_{k}}H^{(1)}_{q_{k,j}}(\rho^{+}_{i};z_{k,j})\prod_{k=1}^{m} \prod_{i=0}^{k-1}\prod_{j=1}^{d_{k}}H^{(1)}_{q_{k,j}}(\rho^{-}_{i};z_{k,1})\prod_{k=1}^{m}\prod_{i=k}^{m}\prod_{j=1}^{d_{k}}H^{(2)}_{q_{k,j}}(\rho^{+}_{i};z_{k,1})
\nonumber\\
\qquad{} \times \prod_{1\leq m_1<m_2\leq m}\prod_{j=1}^{d_{m_1}}\prod_{k=1}^{d_{m_2}}\frac{(q_{m_1,j}z_{m_1,j}-q_{m_2,j}z_{m_2,k})(z_{m_1,j}-z_{m_2,k})}{(z_{m_1,j}-q_{m_2,k}z_{m_2,k})(q_{m_1,j}z_{m_1,j}-z_{m_2,k})}
\nonumber\\
\qquad{} \times \prod_{1\leq m_1<m_2\leq m}\prod_{j=1}^{d_{m_1}}\prod_{k=1}^{d_{m_2}}\frac{(1-q_{m_1,j}z_{m_1,j}z_{m_2,k})(1-q_{m_2,j}z_{m_2,j}z_{m_1,k})}{(1-q_{m_1,j}z_{m_1,j}q_{m_2,k}z_{m_2,k})(1-z_{m_2,k}z_{m_1,j})}.
\end{gather*}
\et

To prove Theorem~\ref{LimitRep}, we need the following lemma.
\bl\label{InnerProductLemma1}
 For any $t\in \NN$, $X^{[1,t]}$ denotes the first $t$ elements of the set $X$. Let $\{a_l\}_{l\in \NN}$ and $\{b_l\}_{l\in \NN}$ be two sequence of elements in the respective $\ZZ_{\geq 0}$-graded algebras $A$ and $B$ with
 \[\lim_{l\to \infty} {\rm ldeg}(a_l)=\infty, \qquad \lim_{l\to \infty} {\rm ldeg}(b_l)=\infty\]
 such that $\sum_{l}l^{-1}p_l(X)a_l\in \overline{{\rm Sym}(X)}\otimes \overline{A}$ and $\sum_{l}l^{-1}p_l(X)b_l\in \overline{{\rm Sym}(X)}\otimes \overline{B}$ where $\overline{A}$ and $\overline{B}$ denote the topological completions of $A$ and $B$ respectively.
 Consider the following functions in $\overline{{\rm Sym}}\otimes \overline{A}$ and $\overline{{\rm Sym}}\otimes \overline{B}$ respectively
\begin{gather}
\psi\big(X^{[1,t]}\big) :=\exp\left(\sum_{l=1}^\infty \frac{p_l(X^{[1,t]})}{l}a_l\right),\nonumber\\
\phi\big(X^{[1,t]}\big) :=\exp\left(\sum_{l=1}^\infty \frac{p_l(X^{[1,t]})}{l}b_l+\sum_{l=1}^\infty \frac{p^2_{l}(X_{[1,t]})-p_{2l}(X_{[1,t]})}{2l}\right).\label{eq:SymPol}
\end{gather}
 Then,
\begin{gather*}
\lim_{t\to \infty} \big\langle \psi\big(X^{[1,t]}\big), \phi\big(X^{[1,t]}\big)\big\rangle^{X^{[1,t]}}_t =\exp\left(\sum_{l=1}^{\infty} \frac{a_lb_l}{l}+\sum_{l=1}^\infty \frac{a^2_l-a_{2l}}{2l}\right).
\end{gather*}
\el
\begin{proof}Let $R$ and $S$ be two alphabets. We first consider the case when we have $p_l(R)$ instead of $a_l$ in $\psi$ and $p_l(S)$ instead of~$b_l$ in~$\phi$. To accommodate this change, we define
\begin{gather*}
\tilde{\psi}\big(X^{[1,t]}\big) := \exp\left(\sum_{l=1}^{\infty}\frac{p_{l}\big(X^{[1,t]}\big)}{l}p_{l}(R)\right), \\
\tilde{\phi}\big(X^{[1,t]}\big) :=\exp\left(\sum_{l=1}^\infty \frac{p_l\big(X^{[1,t]}\big)}{l}p_l(S)+\sum_{l=1}^\infty \frac{p^2_{l}\big(X^{[1,t]}\big)-p_{2l}\big(X^{[1,t]}\big)}{2l}\right).
\end{gather*}
Then, using \eqref{eq:SpecializedH0H1}, one can write
\begin{gather*}
\tilde{\psi}\big(X^{[1,t]}\big) =\sum_{\lambda\in \mathbb{Y}} s_\lambda\big(X^{[1,t]}\big)s_\lambda(R), \qquad
\tilde{\phi}\big(X^{[1,t]}\big) =\sum_{\lambda\in \mathbb{Y}} \tau_\lambda(S)s_\lambda\big(X^{[1,t]}\big).
\end{gather*}
For any $t\in \NN$, the finite inner product $\big\langle \tilde{\psi}\big(X^{[1,t]}\big), \tilde{\phi}\big(X^{[1,t]}\big)\big\rangle_t$ is equal to $\sum_{|\lambda|\leq t}\tau_\lambda(S)s_\lambda(R) $. Thus, the sequence of the finite inner products $\big\langle \tilde{\psi}\big(X^{[1,t]}\big), \tilde{\phi}\big(X^{[1,t]}\big)\big\rangle_t$ converges inside the space $\overline{{\rm Sym}(R)}\otimes \overline{{\rm Sym}(S)}$ (endowed with the graded topology) to the sum $\sum_{\lambda} \tau_\lambda(S)s_\lambda(R)$, namely, the partition function of the Pfaffian Schur measure. Owing to Proposition~\ref{PfaffianPartitionFunctionIdentity}, the limit of $\big\langle \tilde{\psi}\big(X^{[1,t]}\big), \tilde{\phi}\big(X^{[1,t]}\big)\big\rangle_t$ as $t\to \infty$ will be equal to $H_0(R)H(R;S)$. Now, we proceed to find the limit of $\big\langle \psi\big(X^{[1,t]}\big), \phi\big(X^{[1,t]}\big)\big\rangle_t$. Recall that $\{p_\lambda\}_{\lambda}$ is a basis of ${\rm Sym}$. So, the following maps
\begin{gather*}
\Psi_{A}\colon \ \overline{{\rm Sym}(R)} \to \overline{A}, \qquad \Psi_{A}(p_k(R)) = a_k, \qquad \text{and}\\
 \Psi_{B}\colon \ \overline{{\rm Sym}(S)} \to \overline{B}, \qquad \Psi_{B}(p_k(R)) = b_k
\end{gather*}
define algebra homomorphisms (continuous in the graded topology). Now, we define a bilinear algebra homomorphism $\Psi_{A}\otimes \Psi_B$ in $\overline{{\rm Sym}(R)} \otimes \overline{{\rm Sym}(S)}$ as
\begin{gather*}
\Psi_{A}\otimes \Psi_B \colon \ \overline{{\rm Sym}(R)} \otimes \overline{{\rm Sym}(S)}\mapsto \overline{A} \otimes \overline{B}, \qquad (\Psi_{A}\otimes \Psi_{B}) (p_k(R) p_l(S) ) = a_k b_l.
\end{gather*}
Since the inner product preserves an algebra homomorphism, for any $t\in \NN$,
\begin{gather*}
\big\langle \psi\big(X^{[1,t]}\big), \phi\big(X^{[1,t]}\big)\big\rangle^{X^{[1,t]}}_t = \sum_{|\lambda|\leq t} \Psi_{A}(\tau_{\lambda}(S))\Psi_{B}(s_{\lambda}(R)).
\end{gather*}
Using continuity of the map $\Psi_{A}$ and $\Psi_B$, we get
\begin{gather}\label{eq:LimitOfFinitelyMany}
\lim_{t\to \infty}\big\langle \psi\big(X^{[1,t]}\big), \phi\big(X^{[1,t]}\big)\big\rangle^{X^{[1,t]}}_t =\sum_{\lambda\in \mathbb{Y}} \Psi_{A}(\tau_{\lambda}(S))\Psi_{B}(s_\lambda(R)).
\end{gather}

Owing to the bilinearity of the map $\Psi_{A}\otimes \Psi_{B}$, one can write
\begin{gather}
\sum_{\lambda\in \mathbb{Y}} \Psi_{A}(\tau_{\lambda}(R))\Psi_{B}(s_\lambda(S)) =(\Psi_{A}\otimes \Psi_{B})\left(\sum_{\lambda\in \mathbb{Y}} \tau_{\lambda}(R)s_\lambda(S)\right)\nonumber\\
\hphantom{\sum_{\lambda\in \mathbb{Y}} \Psi_{A}(\tau_{\lambda}(R))\Psi_{B}(s_\lambda(S))}{} =(\Psi_{A}\otimes \Psi_{B})\left(H_0(R)H(R;S)\right).\label{eq:Linearity}
\end{gather}
Thus, the limit of those finite inner products is given by
\begin{gather}
\lim_{t\to \infty}\big\langle \psi\big(X^{[1,t]}\big), \phi\big( X^{[1,t]}\big)\big\rangle^{X^{[1,t]}}_t = (\Psi_{A}\otimes \Psi_{B})(H_0(R)H(R;S))\nonumber\\
 \qquad{} =\exp\left(\sum_{l=1}^{\infty} \frac{\Psi_{A}\otimes \Psi_{B}\big(p_l(R)p_l(S)\big)}{l}+\sum_{l=1}^\infty \frac{\Psi_A\otimes \Psi_B\big(p^2_l(R)-p_{2l}(R)\big)}{2l}\right),\label{eq:InnerProductIdentity}
\end{gather}
 where the equality in first line follows by combining \eqref{eq:LimitOfFinitelyMany} with \eqref{eq:Linearity} and the equality of the second line follows by combining the formula of $H_0(R)$ and $H(R;S)$ given in \eqref{eq:SpecializedH0H1} with the fact that $\Psi$ is an algebra homomorphism. Finally, substituting the image of $p_l(R)p_l(S)$ and $p_{l}(R)$ under the map $\Psi_{A}\otimes \Psi_{B}$ in the last line of~\eqref{eq:InnerProductIdentity} completes the proof.
\end{proof}

\begin{rem} It is worth noting that Lemma~\ref{InnerProductLemma1} is similar to Lemma~2.1.3 of \cite{A14} (see also \cite[Proposition~2.3]{BCGS16}). The main difference between these two results is that unlike Lemma~\ref{InnerProductLemma1}, Lemma~2.1.3 of \cite{A14} finds a formula of the inner product between the functions $\exp\big(\sum_{l=1}^{\infty}p_l(X)a_l\big)$ and $\exp\big(\sum_{l=1}^{\infty}p_l(X)b_l\big)$ where $a_l,b_l$ are arbitrary elements of a $\ZZ_{\geq 0}$-graded algebra.
\end{rem}

 \begin{proof}[Proof of Theorem~\ref{LimitRep}] Here, we demonstrate how to compute the limit of the inner products when $m=2$. This requires computation of two successive inner products and their limits. Furthermore, the computation in $m=2$ case is very similar to computing the two innermost inner products of~\eqref{eq:ULimit} and then, taking $u_{1,1}$ and $u_{1,2}$ to $\infty$. The computation of the successive inner products for $m>2$ case can be also done by taking similar route and thus, will be skipped for avoiding repetitions.
 We now claim and prove that
 \begin{gather}
 \lim_{u_{1,2},u_{1,1}\to \infty} \bigg\langle \bigg\langle Z\big(\rho^{+}_{1},\xi^{-,[1,u_{1,1}]}_1; \rho^{-}_{0}\big)\prod_{j=1}^{d_1} H^{(1)}_{ q_{1,j}}\big(\rho^{+}_{1},\xi^{-,[1,u_{1,1}]}_1, \rho^{-}_{0}; z_{1,j}\big)\nonumber\\
 \qquad{}\times H^{(2)}_{q_{1,j}}\big(\rho^{+}_{1},\xi^{-,[1,u_{1,1}]}_1; z_{1,j}\big),F\big(\xi^{-,[1,u_{1,1}]}_{1};\xi^{+,[1,u_{1,2}]}_{1}\big)\bigg\rangle^{\xi^{-,[1,u_{1,1}]}_1}_{u_{1,1}}, \nonumber \\
 \qquad {} F\big(\rho^{+}_{2}; \rho^{-}_{1},\xi^{+,[1,u_{1,2}]}_{1}\big) \prod_{j=1}^{d_{2}}H^{(1)}_{q_{2,j}}\big(\rho^{-}_{1},\xi^{+,[1,u_{1,2}]}_{1};z_{2,j}\big)\prod_{j=1}^{d_2}H^{(2)}_{q_{2,j}}\big(\rho^{+}_{2}; z_{2,j}\big)\bigg\rangle^{\xi^{+,[1,u_{1,2}]}_1}_{u_{1,2}}\label{eq:ULimitAlt}
 \end{gather}
 is equal to
\begin{gather}
 Z\big(\rho^{+}_1,\rho^{+}_{2},\xi^{-,[1,u_{2,1}]}_{2};\rho^{-}_{0}\big)F\big(\rho^{+}_{2},\xi^{-,[1,u_{2,1}]}_{2}; \rho^{-}_{1}\big)\prod_{j=1}^{d_1}H^{(1)}_{ q_{1,j}}\big(\rho^{+}_{1},\rho^{-}_0,\rho^{+}_{2},\xi^{-,[1,u_{2,1}]}_{2}; z_{1,j}\big)\nonumber\\
 \qquad{} \times \prod_{j=1}^{d_{2}}H^{(1)}_{q_{2,j}}\big(\rho^{+}_{1},\rho^{-}_0,\rho^{-}_1,\rho^{+}_{2}; z_{2,j}\big)\prod_{j=1}^{d_1} H^{(2)}_{ q_{1,j}}\big(\rho^{+}_1,\rho^{+}_2; z_{1,j}\big)\nonumber\\
 \qquad{} \times \prod_{j=1}^{d_{2}} H^{(2)}_{q_{2,j}}\big(\rho^{+}_2; z_{1,j}\big)\prod_{j=1}^{2d_2}\frac{1-q_{2,j}z^2_{2,j}}{1-z^2_{2,j}} \prod_{j=1}^{d_1}\prod_{k=1}^{d_{2}}\frac{(q_{1,j}z_{1,j}-q_{2,j}z_{2,k})(z_{1,j}-z_{2,k})}{(z_{1,j}-q_{2,k}z_{2,k})(q_{1,j}z_{1,j}-z_{2,k})}\label{eq:ContourIntegralAlt}\\
\qquad{} \times \prod_{j=1}^{d_1}\prod_{k=1}^{d_{2}}\frac{(1-q_{1,j}z_{1,j}z_{2,k})(1-q_{2,j}z_{2,j}z_{1,k})}{(1-q_{1,j}z_{1,j}q_{2,k}z_{2,k})(1-z_{1,j}z_{2,k})}\!\prod_{1\leq j<k\leq d_2}\!\frac{(1-q_{2,j}z_{2,j}z_{2,k})(1-q_{2,k}z_{2,k}z_{2,j})}{(1-q_{2,j}z_{2,j}q_{2,k}z_{2,k})(1-z_{2,j}z_{2,k})}.\nonumber
\end{gather}

We first analyze the innermost inner product in \eqref{eq:ULimitAlt}. The key tool that we use is Lemma~\ref{InnerProductLemma1}. We start with writing the partition functions $Z(\cdot;\cdot)$, $F(\cdot;\cdot)$ and $H^{(1)}_q$ in terms of the power sum symmetric functions
\begin{gather}
Z\big(\rho^{+}_{1},\xi^{-,[1,u]}_1; \rho^{-}_{0}\big) = \exp\left(\sum_{l=1}^{\infty}\frac{p_l(\rho^{-}_{0})\big(p_l\big(\rho^{+}_{1}\big)+p_l\big(\xi^{-,[1,u]}_{1}\big)\big)}{l}\right.\nonumber\\
\left.\phantom{Z\big(\rho^{+}_{1},\xi^{-,[1,u]}_1; \rho^{-}_{0}\big) =}{} +\sum_{l=1}^\infty\frac{p^2_{l}\big(\rho^{+}_{1},\xi^{-,[1,u]}_1\big)-p_{2l}\big(\rho^{+}_{1},\xi^{-,[1,u]}_1\big)}{2l}\right),\label{eq:RecallSomeRep1}\\
F\big(\xi^{-,[1,u]}_{1};\xi^{+,[1,u]}_{1}\big) = \exp\left(\sum_{l=1}\frac{p_l\big(\xi^{-,[1,u]}_{1}\big)p_l\big(\xi^{+,[1,u]}_{1}\big)}{l}\right),\label{eq:RecallSomeRep2}\\
H^{(1)}_{ q_{m,j}}\big(\rho^{+}_{1},\xi^{-,[1,u]}_1, \rho^{-}_{0}; z_{1,j}\big) = \exp\left(\sum_{l=1}^\infty \frac{\big(p_l\big(\rho^{+}_1,\rho^{-}_0\big)+ p_l\big(\xi^{-,[1,u]}_{1}\big)\big)z_{1,j}^l\big[q^l_{1,j}-1\big]}{l}\right).\label{eq:RecallSomeRep3}
\end{gather}

We compute the innermost inner product of \eqref{eq:ULimitAlt} and let $u_{1,1}\to \infty$. On doing so, we get
\begin{gather}
\lim_{u_{1,1}\to \infty} \bigg\langle Z\big(\rho^{+}_{1},\xi^{-,[1,u_{1,1}]}_1; \rho^{-}_{0}\big)\prod_{j=1}^{d_1} H^{(1)}_{ q_{1,j}}\big(\rho^{+}_{1},\xi^{-,[1,u_{1,1}]}_1 \rho^{-}_{0}; z_{1,j}\big)H^{(2)}_{q_{1,j}}\big(\rho^{+}_{1},\xi^{-,[1,u_{1,1}]}_1; z_{1,j}\big), \nonumber \\
\qquad\quad{} F\big(\xi^{-,[1,u_{1,1}]}_{1};\xi^{+,[1,u_{1,2}]}_{1}\big)\bigg\rangle^{\xi^{-,[1,u_{1,1}]}_1}_{u_{1,1}}\label{eq:InnermostSimplified}\\
 \qquad{} = Z\big(\rho^{+}_1,\xi^{+,[1,u_{1,2}]}_{1}; \rho^{-}_{0}\big)\prod_{j=1}^{d_1} H^{(1)}_{ q_{1,j}}\big(\rho^{+}_{1},\xi^{+,[1,u_{1,2}]}_1 ,\rho^{-}_{0}; z_{1,j}\big)H^{(2)}_{q_{1,j}}\big(\rho^{+}_{1},\xi^{+,[1,u_{1,2}]}_1; z_{1,j}\big).\nonumber
\end{gather}
To see how one computes the limit in \eqref{eq:InnermostSimplified}, note that the two terms inside the inner product in~\eqref{eq:InnermostSimplified} can be written as $\phi\big(\xi^{-,[1,u_{1,1}]}_1\big)$ and $\psi\big(\xi^{-,[1,u_{1,1}]}_1\big)$ as in~\eqref{eq:SymPol}. So, the limiting value in~\eqref{eq:InnermostSimplified} now follows from Lemma~\ref{InnerProductLemma1}. Now, we turn to compute the second inner product of~\eqref{eq:ULimitAlt} and let $u_{1,2}\to \infty$.
To this end, let us define
\begin{gather*}
A :=Z\big(\rho^{+}_1,\xi^{+,[1,u_{1,2}]}_{1}; \rho^{-}_{0}\big)\prod_{j=1}^{d_1} H^{(1)}_{ q_{1,j}}\big(\rho^{+}_{1},\xi^{+,[1,u_{1,2}]}_1 ,\rho^{-}_{0}; z_{1,j}\big)H^{(2)}_{q_{1,j}}\big(\rho^{+}_{1},\xi^{+,[1,u_{1,2}]}_1; z_{1,j}\big),\\
B :=F\big(\rho^{+}_{2}; \rho^{-}_{1},\xi^{+,[1,u_{1,2}]}_{1}\big)\prod_{j=1}^{d_{2}}H^{(1)}_{q_{2,j}}\big(\rho^{-}_{1},\xi^{+,[1,u_{1,2}]}_{1};z_{2,j}\big)H^{(2)}_{q_{2,j}} \big(\rho^{+}_{2}; z_{2,j}\big).
\end{gather*}
Simplifying the second most inner product in \eqref{eq:ULimitAlt} boils down to taking the inner product of~$A$ and~$B$ and letting $u_{1,2}$ to converge to $\infty$.
Let us define
\begin{gather*} Y := \big(\rho^{+}_1, \rho^{-}_0, z_{1,j}, q_{1,j}z_{1,j}\big), \qquad W_1 :=(\rho^{+}_2, z_{2,j}, q_{2,j}z_{2,j}) , \qquad W_2 := (\rho^{-}_1).\end{gather*}
 Owing to the relations \eqref{eq:RecallSomeRep1}--\eqref{eq:RecallSomeRep3}, one can write
 \begin{gather*} A = \Psi_1(Y)\exp\left(\sum_{l=1}^\infty \frac{p_l\big(\xi^{+,[1,u_{1,2}]}_1\big)}{l}a_l(Y) + \sum_{l=1}^{\infty}\frac{p^2_l\big(\xi^{+,[1,u_{1,2}]}_1\big)- p_{2l}\big(\xi^{+,[1,u_{1,2}]}_1\big)}{2l}\right)\end{gather*} and
 \begin{gather*} B = \Psi_2(W_1, W_2)\exp\left(\sum_{l=1}^\infty \frac{p_l\big(\xi^{+,[1,u_{1,2}]}_1\big)}{l}b_l(W_1)\right)\end{gather*} for some $\{a_k\}_{k\in \NN}$, $\{b_k\}_{k\in \NN}$ and $\Psi_1$, $\Psi_2$ of a $\ZZ_{\geq 0}$-graded algebra. Lemma~\ref{InnerProductLemma1} implies that
\begin{gather}
\lim_{u_{1,2}\to \infty} \langle A, B\rangle^{\xi^{+,[1,u_{1,2}]}_1}_{u_{1,2}} = Z\big(\rho^{+}_1,\rho^{+}_{2},\xi^{-,[1,u_{2,1}]}_{2};\rho^{-}_{0}\big)F\big(\rho^{+}_{2},\xi^{-,[1,u_{2,1}]}_{2}; \rho^{-}_{1}\big) \nonumber\\
\qquad{}\times \prod_{j=1}^{d_1}H^{(1)}_{ q_{1,j}}\big(\rho^{+}_{1},\rho^{-}_0,\rho^{+}_{2}; z_{1,j}\big) \prod_{j=1}^{d_{2}}H^{(1)}_{q_{2,j}}\big(\rho^{+}_{1},\rho^{-}_0,\rho^{-}_1,\rho^{+}_{2}; z_{2,j}\big)\prod_{j=1}^{d_1} H^{(2)}_{ q_{1,j}}\big(\rho^{+}_1,\rho^{+}_2; z_{1,j}\big)\nonumber\\
\qquad{}\times
\prod_{j=1}^{d_{2}} H^{(2)}_{q_{2,j}}\big(\rho^{+}_2; z_{1,j}\big)\prod_{j=1}^{2d_2}\frac{1-q_{2,j}z^2_{2,j}}{1-z^2_{2,j}}
 \prod_{j=1}^{d_1}\prod_{k=1}^{d_{2}}H^{(1)}_{q_{2,k}}\big(z^{-1}_{1,j};z_{2,k}\big)\big(H^{(1)}_{q_{2,k}}\big((q_{1,j}z)^{-1}_{1,j};z_{2,k}\big)\big)^{-1}\nonumber\\
\qquad{} \times \prod_{j=1}^{d_1}\prod_{k=1}^{d_{2}}H^{(1)}_{q_{2,k}}(q_{1,j}z_{1,j};z_{2,k})\big(H^{(1)}_{q_{2,k}}(z_{1,j};z_{2,k})\big)^{-1}\nonumber\\
\qquad{} \times \exp\left(-\sum_{l=1}^\infty \frac{ (-p_{l}(\{z_{2,j}\}_{j} )+ p_{l}(\{q_{2,j}z_{2,j}\}_j) )^2-p_{2l} (\{z_{2,j}\}_{j}, \{q_{2,j}z_{2,j}\}_j )}{2l}\right).\label{eq:2ndInnerInnerProd}
\end{gather}
In \eqref{eq:2ndInnerInnerProd}, using the definitions of $H^{(1)}$ and $H^{(2)}$ from \eqref{eq:DefineH1}--\eqref{eq:DefineH2}, one can do the following simplifications
\begin{gather}
H^{(1)}_{q_{2,k}}(z^{-1}_{1,j};z_{2,k})\big(H^{(1)}_{q_{2,k}}\big((q_{1,j}z)^{-1}_{1,j};z_{2,k}\big)\big)^{-1} =\frac{(q_{1,j}z_{1,j}-q_{2,j}z_{2,k})(z_{1,j}-z_{2,k})}{(z_{1,j}-q_{2,k}z_{2,k})(q_{1,j}z_{1,j}-z_{2,k})},\label{eq:ComponentOfZ1}\\
H^{(1)}_{q_{2,k}}(q_{1,j}z_{1,j};z_{2,k})\big(H^{(1)}_{q_{2,k}}(z_{1,j};z_{2,k})\big)^{-1} =\frac{(1-q_{1,j}z_{1,j}z_{2,k})(1-q_{2,j}z_{2,j}z_{1,k})}{(1-q_{1,j}z_{1,j}q_{2,k}z_{2,k})(1-z_{1,j}z_{2,k})}.\label{eq:ComponentOfZ2}
\end{gather}

Moreover, owing to \eqref{eq:SpecializedH0H1}, the last term of the r.h.s.\ in \eqref{eq:2ndInnerInnerProd} can be further simplified to
\begin{gather}\label{eq:NeatRepresentationTrick}
\prod_{1\leq j<k\leq d_2}\frac{(1-q_{2,j}z_{2,j}z_{2,k})(1-q_{2,k}z_{2,k}z_{2,j})}{(1-q_{2,j}z_{2,j}q_{2,k}z_{2,k})(1-z_{2,j}z_{2,k})}.
\end{gather}
Substituting \eqref{eq:ComponentOfZ1}, \eqref{eq:ComponentOfZ2} and \eqref{eq:NeatRepresentationTrick} into the right hand side of \eqref{eq:2ndInnerInnerProd} yields \eqref{eq:ContourIntegralAlt}. This completes the proof.
\end{proof}

\subsection[Contour integral formula of the limit of $C(\rho^{+};\rho^{-};u)$]{Contour integral formula of the limit of $\boldsymbol{C(\rho^{+};\rho^{-};u)}$}\label{sec:CIFofC}

The main goal of this section is to complete the proof of Theorem~\ref{ContourIntegralRepOfC} by employing the computations presented in Proposition~\ref{eq:NestedInnerProductSummation} and Theorem~\ref{LimitRep}.

\begin{proof}[Proof of Theorem~\ref{ContourIntegralRepOfC}]
We start with simplifying the innermost inner product in \eqref{eq:ProofOfRep}. Substituting the contour integral formula of Theorem~\ref{OneMainTheorem} inside the inner product yields
\begin{gather*}
 \left\langle \prod_{j=1}^{d_1} D^{1, q_{1,j}}_{(n_1+u_{1,1}),[\rho^{+}_{1},\xi^{-,[1,u_{1,1}]}_1]} Z\big(\rho^{+}_{1},\xi^{-,[1,u_{1,1}]}_1; \rho^{-}_{0}\big),F\big(\xi^{-,[1,u_{1,1}]}_{1};\xi^{+,[1,u_{1,2}]}_{1}\big)\right\rangle^{\xi^{-,[1,u_{1,1}]}_1}_{u_{1,1}}\\
\qquad{} = \left\langle \oint_{\mathcal{C}_{1,d_1}}\cdots \oint_{\mathcal{C}_{1,1}} \prod_{j=1}^{d_1}\frac{\prod\limits_{j=1}^{d_1} q^{n_1+u_{1,1}}_{1,j}}{q_{1,j}z_{1,j}-z_{1,j}}\prod_{1\leq j<k\leq d_1}\frac{(q_{1,j}z_{1,j}-q_{1,k}z_{1,k})(z_{1,j}-z_{1,k})}{(z_{1,j}-q_{1,k}z_{1,k})(q_{1,j}z_{1,j}-z_{1,k})} \right.\\
\qquad\quad{} \times \prod_{j=1}^{d_1}\frac{1-q_{1,j}z^2_{1,j}}{1-z^2_{1,j}}\prod_{1\leq j<k\leq d_1}\frac{(1-q_{1,j}z_{1,j}z_{1,k})(1-q_{1,k}z_{1,k}z_{1,j})}{(1-q_{1,j}z_{1,j}q_{1,k}z_{1,k})(1-z_{1,j}z_{1,k})}Z\big(\rho^{+}_{1},\xi^{-,[1,u_{1,1}]}_1; \rho^{-}_{0}\big)\\
\qquad\quad{} \times \prod_{j=1}^{d_1} H^{(1)}_{ q_{1,j}}\big(\rho^{+}_{1},\xi^{-,[1,u_{1,1}]}_1 ,\rho^{-}_{0}; z_{1,j}\big)H^{(2)}_{q_{1,j}}\big(\rho^{+}_{1},\xi^{-,[1,u_{1,1}]}_1; z_{1,j}\big)\\
\left.\qquad\quad{}\times \prod_{j=1}^{d_1}dz_{1,j}, F\big(\xi^{-,[1,u_{1,1}]}_{1};\xi^{+,[1,u_{1,2}]}_{1}\big) \right\rangle^{\xi^{-,[1,u_{1,1}]}_1}_{u_{1,1}},
\end{gather*}
where the contours of the variables $\{z_{1,1},\dots ,z_{1,d_1}\}$ contains the poles at $\rho^{+}_{1}\cup \big\{\xi^{-,[1,u_{1,1}]}\big\}$. To use Theorem~\ref{LimitRep}, we need $\mathbf{(i)}$ to bring the inner product inside the integral, $(ii)$ then to interchange the limit and the integral and~$(iii)$ lastly, we have to look at how the behavior of the poles changes if we somehow manage to do steps~$(i)$ and~$(ii)$. In what follows, we address those issues one at a time.

Before going into the details of the justifications of $(i)$, $(ii)$ and $(iii)$, we perform a deformation of the contours for the variables $\{z_{1,j}\}_{j\in [1,d_1]}$. Recall that the contours $\{\mathcal{C}_{1,j}\}_j$ are composed of circles of very small radii around the points $\big\{\rho^{+}_1,\xi^{-,[1,u_{1,1}]}_{1}\big\}$. We claim that one can deform each of them to a contour which is composed of the negatively oriented circle $|z|=r$ and positively oriented circle $|z|=R$ where $r$ and $R$ must satisfy the conditions
\begin{gather}\label{eq:ConstrainOnRadii}
r<\min |\rho^{+}_1|\wedge \min \big|\xi^{-,[1,u_{1,1}]}_1\big| \leq \max |\rho^{+}_1|\vee \max \big|\xi^{-,[1,u_{1,1}]}_1\big| <R<1,
\end{gather}
where $|\rho^{+}_1|$ (or $\big|\xi^{-,[1,u_{1,1}]}_1\big|$) denotes the set of the absolute values of all the variables in $\rho^{+}_1$ (or $\big|\xi^{-,[1,u_{1,1}]}_1\big|$). This claim is in the same spirit of \cite[Proposition~2.2.7]{A14}. We briefly prove this claim here. While deforming the contours, we have to take into account the residues at the poles of types $(a)$~$z_{1,j}=q_{1,k}z_{1,k}$, $(b)$~$z_{1,j}=z^{-1}_{1,k}$ and $(c)$ $z_{1,j}=(q_{1,j}\zeta)^{-1}$ where $\zeta \in \xi^{-,[1,u_{1,1}]}_1\cup \rho^{+}_1\cup\rho^{-}_1$. Note that the residues of the poles of type $(a)$ are always zero due to the presence of $1-(q_{1,j}z_{1,j})^{-1}\theta$ for $\theta \in \xi^{-,[1,u_{1,1}]}_1\cup \rho^{+}_{1}$ in the numerator of the factor $H^{(2)}_{q_{1,j}}\big(\rho^{+}_{1},\xi^{-,[1,u_{1,1}]}_1; z_{1,j}\big)$. Similarly, due to the presence of the factor $1-z_{1,j}\theta$ for $\theta\in \xi^{-,[1,u_{1,1}]}_1\cup \rho^{+}_{1}$ in the numerator of $H^{(1)}_{ q_{1,j}}\big(\rho^{+}_{1},\xi^{-,[1,u_{1,1}]}_1 ,\rho^{-}_{0}; z_{1,j}\big)$, the residues at the poles of type $(b)$ are always zero. Moreover, the poles of type~$(c)$ contribute nothing while deforming the contours because of the constraints on~$\{q_{1,j}\}_j$ in~\eqref{eq:ConstrainONQ}. Now, we turn to justifying the heuristics $(i)$, $(ii)$ and $(iii)$ in details.

$(i)$ Since the integrand is uniformly bounded over the circles $|z|=r$ and $|z|=R$, we can interchange the inner product and the integral using bilinearity of the inner product.

$(ii)$ Now, we justify the interchange of the limit and the integral. First, we fix $u_{1,1}$. We denote \[Y:={\rm Sym}\big(\rho^{-}_{0}, \rho^{+}_1, \{q_{1,j}z_{1,j}\}_j, \big\{(q_{1,j}z_{1,j})^{-1}\big\}_j, \{z_{1,j}\}_j, \big\{z^{-1}_{1,j}\big\}_j\big)\] and then define a homomorphism
\begin{gather*}
\Psi \colon \ {\rm Sym}(Y)\mapsto {\rm Sym}(Y),\\
\Psi(p_k(Y)) = p_k\big(\rho^{-}_{0}, \rho^{+}_1,\{q_{1,j}z_{1,j}\}_j,\big\{z^{-1}_{1,j}\big\}_j\big)- p_k\big(\big\{{-}(q_{1,j}z_{1,j})^{-1}\big\}_j, \{-z_{1,j}\}_j\big),
\end{gather*}
where $\{p_k\}_k$ is the set of power sum symmetric polynomials.
Owing to the representations in \eqref{eq:RecallSomeRep1}--\eqref{eq:RecallSomeRep3}, one can write the inner product after taking it inside the integrand in the following way
\begin{gather}
\left\langle \sum_{\lambda^{(1)}\in \mathbb{Y}} s_{\lambda^{(1)}}\big(\xi^{-,[1,u_{1,1}]}_1\big) \Psi\big(\tau_{\lambda^{(1)}}(Y)\big), \sum_{\lambda^{(1)}\in \mathbb{Y}} s_{\lambda^{(1)}}\big(\xi^{-,[1,u_{1,1}]}_1\big) s_{\lambda^{(1)}}\big(\xi^{+,[1,u_{1,2}]}_{1}\big)\right\rangle^{\xi^{-,[1,u_{1,1}]}_1}_{u_{1,1}}\nonumber\\
\qquad{} =\sum_{\lambda^{(1)}\in \mathbb{Y}} s_{\lambda^{(1)}}\big(\xi^{+,[1,u_{1,2}]}_1\big) \Psi\big(\tau_{\lambda^{(1)}}(Y)\big)\mathbbm{1}\big(\big|\lambda^{(1)}\big|\leq u_{1,1}\big).\label{eq:InteriorInnerProductAtWork}
\end{gather}
 Now, we seek to bound $\psi(\tau_{\lambda^{(1)}})$ using the fact $\Psi$ is a algebra homomorphism. For the sake of notational convenience, we denote $\lambda^{(1)}$ just by $\lambda$ for the rest of this discussion. Recall that we have $s_\lambda=\det(h_{\lambda_i-i+j})$ and henceforth, $\Psi(s_\lambda)= \det(\Psi(h_{\lambda_i-i+j}))$. Thus, using the identity $|\det(A)|=\sqrt{\det(AA^*)}$ and the inequality $\det(D)\leq \prod_{i}D_{ii}$ for any positive semi-definite mat\-rix~$D$, one can write $|\Psi(s_\lambda)|\leq \prod_{i}\big(\sum_{j}\Psi(h_{\lambda_i-i+j})^2\big)^{1/2}$. Furthermore, $h_r$ is the coefficient of $t^r$ in the expansion of $\prod_{k\geq 1}\big(1+t^kp_k\big)$. Therefore, $\Psi(h_r)$ is given by the coefficient of~$t^r$ in $\prod_{k\geq 1}\big(1+t^k\Psi(p_k)\big)$. Moreover, one can bound $|\Psi(p_r(Y))|$ by $p_r(|Y|)$ where~$|Y|$ denotes the set of the absolute values of the variables in $Y$. This shows $|\Psi(h_r(Y))|$ is also bounded above by~$h_r(|Y|)$. Consequently, one can bound~$|\Psi(s_\lambda)|$ by $\prod_{i}\big(\sum_{j}h_{\lambda_i-i+j}^2\big)^{1/2}$. But, $h_{\lambda_i-i+j}(|Y|)$ is bounded above $(1+|\lambda|)^{l(\lambda)}(\max |Y|)^{\lambda_i-i+j}$. This implies one can bound $|\Psi(s_\lambda)|$ by $(1+|\lambda|)^{l^2(\lambda)}(\max |Y|)^{|\lambda|}$. Similarly, one can show that $|\Psi(s_{\lambda/\mu})|$ can be bounded by $(1+|\lambda|)^{l^2(\lambda)}(\max |Y|)^{|\lambda|-|\mu|}$. To that effect, we get $|\Psi(\tau_\lambda(Y))|\leq (1+|\lambda|)^{l^2(\lambda)+1}(\max |Y|)^{|\lambda|}$. Now note that when $l(\lambda)$ is large in~\eqref{eq:InteriorInnerProductAtWork}, $s_\lambda\big(\xi^{+,[1,u_{1,2}]}_1\big)$ is zero. Henceforth, one can bound the right side of~\eqref{eq:InteriorInnerProductAtWork} by $\sum_{\lambda\in \mathbb{Y}} \big|s_\lambda\big(\xi^{+,u_{1,2}}_1\big)\big|(1+|u_{1,1}|)^{u_{1,2}}(\max |Y|)^{|\lambda|}$. Furthermore, $\max |Y|$ is bounded above by $\max {|z_{1,j}|}$. Combining all these observations with the fact that $\big|\max \xi^{+,[1,u_{1,2}]}_1\big|\times \max |Y|$ is less than $1$, we get
\begin{gather}
 \left|\sum_{\lambda\in \mathbb{Y}} s_\lambda\big(\xi^{+,u_{1,2}}_1\big)\Psi(\tau_\lambda(Y))\mathbbm{1}(|\lambda|\leq u_{1,1}) \right|\nonumber\\
 \qquad{} \leq \sum_{\substack{\lambda\in \mathbb{Y}}} (1+|\lambda|)^{2u_{1,2}} \big(\max |\xi^{+,[1,u_{1,2}]}_1|^{|\lambda|}\big) \big(\max |Y|^{|\lambda|}\big)<\infty .\label{eq:InnerProductWithMultiplierBound}
\end{gather}
Note that the bound on the right side in \eqref{eq:InnerProductWithMultiplierBound} does not depend on the value of $u_{1,1}$. This show after taking the inner product in~\eqref{eq:InteriorInnerProductAtWork} the integrand is uniformly bounded over all values of~$u_{1,1}$. Therefore, we can interchange the limit $u_{1,1}\to \infty$ and the integral.

$(iii)$ Here, we discuss the consequences of taking the limit of the following term
\begin{gather*}
\prod_{j=1}^{d_1} q^{-u_{1,1}}_{1,j} \!\left\langle \prod_{j=1}^{d_1} D^{1, q_{1,j}}_{(n_1+u_{1,1}),[\rho^{+}_{1},\xi^{-,[1,u_{1,1}]}_1]} Z\big(\rho^{+}_{1},\xi^{-,[1,u_{1,1}]}_1; \rho^{-}_{0}\big),F\big(\xi^{-,[1,u_{1,1}]}_{1};\xi^{+,[1,u_{1,2}]}_{1}\big)\right\rangle^{\xi^{-,[1,u_{1,1}]}_1}_{u_{1,1}}\!\!\!.
\end{gather*}
 as $u_{1,1}$ tends to $\infty$. Using~\eqref{eq:InnermostSimplified}, one can write the limiting value as
\begin{gather}
\oint_{\mathcal{C}_{1,d_1}}\cdots \oint_{\mathcal{C}_{1,1}} \prod_{j=1}^{d_1}\frac{\prod\limits_{j=1}^{d_1} q^{n_1}_{1,j}}{q_{1,j}z_{1,j}-z_{1,j}}\prod_{1\leq j<k\leq d_1}\frac{(q_{1,j}z_{1,j}-q_{1,k}z_{1,k})(z_{1,j}-z_{1,k})}{(z_{1,j}-q_{1,k}z_{1,k})(q_{1,j}z_{1,j}-z_{1,k})}\prod_{j=1}^{d_1}\frac{1-q_{1,j}z^2_{1,j}}{1-z^2_{1,j}}\nonumber\\
\qquad{} \times \prod_{1\leq j<k\leq d_1}\frac{(1-q_{1,j}z_{1,j}z_{1,k})(1-q_{1,k}z_{1,k}z_{1,j})}{(1-q_{1,j}z_{1,j}q_{1,k}z_{1,k})(1-z_{1,j}z_{1,k})}Z\big(\rho^{+}_1,\xi^{+,[1,u_{1,2}]}_{1}; \rho^{-}_{0}\big)\nonumber\\
\qquad {} \times \prod_{j=1}^{d_1} H^{(1)}_{ q_{1,j}}\big(\rho^{+}_{1},\xi^{+,[1,u_{1,2}]}_1 ,\rho^{-}_{0}; z_{1,j}\big) H^{(2)}_{q_{1,j}}\big(\rho^{+}_{1},\xi^{+,[1,u_{1,2}]}_1; z_{1,j}\big)\prod_{j=1}^{d_1}{\rm d}z_{1,j}, \label{eq:AfterTakingFirstLimit}
\end{gather}
where the contours $\{\mathcal{C}_{1,d_1}\}$ are now composed of two circles $|z|=r$ and $|z|=R$ satis\-fying~\eqref{eq:ConstrainOnRadii}. Thus, the variables $\{z_{1,j}\}$ acquire some new poles at $\xi^{+,[1,u_{1,2}]}_1$ coming out from the factor $H^{(2)}_{q_{1,j}}\big(\rho^{+}_{1},\xi^{+,[1,u_{1,2}]}_1; z_{1,j}\big)$ in~\eqref{eq:AfterTakingFirstLimit}.

Following the description of $C(\rho^{+},\rho^{-})$ in Proposition~\ref{eq:NestedInnerProductSummation}, we now consider taking the inner product of the integral in \eqref{eq:AfterTakingFirstLimit} with the following expression \[\prod_{j=1}^{d_{2}}D^{1,q_{2,j}}_{(n_2+u_{2,1}),[\rho^{+}_{2},\xi^{-,[1,u_{2,1}]}_{2}]} F\big(\rho^{+}_{2},\xi^{-,[1,u_{2,1}]}_{2}; \rho^{-}_{1},\xi^{+,[1,u_{1,2}]}_{1}\big).\] Using the steps $(i)$, $(ii)$ and $(iii)$ as described above, one can again bring the limit $u_{1,2}\to \infty$ inside the integral and obtain the limit using \eqref{eq:2ndInnerInnerProd}. Note that those steps also have to be accompanied by suitable deformation of the contours. Consequently, it changes the set of valid poles for the variables $z_{1,j}$ to $\rho^{+}_1\cup \rho^{+}_2\cup \xi^{-,[1,u_{2,1}]}_2\cup \{q_{2,j}z_{2,j}\}_{j\leq d_2}$. Now, using induction over $i$ (in conjunction with Theorem~\ref{LimitRep}), we see that the deformed contours of the variables $\{z_{i,j}\}_{j}$ now contain the poles at \[\rho^{+}_{[i,m]}\cup \{q_{i+1,j}z_{i+1,j}\}_{j\in [1,d_{i+1}]}\cup \cdots \cup \{q_{m,j}z_{m,j}\}_{j\in [1,d_{i+1}]}\]
for all $i\in [1,m-1]$. Furthermore, applying Theorem~\ref{LimitRep} (in a similar way as in \eqref{eq:AfterTakingFirstLimit}) at every step of the induction, we observe that the limiting value of $C(\rho^+;\rho^-;u)$ (after taking $u\to \infty$ as in \eqref{eq:SuccessiveLimitingNotations}) matches with the limiting expression of \eqref{eq:ContourIntegral12}. This completes the proof.
\end{proof}

\subsection*{Acknowledgements}
The author is grateful to his adviser Professor Ivan Corwin for suggesting the problem considered in this paper. We also thankfully acknowledge his numerous critical comments on the earlier versions of this paper. We are thankful to Guillaume Barraquand for his various remarks on the proofs and especially, suggesting a substitution in Theorem~\ref{TheoremForSinglePartition}. We thank the anonymous referees whose comments have greatly improved this manuscript. We also like to thank Professor Alexei Borodin for pointing out some mistakes in the citations which have been corrected.

\pdfbookmark[1]{References}{ref}
\LastPageEnding

\end{document}